\documentclass[12pt,leqno,twoside]{amsart}
\usepackage{amsmath,amscd,amssymb,amsfonts,latexsym,mathrsfs,bm,extarrows,graphicx,stackrel}
\usepackage[all,cmtip]{xy}
\usepackage{ulem}
\usepackage{tikz}

\usepackage{srcltx}
\usepackage[colorlinks=true, pdfstartview=FitV, linkcolor=blue, citecolor=blue, urlcolor=blue]{hyperref}
\usepackage{subcaption} 

\usepackage{amsbsy}
\usepackage{amsthm,alltt,dsfont}

\usepackage[margin=1.25in]{geometry}

\usepackage{mathptmx}

\newcommand{\Z}{\mathbb{Z}}
\newcommand{\C}{\mathbb{C}}

\newcommand{\Q}{\mathbb{Q}}

\allowdisplaybreaks

\newcommand{\bC}{{\mathbb C}}
\newcommand{\bD}{{\mathbb D}}

\newcommand{\bN}{{\mathbb N}}
\newcommand{\bP}{{\mathbb P}}
\newcommand{\bR}{{\mathbb R}}
\newcommand{\bQ}{{\mathbb Q}}
\newcommand{\bZ}{{\mathbb Z}}
\newcommand{\bW}{{\mathbb W}}

\newcommand{\bx}{{\bb x}}

\newcommand{\cD}{{\mathscr D}}

\newcommand{\cF}{{\mathscr F}}
\newcommand{\cG}{{\mathscr G}}
\newcommand{\cH}{{\mathscr H}}

\newcommand{\cM}{{\mathscr M}}

\newcommand{\cO}{{\mathscr O}}

\newcommand{\cS}{{\mathscr S}}
\newcommand{\cT}{{\mathscr T}}
\newcommand{\cL}{{\mathscr L}}

\newcommand\sV{{\mathcal V}}
\newcommand\sM{{\mathcal M}}
\newcommand{\cAb}{{\mathscr A}^{\bullet}}
\newcommand{\cBb}{{\mathscr B}^{\bullet}}

\newcommand{\cEb}{{\mathscr E}^{\bullet}}
\newcommand{\cFb}{{\mathscr F}^{\bullet}}
\newcommand{\cGb}{{\mathscr G}^{\bullet}}

\def\Z{\mathbb{Z}}
\def\C{\mathbb{C}}
\def\Q{\mathbb{Q}}

\newcommand{\im}{\mathrm{Im~}}
\renewcommand{\ker}{\mathrm{Ker}~}
\newcommand{\coker}{\mathrm{Coker}~}

\newcommand{\Tor}{{\rm{Tor}}}

\newcommand{\supp}{\mathrm{supp}}
\newcommand{\Sing}{\mathrm{Sing}} 
\newcommand{\wti}{\widetilde}

\newcommand{\Hom}{{\rm{Hom}}}

\newcommand{\codim}{\hbox{\rm codim}\,}

\newcommand{\lra}{\longrightarrow}

\theoremstyle{plain}
\newtheorem{thm}{Theorem}[section]
\newtheorem{cor}[thm]{Corollary}
\newtheorem{conj}[thm]{Conjecture}
\newtheorem{lem}[thm]{Lemma}
\newtheorem{prop}[thm]{Proposition}

\newtheorem{problem}[thm]{Problem}

\theoremstyle{definition}
\newtheorem{df}[thm]{Definition}
\newtheorem{rem}[thm]{Remark}
\newtheorem{example}[thm]{Example}
\newtheorem{exercise}[thm]{Exercise}

\def\be{\begin{equation}}
\def\ee{\end{equation}}

\def\bt{\begin{thm}}
\def\et{\end{thm}}

\def\bc{\begin{cor}}
\def\ec{\end{cor}}

\def\br{\begin{rem}}
\def\er{\end{rem}}

\def\bp{\begin{prop}}
\def\ep{\end{prop}}

\def\bpr{\begin{problem}}
\def\epr{\end{problem}}

\def\bl{\begin{lem}}
\def\el{\end{lem}}

\def\bn{\begin{enumerate}}
\def\en{\end{enumerate}}

\def\bex{\begin{example}}
\def\eex{\end{example}}

\def\bd{\begin{df}}
\def\ed{\end{df}}

\def\bx{\begin{exercise}}
\def\ex{\end{exercise}}

\def\ben{\begin{enumerate}}
\def\een{\end{enumerate}}

\DeclareMathOperator{\Char}{Char}
\DeclareMathOperator{\Alb}{Alb}
\DeclareMathOperator{\alb}{alb}
\DeclareMathOperator{\spec}{Spec}

\DeclareMathOperator{\rhomo}{RHom}
\renewcommand{\r}{{\mathrm{rHd}}}

\newcommand{\cosupp}{\mathrm{cosupp}}

\newcommand{\leftrarrows}{\mathrel{\raise.75ex\hbox{\oalign{%
  $\scriptstyle\leftarrow$\cr
  \vrule width0pt height.5ex$\hfil\scriptstyle\relbar$\cr}}}}
\newcommand{\lrightarrows}{\mathrel{\raise.75ex\hbox{\oalign{%
  $\scriptstyle\relbar$\hfil\cr
  $\scriptstyle\vrule width0pt height.5ex\smash\rightarrow$\cr}}}}
\newcommand{\Rrelbar}{\mathrel{\raise.75ex\hbox{\oalign{%
  $\scriptstyle\relbar$\cr
  \vrule width0pt height.5ex$\scriptstyle\relbar$}}}}

\makeatletter
\def\leftrightarrowsfill@{\arrowfill@\leftrarrows\Rrelbar\lrightarrows}
\newcommand{\xleftrightarrows}[2][]{\ext@arrow 3399\leftrightarrowsfill@{#1}{#2}}
\makeatother

\title[Constructible sheaf complexes]{Constructible sheaf complexes in complex geometry \\ and Applications}

\author[L. Maxim]{Lauren\c{t}iu G. Maxim}
\address{L. Maxim: Department of Mathematics, University of Wisconsin-Madison, USA.}
\email {maxim@math.wisc.edu}
\author[J. Sch\"urmann ]{J\"org Sch\"urmann}
\address{J.  Sch\"urmann : Mathematische Institut, Universit\"at M\"unster, Einsteinstr. 62, 48149 M\"unster, Germany.}
\email {jschuerm@math.uni-muenster.de}

\date{\today}



\thanks{L. Maxim is partially supported by the Simons Foundation Collaboration Grant \#567077. J. Sch\"{u}rmann is funded by the Deutsche Forschungsgemeinschaft (DFG, German Research Foundation) Project-ID 427320536 -- SFB 1442, as well as under Germany’s Excellence Strategy EXC 2044 390685587, Mathematics M\"{u}nster: Dynamics -- Geometry -- Structure.}

\begin{document}

\maketitle

\begin{abstract}
We present a detailed introduction of the theory of  constructible sheaf complexes in the complex algebraic and analytic setting.  All concepts are illustrated by many interesting examples and relevant applications, while some important results are  presented with complete proofs. This paper is intended as a broadly accessible user's guide to these topics, providing the readers with a taste of the subject, reflected by concrete examples and applications that motivate the general theory. We discuss the stability of  constructible sheaf complexes under the standard functors, and explain the relation of these functors to perverse sheaves and the perverse $t$-structure. We introduce the main results of stratified Morse theory in the framework of  constructible sheaves, for proving the basic vanishing and finiteness results. Applications are given to various index theorems,  the functorial calculus of characteristic cycles of constructible functions, and  to  weak Lefschetz and Artin-Grothendieck type theorems. We recall the construction of Deligne's nearby and vanishing cycle functors, prove  that they preserve  constructible complexes, and discuss  their relation with the  perverse $t$-structure. We finish this paper with a description and applications of  the K\"ahler package for intersection cohomology of complex algebraic varieties, and the recent study of perverse sheaves on semi-abelian varieties.
\end{abstract}

\tableofcontents


\section{Introduction}
\label{sec:1}

The main goal of this paper is to provide a user's guide, both for the novice and the expert, for the theory of (weakly) constructible sheaf complexes in complex geometry and their many applications.  Our guiding principle for writing these notes was to provide an explicit and geometric introduction of the mathematical concepts, while also explaining some of the most important examples of the general theory. For this reason, we aim to present as many of the basic results as possible, sometimes even with complete proofs in some special important cases. Moreover, these results and definitions are then always explained and illustrated by many examples.

{\it Constructible (complexes of) sheaves} are the algebraic counterpart of the decomposition of a variety into  manifolds pieces (strata), and they are, roughly speaking, obtained by gluing local systems defined along strata of a Whitney stratification $\cS$. {\it Perverse sheaves} are an important class of constructible complexes, introduced in \cite{BBD} as a formalization of the celebrated Riemann-Hilbert correspondence of Kashiwara \cite{Ka}, which relates the topology of algebraic, resp., analytic  varieties (intersection homology) and the algebraic, resp., analytic theory of differential equations (holonomic $D$-modules).

In recent years, constructible sheaf complexes and especially perverse sheaves have become indispensable tools for studying complex algebraic and analytic varieties.
They have seen spectacular applications in geometry and topology (e.g., the decomposition theorem \cite{BBD} and the topology of complex algebraic maps), but also in fields such as representation theory (e.g., proof of the Kazhdan-Lusztig conjecture, proof of the geometrization of the Satake isomorphism, and proof of the fundamental lemma in the Langlands program) or combinatorics (e.g., Stanley's proof of McMullen's conjecture, or the resolution of the Dowling-Wilson top-heavy conjecture);  see, e.g.,  \cite{CM2, M} for more recent surveys of such applications. 
Furthermore, perverse sheaves and the {\it nearby and vanishing cycle functors} of Deligne \cite{Gr2} are the backbone of Saito's mixed Hodge module theory \cite{Sa1, Sa2}, a far-reaching generalization of Deligne's mixed Hodge theory. 

However, despite their fundamental importance, perverse sheaves as special constructible complexes of sheaves remain rather mysterious objects. It is our hope that the present paper will help readers become better acquainted with various aspects of the general theory.
Those looking to delve further into more specialized topics or wishing to explore problems of current research will find ample references to facilitate navigation of both classic and recent literature.

\medskip

Let us next give a brief summary of the content of the paper.

\medskip

In Section \ref{sec:co}, we define the notion of (weak) constructibility, and discuss the stability of (weakly) constructible sheaf complexes under the standard functors. As far as possible, we allow also weakly constructible sheaf complexes, where one does not impose any finiteness conditions for the stalks
(and which is sometimes even more natural or simpler to work with, especially as long as no duality is used). Similarly, we try to work in such a way that it applies to the complex algebraic as well as complex analytic context (sometimes only under suitable compactness assumptions in the complex analytic context). The presented calculus includes {\it external products, K\"{u}nneth isomorphisms, Verdier duality} and the relation to the {\it Euler characteristic calculus of constructible functios}.

We also introduce here the perverse t-structure and perverse sheaves (with respect to middle perversity), and explain their relation with the standard functors
(following  their counterparts in $l$-adic cohomology as presented in \cite[Chapter 4]{BBD}). Several aspects of the general theory are worked out in detail for intersection cohomology complexes, which provide some of the main examples of perverse sheaves. Furthermore, for constructible sheaf complexes of 
$R$-modules with $R$ a {\it Dedekind domain}, we also consider the {\it dual perverse t-strucure} and its relation to the {\it rectified homological depth} of Grothendieck, as studied by Hamm and L\^e \cite{HL} (also for the corresponding homotopical notion).

\medskip

In Section \ref{smt}, we explain  the basic results from \cite{Sc} about {\it stratified Morse theory} in the framework of (weakly)  constructible sheaves in the {\it complex context}, continuing  and extending the recent survey of Goresky \cite {Gor} in this handbook series, as well as  Massey's survey  \cite{Ma1}. We follow here the notions of the geometric {stratified Morse theory} of Goresky-Mac\-Pherson \cite{GMs}, so that one can easily compare the results of our paper with those of {\it loc.cit.}. We introduce, for example, the sheaf theoretic counterparts of the {\it local and normal Morse data}, 
as well as their relations for a {\it stratified Morse critical point} of a $C^{\infty}$-function on a complex algebraic (or analytic) variety.
The normal Morse data of (weakly) constructible sheaf complexes are studied via the
{\it complex link} of a stratum of a Whitney stratification, which allows to prove the basic {\it vanishing and finiteness theorems} by induction on the dimension of the underlying complex analytic variety. In particular, we  get a description of the (dual) perverse t-structure in terms of properties of the normal Morse data.
 Moreover, we also explain some relations to the general {\it micro-local sheaf theory} of Kashiwara-Schapira \cite{KS}, e.g. like a description of the {\it micro-support} of a (weakly) constructible sheaf complex in terms of the normal Morse data.

We use the language of {\it stratified Morse theory for constructible functions and sheaves}  to also give in this section an  introduction to the functorial theory of {\it Lagrangian cycles} in the complex analytic and algebraic context. We discuss from this viewpoint the {\it Euler isomorphism} between  constructible functions and Lagrangian cycles given by the {\it characteristic cycle} of a constructible function, together with some index theorems. 
Examples are given to {\it Poincar\'{e}-Hopf index theorems for singular spaces}, {\it effective characteristic cycles on abelian varieties} \cite{AMSS2} and the
{\it global Euler obstruction for affine varieties} \cite {STV}, as well as the famous {\it local Euler obstruction} of MacPherson \cite{Mch}.
We also  explain (using  this language of {\it stratified Morse theory for constructible functions}) the translation into the context of Lagrangian cycles of the following operations for constructible functions and sheaves: external product, proper direct image, non-characteristic pullback and specialization (i.e., nearby cycles),
together with an intersection formula for vanishing cycles. 

Finally, the last part of Section \ref{smt} deals with applications of the {\it stratified Morse theory for constructible sheaves} to {\it vanishing and weak Lefschetz theorems} in the complex algebraic and analytic context. This includes different versions of the {\it Artin vanishing theorem} for complex algebraically
(weakly) constructible complexes on an {\it affine} variety, and {\it vanishing theorems} for (weakly) constructible complexes on complex analytic 
{\it Stein and $q$-complete varieties}, as well as relative counterparts for morphisms given by {\it Artin-Grothendieck types theorems} in the complex algebraic and analytic context.

\medskip

In Section \ref{sec:nvc}, we recall the construction of Deligne's {\it nearby and vanishing cycle functors} \cite{Gr2}, and prove  that they preserve (weakly) constructible complexes. For the constructible context we also need and explain  the calculation of their {\it (co)stalks} in terms of the {\it Milnor fibers} for the corresponding local Milnor fibrations,
\cite{Le2} based on the existence of an adopted Whitney stratifications satisfying the {\it $a_f$-condition of Thom} \cite{BMM, Hir, LeTe}.
Then we state the relation of  {\it nearby and vanishing cycle functors} with duality  \cite{Ma2}, and prove their relation with the (dual) perverse $t$-structure and perverse sheaves (from the point of view of stratified Morse theory, as developed in the previous section), i.e.,
$$\psi_f[-1]\quad \text{and} \quad  \varphi_f[-1] \quad \text {are {\it t-exact} for the (dual) perverse t-structure.}$$ 
We also include here a discussion on the {\it Thom-Sebastiani theorem} for vanishing cycles \cite{Ma3, Sc}, and give a brief description of Beilinson's
and Deligne-Verdier's procedure for {\it gluing perverse sheaves} via vanishing cycles \cite{Bei, Ve2}.
As a final application we explain the equality $Rf_!=Rf_*$ on the level of Grothendieck groups of algebraically constructible complexes \cite{La, Vir},
as well as analytically constructible functions in a compactifiable complex analytic context \cite{Sc}.

\medskip

In Section \ref{sec:ih}, we overview properties of the {\it intersection cohomology} groups of complex algebraic varieties, which generalize the corresponding features of the cohomology groups of smooth varieties. These properties, consisting of {\it Poincar\'e duality, (weak and hard) Lefschetz theorems} and the {\it decomposition theorem}, are collectively termed the {\it K\"ahler package} for intersection cohomology \cite{BBD, CM, CM2, Sa1, Sa2}.

We also mention briefly a recent combinatorial application of the K\"ahler package for intersection cohomology, namely a proof (for realizable matroids) by Huh-Wang \cite{HW} of the {\it Dowling-Wilson top-heavy conjecture} \cite{DW1, DW2}.

\medskip

Finally, in Section \ref{semi-ab}, we survey recent developments in the study of {\it perverse sheaves on semi-abelian varieties}. We also include several concrete applications of this theory, e.g., to the study of homotopy types of complex algebraic manifolds (formulated in terms of their cohomology jump loci), as well as new topological characterizations of semi-abelian varieties \cite{LMW3, LMW4, LMW5}.

\medskip

We assume reader's familiarity with derived categories and the derived calculus; for a quick refresher on these topics the interested reader may consult  \cite{Di},  \cite{M} or \cite[Chapter I-III]{KS}.

We work in the complex algebraic or analytic setting with reduced  Hausdorff spaces, e.g., in the complex algebraic context we are working with the complex analytic space associated to a reduced separated scheme of finite type over $\spec(\bC)$. Unless otherwise specified, all dimensions are taken to be complex dimensions. 
In the complex analytic setting we always assume that our spaces have bounded dimension and a countable topology, e.g., the disjoint union of all 
$\bC^n$'s ($n\in \bN_0$) is not allowed.


\section{Constructible and perverse sheaf complexes}\label{sec:co}

\subsection{Constructibility}

Let $X$ be a complex algebraic (or analytic) variety (here a variety does not need to be irreducible). It is well known that such a variety can be endowed with a \index{Whitney stratification} {\it Whitney stratification}, i.e., a (locally) finite partition $\mathscr{S}$ into non-empty, connected, locally closed nonsingular subvarieties $S$ of $X$ (called ``strata'') which satisfy the following properties: 
\begin{enumerate}
\item[(i)] {\it Frontier condition}: for any stratum $S \in \mathscr{S}$, the frontier $\partial S:={\bar S} \setminus S$ is a union of strata of $\mathscr{S}$, where ${\bar S}$ denotes the closure of $S$.
\item[(ii)] {\it Constructibility}: the closure ${\bar S}$ and the frontier $\partial S$ of any stratum $S \in \mathscr{S}$ are closed complex algebraic (respectively, analytic) subspaces in $X$.
\end{enumerate}
These conditions already imply that $X$ gets an induced dimension filtration
$$X_{\bullet}: \emptyset:=X_{-1}\subset X_0\subset \cdots \subset X_n=X$$
by the closed algebraic (or analytic) subsets $X_i:=\bigcup_{\dim S\leq i} S$ ($0\leq i\leq n=\dim X$), with the strata $S$ of dimension $i$ given by the connected components of $X_i\backslash X_{i-1}$.

In addition, whenever two strata $S_1$ and $S_2$ are such that $S_2 \subseteq {\partial S_1}$, the pair $(S_2, S_1)$ is required to satisfy the following Whitney b-regularity conditions that guarantee that the variety $X$ is topologically or cohomologically equisingular along each stratum
(as in \cite[Section 4.2]{Sc}):
\begin{enumerate}
\item[(iii)] {\it Whitney b-condition}: If $x_n\in S_1$ and $y_n\in S_2$ are sequences converging to $x\in S_2$ such that the tangent planes $T_{x_n}S_1$ converge to some limiting plane $\tau$ and the secant lines $l_n=\overline{x_n,y_n}$ converge to some limiting line $l$ (in some local coordinates), then $l\subset \tau$.
\end{enumerate}
Note that the Whitney b-condition implies the following.
\begin{enumerate}
\item[(iv)] {\it Whitney a-condition}: If $x_n\in S_1$ is a sequence converging to $x\in S_2$ such that the tangent planes $T_{x_n}S_1$ converge to some limiting plane $\tau$, then $T_xS_2\subset \tau$.
\end{enumerate}
These conditions are independent of the choice of local coordinates, and any algebraic (or analytic) partition $\mathscr{S}$
as above has a refinement to an algebraic (or analytic) Whitney stratification (see, e.g., \cite[ Theorem 1.2 and Proposition 2.1]{Te}, \cite{Ve}, as well as the references given in \cite[Section 1.7]{GMs}). Moreover, if the Whitney b-conditions hold for all
connected components $S$ of the $X_i\backslash X_{i-1}$ and a filtration $X_{\bullet}$ by closed algebraic (or analytic) subsets as above,
then the partition $\mathscr{S}$ with these strata $S$ is (locally) finite and satisfies the frontier condition (see, e.g., 
\cite[Section 4.2.1]{Sc}) as well as the constructibility condition.

\bex[Whitney umbrella]\label{Wa}
The singular locus of the Whitney umbrella $$X=\{z^2=xy^2\}\subset \bC^3$$  is the $x$-axis, but the origin is ``more singular'' than any other point on the $x$-axis. A Whitney stratification of $X$ is given by the strata
\[
S_1=X \setminus \{x-\text{axis}\},\quad 
S_2=\{(x,0,0) \mid x \neq 0\},\quad
S_3=\{(0,0,0)\}.
\]
\eex
Another example is given by $X$ a complex manifold, with strata $S$ given by its connected components.\\

In the following, let $R$ be a noetherian and commutative ring of finite global dimension. Let $X$ be a complex algebraic (or analytic) variety, and denote by $D^b(X;R)$ the derived category of bounded complexes of sheaves of $R$-modules. By our assumptions on $X$ and $R$, these bounded derived categories
$D^b(-;R)$ are closed under Grothendieck's six operations: $Rf_*$, $Rf_!$, $f^*$, $f^!$, $R\cH om^{\bullet}$ and $\overset{L}{\otimes}$, with $f$ an algebraic (or analytic) morphism.

\bd  A sheaf $\cF$ of $R$-modules on $X$ is said to be \index{weakly constructible} {\it weakly constructible} if there is a Whitney stratification $\cS$ of $X$ so that the restriction 
$\cF\vert_S$ of $\cF$ to every stratum $S\in \cS$ is an $R$-local system (i.e., a locally constant sheaf). In this case we also say that $\cF$ is
{\it  $\cS$-weakly constructible}. A ($\cS$-)weakly constructible sheaf $\cF$ on $X$ is called \index{constructible} {\it ($\cS$-)constructible}, if all stalks 
$\cF_x$ for $x\in X$
are  finitely generated $R$-modules. A bounded complex $\cFb \in D^b(X;R)$ is called (weakly) constructible if all its cohomology sheaves $\cH^j(\cFb)$ are (weakly) constructible. Similarly for  $\cS$-(weakly) constructible complexes in case one  works with a fixed Whitney stratification $\cS$.
\ed
Note that  the category $Sh_{(\cS-)wc}(X)$   of ($\cS$-)weakly constructible sheaves on $X$ is an abelian subcategory of the category of all sheaves of $R$-modules on $X$, which is stable under extensions (like the subcategory of $R$-local systems). Similarly for the category $Sh_{(\cS-)c}(X)$  of
($\cS$-)constructible sheaves (like the subcategory of $R$-local systems with finitely generated stalks), since we assume $R$ to be noetherian.

\bex On a point space $X=\{pt\}$ any sheaf is weakly constructible. For general  $X$, the constant sheaf ${R}_X$ (resp., an $R$-local system $\cL$ on $X$)  is (weakly) constructible on $X$ with respect to any Whitney stratification. On the other hand, if $i\colon X \hookrightarrow \bC^*$ denotes the inclusion of the closed  analytic subset $X:=\{\frac{1}{n} \mid n\in \bN\}\subset \bC^*$, then the direct image sheaf $i_*{R}_X$ is constructible on 
$\bC^*$. But if $k\colon \bC^*\hookrightarrow \bC$ is the open inclusion, then the closure $\bar{X}=X\cup\{0\}$ in $\bC$ is not analytic and
 the direct image sheaf $k_*i_*{R}_X=(k\circ i)_*R_X$ is not weakly constructible on $\bC$.
\eex

We let $D^b_{(w)c}(X;R)$ be the full triangulated subcategory of $D^b(X;R)$ consisting of (weakly) constructible complexes (that is, complexes which are (weakly) constructible with respect to {\it some} Whitney stratification). We will also use the simpler notation $D^b_{(w)c}(X)$ if the coefficient ring $R$ is understood from the context. Similarly for $D^b_{\cS-(w)c}(X;R)$ in case we work with (weakly) constructible complexes with respect to a {\it fixed} Whitney stratification $\cS$. By viewing a sheaf as a complex concentrated in degree zero, one gets an isomorphism of Grothendieck groups (and similarly for weakly constructible sheaves)  \index{Grothendieck group}
\be \begin{CD}
K_0\left(Sh_{(\cS-)c}(X)\right) @> \sim >>K_0\left(D^b_{(\cS-)c}(X;R)\right)\:,
\end{CD}
\ee
whose inverse is given by taking the alternating sum of (classes of) cohomology sheaves (see \cite[Lemma 3.3.1]{Sc}).

\br Even in the complex analytic setting a sheaf $\cF$ (resp.,  a bounded sheaf complex $\cF^{\bullet}$) on $X$ is (weakly) constructible iff this is 
\index{locally constructible} locally the case, in the sense that there is an open covering $(U_i)$ of $X$ such that $\cF|_{U_i}$ (resp., $\cF^{\bullet}|_{U_i}$) is (weakly) constructible for all $i$
(see \cite[Proposition 4.1.13]{Di} and compare also with \cite[Corollary 3.4]{HS}). In fact, the complement of the largest open $U$ of $X$ such $\cF|_U$ is locally constant (resp., $\cF^{\bullet}|_U$ has locally constant cohomology sheaves) is then an  analytic subset of $X$ (see \cite[Proposition 4.1.12]{Di}).
\er

The derived category $D^b_{(w)c}(X;R)$ of bounded (weakly) constructible complexes is closed under Grothendieck's six operations: $Rf_*$, $Rf_!$, $f^*$, $f^!$, $R\cH om^{\bullet}$ and $\overset{L}{\otimes}$.
More precisely, one has the following (e.g., see \cite{Bo} or the unified treatment of \cite[Theorem 4.0.2 and Proposition 4.0.2]{Sc}):
\bt\label{cpres} Let $f\colon X \to Y$ be a morphism of complex algebraic (or analytic) varieties, with $j\colon  U\hookrightarrow Y$ the inclusion of the open complement of a closed algebraic (or analytic) subset.
\ben
\item[(a)] If $\cG^{\bullet} \in D^b_c(Y;R)$, then $f^*\cG^{\bullet}, f^!\cG^{\bullet} \in D^b_c(X;R)$.
\item[(b)] If $\cF^{\bullet} \in D^b_c(X;R)$ and $f$ is an algebraic map, then $Rf_*\cF^{\bullet}, Rf_!\cF^{\bullet} \in D^b_c(Y;R)$.\\ If $\cF^{\bullet} \in D^b_c(X;R)$ and $f$ is an analytic map so that the restriction of $f$ to ${\rm supp}(\cF^{\bullet})$ is proper (e.g., $f$ is proper), then $Rf_*\cF^{\bullet} \simeq Rf_!\cF^{\bullet} \in D^b_c(Y;R)$.
\item[(c)]  If $\cG^{\bullet} \in D^b_c(Y;R)$, then $Rj_*j^*\cG^{\bullet}, Rj_!j^*\cG^{\bullet} \in D^b_c(Y;R)$.
\item[(d)] If $\cF^{\bullet}, \cG^{\bullet} \in D^b_c(X;R)$, then $\cF^{\bullet} \overset{L}{\otimes} \cG^{\bullet},  R\cH om^{\bullet}(\cF^{\bullet},\cG^{\bullet}) \in D^b_c(X;R)$.
\een
Similarly for weakly constructible instead of constructible complexes.
\et

In fact, the above theorem in an application of the following more precise version, where one considers a stratified morphism with respect to given Whitney stratifications (e.g., see \cite{Bo} or \cite[Proposition 4.0.2 and Corollary 4.2.1]{Sc}):
\bt\label{cpres-strat} Let $f\colon X \to Y$ be a stratified morphism of complex algebraic (or analytic) varieties mapping all strata $S\in \cS$ of a Whitney stratification of $X$ into strata $T\in\cT$ of a Whitney stratification of $Y$, with $j\colon  U\hookrightarrow Y$ the inclusion of the open complement of a closed union of strata of $\cT$ (with its induced Whitney stratification $\cT|_U$).
\ben
\item[(a)] If $\cG^{\bullet} \in D^b_{\cT-c}(Y;R)$, then $f^*\cG^{\bullet}, f^!\cG^{\bullet} \in D^b_{\cS-c}(X;R)$.
\item[(b)] If $\cF^{\bullet} \in D^b_{\cS-c}(X;R)$ and $f$ is a \index{stratified submersion}
stratified submersion (i.e., it maps all strata $S\in \cS$ submersively to a stratum $T\in \cT$) so that the restriction of $f$ to ${\rm supp}(\cF^{\bullet})$ is proper (e.g., $f$ is proper), then $Rf_*\cF^{\bullet} \simeq Rf_!\cF^{\bullet} \in D^b_{\cT-c}(Y;R)$.
\item[(c)]  If $\cG^{\bullet} \in D^b_{\cT|_U-c}(U;R)$, then $Rj_*\cG^{\bullet}, Rj_!\cG^{\bullet} \in D^b_{\cT-c}(Y;R)$.
\item[(d)] If $\cF^{\bullet}, \cG^{\bullet} \in D^b_{\cS-c}(X;R)$, then $\cF^{\bullet} \overset{L}{\otimes} \cG^{\bullet},  R\cH om^{\bullet}(\cF^{\bullet},\cG^{\bullet}) \in D^b_{\cS-c}(X;R)$.
\een
Similarly for weakly constructible instead of constructible complexes.
\et

By considering  the constant map $c\colon  X\to pt$ to a point, one gets the following.
\bex
 Let $X$ be a complex algebraic (or analytic) variety with a Whitney stratification $\cS$, and let $i\colon  Z\hookrightarrow X$ be the inclusion of a closed algebraic (or analytic) subset given as a union of strata of $\cS$ with its induced Whitney stratification $\cS|_Z$.
\begin{enumerate}
\item The constant sheaf $R_X=c^*R$ is constructible with respect to $\cS$, with
$$H^k_c(X;R)=H^k(Rc_!c^*R) \quad \text{resp.,} \quad H^k(X;R)=H^k(Rc_*c^*R) \quad (k\in \bZ)$$
the corresponding \index{cohomology} \index{cohomology with compact support}
 cohomology (with compact support) of $X$.
\item The sheaf complex $i^!R_X=i^!c^*R$ is constructible with respect to $\cS|_Z$, with
$$H^k_Z(X;R)=H^k(Rc_*i^!c^*R)=H^k(X,X\backslash Z;R) \quad  (k\in \bZ)$$
the corresponding cohomology with support in $Z$ of $X$ (or relative cohomology of the pair $(X,X\backslash Z)$).
\item  \index{dualizing complex}  The {\it dualizing complex} $\bD_X^{\bullet}:=c^!R$ of $X$  is constructible with respect to $\cS$, with
$$H_k(X;R)=H^{-k}_c(X;\bD_X^{\bullet})=H^{-k}(Rc_!c^!R)\quad (k\in \bZ)$$
resp.,
$$H_k^{BM}(X;R)=H^{-k}(X;\bD_X^{\bullet})=H^{-k}(Rc_*c^!R) \quad (k\in \bZ)$$
the corresponding  \index{homology} \index{Borel-Moore homology}
(Borel-Moore) homology of $X$. If $X$ is smooth (or more generally an $R$-homology manifold) of pure dimension $d$, then 
$\bD_X^{\bullet}=c^!R\simeq c^*R[2d]$, which implies the \index{Poincar\'{e} duality}
{\it Poincar\'{e} duality}:
$$H_{2d-k}(X;R) \simeq H^k_c(X;R) \quad \text{and} \quad H_{2d-k}^{BM}(X;R) \simeq H^k(X;R) \quad (k\in \bZ)\:.$$
\end{enumerate}
\eex

\bc[External tensor product]   \index{external tensor product} Let $X_i$ be complex algebraic (or analytic) varieties, with 
 $\cF_i^{\bullet} \in D^b_{(w)c}(X_i;R)$ for $i=1,2$. Consider the external tensor product
$$\cF_1^{\bullet}\overset{L}{\boxtimes}  \cF_2^{\bullet}:=p_1^*\left(\cF_1^{\bullet}\right)\overset{L}{\otimes} 
p_2^*\left(\cF_2^{\bullet}\right)\in D^b_{(w)c}(X_1\times X_2;R)\:,$$
with $p_i\colon  X_1\times X_2\to X_i$ the projection on the corresponding factor for $i=1,2$. If $\cF_i^{\bullet}$ is (weakly) constructible
with respect to a Whitney stratification $\cS_i$ of $X_i$ for $i=1,2$, then $\cF_1^{\bullet}\overset{L}{\boxtimes}  \cF_2^{\bullet}$ is
  (weakly) constructible with respect to the product Whitney stratification $\cS_1\times \cS_2$ of $X_1\times X_2$ (with strata $S_1\times S_2$ for
$S_1\in \cS_1$ and $S_2\in \cS_2$).
\ec
The external tensor product of (weakly) constructible complexes behaves nicely with respect to products of two morphisms, as the following result shows
(see, e.g., \cite[Eq. (1.16) on p.78, Proposition 2.0.1 and Corollary 2.0.4]{Sc} for more general versions).

\bp[K\"{u}nneth isomorphisms] \index{K\"{u}nneth isomorphisms} Let $f_k\colon  X_k\to Y_k$ be two morphisms of complex algebraic (or analytic) varieties, with
$i_k\colon  Z_k\hookrightarrow X_k$ the inclusion of a closed complex algebraic (or analytic) subset, and $j_k\colon  U_k:=X_k\backslash Z_k\hookrightarrow X_k$
the inclusion of the open complement ($k=1,2$).
\begin{enumerate}
\item For any $\cF_k^{\bullet} \in D^b(Y_k;R)$ for $k=1,2$ one has:
$$(f_1\times f_2)^*\left(\cF_1^{\bullet}\overset{L}{\boxtimes}  \cF_2^{\bullet}\right) \simeq
f_1^*\left(\cF_1^{\bullet}\right)\overset{L}{\boxtimes}  f_2^*\left(\cF_2^{\bullet}\right)\:.$$
\item For any $\cF_k^{\bullet} \in D^b(X_k;R)$ for $k=1,2$ one has:
$$R(f_1\times f_2)_!\left(\cF_1^{\bullet}\overset{L}{\boxtimes}  \cF_2^{\bullet}\right) \simeq
Rf_{1!}\left(\cF_1^{\bullet}\right)\overset{L}{\boxtimes}  Rf_{2!}\left(\cF_2^{\bullet}\right)\:.$$
\item For any $\cF_k^{\bullet} \in D^b_{wc}(Y_k;R)$ for $k=1,2$ one has:
$$(i_1\times i_2)^!\left(\cF_1^{\bullet}\overset{L}{\boxtimes}  \cF_2^{\bullet}\right) \simeq
i_1^!\left(\cF_1^{\bullet})\right)\overset{L}{\boxtimes}  i_2^!\left(\cF_2^{\bullet}\right)\:.$$
\item For any $\cF_k^{\bullet} \in D^b_{wc}(X_k;R)$ for $k=1,2$ one has in the complex algebraic context:
$$R(f_1\times f_2)_*\left(\cF_1^{\bullet}\overset{L}{\boxtimes}  \cF_2^{\bullet}\right) \simeq
Rf_{1*}\left(\cF_1^{\bullet}\right)\overset{L}{\boxtimes}  Rf_{2*}\left(\cF_2^{\bullet}\right)\:.$$
\item For any $\cF_k^{\bullet} \in D^b_{wc}(X_k;R)$ for $k=1,2$ one has:
$$R(j_1\times j_2)_*(j_1\times j_2)^*\left(\cF_1^{\bullet}\overset{L}{\boxtimes}  \cF_2^{\bullet}\right) \simeq
Rj_{1*}j_1^*\left(\cF_1^{\bullet}\right)\overset{L}{\boxtimes}  Rj_{2*}j_2^*\left(\cF_2^{\bullet}\right)\:.$$
\end{enumerate}
\ep
By taking for both morphisms $f_k$ the constant map $f_k=c: X_k\to pt$ to a point space ($k=1,2$) one gets the following.
\bex[Classical K\"{u}nneth formulae] \index{K\"{u}nneth formulae}
For any $\cF_k^{\bullet} \in D^b_{wc}(X_k;R)$ for $k=1,2$ one has:
$$R\Gamma_c\left(X_1\times X_2;\cF_1^{\bullet}\overset{L}{\boxtimes}  \cF_2^{\bullet}\right) \simeq
R\Gamma_c\left(X_1;\cF_1^{\bullet}\right)\overset{L}{\otimes}   R\Gamma_c\left(X_2;\cF_2^{\bullet}\right) \:,$$
$$R\Gamma_{Z_1\times Z_2}\left(X_1\times X_2;\cF_1^{\bullet}\overset{L}{\boxtimes}  \cF_2^{\bullet}\right) \simeq
R\Gamma_{Z_1}\left(X_1;\cF_1^{\bullet}\right)\overset{L}{\otimes}   R\Gamma_{Z_2}\left(X_2;\cF_2^{\bullet}\right) \:,$$
and in the complex algebraic context:
$$R\Gamma\left(X_1\times X_2;\cF_1^{\bullet}\overset{L}{\boxtimes}  \cF_2^{\bullet}\right) \simeq
R\Gamma\left(X_1;\cF_1^{\bullet}\right)\overset{L}{\otimes}   R\Gamma\left(X_2;\cF_2^{\bullet}\right) \:.$$
In the special case when $R$ is a field, one further gets
$$H^*_c\left(X_1\times X_2;\cF_1^{\bullet}\boxtimes  \cF_2^{\bullet}\right) \simeq
H^*_c\left(X_1;\cF_1^{\bullet}\right)\otimes   H^*_c\left(X_2;\cF_2^{\bullet}\right) \:,$$
and in the complex algebraic context also
$$H^*\left(X_1\times X_2;\cF_1^{\bullet}\boxtimes  \cF_2^{\bullet}\right) \simeq
H^*\left(X_1;\cF_1^{\bullet}\right)\otimes   H^*\left(X_2;\cF_2^{\bullet}\right) \:.$$
\eex

Another important application of the general calculus of (weakly) constructible complexes deals with {\it Verdier duality}
(see, e.g.,  \cite[Corollary 4.2.2]{Sc}):
\bc Let $X$ be a complex algebraic (or analytic) variety with a Whitney stratification $\cS$, with  morphisms $f\colon Z\to X$ and $g\colon X\to Y$.
\ben
\item  If $\cF^{\bullet} \in D^b(X;R)$ is $\cS$-(weakly) constructible, then its \index{Verdier dual}  {\it Verdier dual} 
$$\cD \cF^{\bullet}:= R\cH om^{\bullet}(\cF^{\bullet},\bD_X^{\bullet}) $$ 
is also $\cS$-(weakly) constructible, with
$$f^!\left(\cD \cF^{\bullet}\right) \simeq \cD\left( f^*\cF^{\bullet}\right) \quad \text{resp.,} \quad
Rg_*\left(\cD \cF^{\bullet}\right) \simeq \cD\left( Rg_!\cF^{\bullet}\right) \:.$$
\item If $\cF^{\bullet} \in D^b(X;R)$ is constructible, then \index{biduality} biduality holds:
\be \cF^{\bullet} \simeq \cD \cD \cF^{\bullet}
\ee
so that $\cF^{\bullet} \in D^b(X;R)$ is ($\cS$-)constructible iff its Verdier dual  $\cD \cF^{\bullet} \in D^b(X;R)$ is ($\cS$-)constructible.
Moreover
$$f^*\left(\cD \cF^{\bullet}\right) \simeq \cD\left( f^!\cF^{\bullet}\right) \quad \text{resp.,} \quad
Rg_!\left(\cD \cF^{\bullet}\right) \simeq \cD\left( Rg_*\cF^{\bullet}\right) \:.$$
Here, for the last isomorphism one has to assume that $Rg_!\left(\cD \cF^{\bullet}\right) $ is constructible (e.g., as in the algebraic context).
\een
\ec

Note that already  for a point space $X=\{pt\}$ the biduality result uses our assumption that  $R$ is a noetherian and commutative ring of finite global dimension (see also \cite[Exercise I.30]{KS}). Similarly, Verdier duality commutes for constructible sheaf complexes with external tensor products
(see, e.g., \cite[Corollary 2.0.4]{Sc} in a more general context).

\bp Let $X_i$ be  a complex algebraic (or analytic) variety with  $\cF^{\bullet}_i \in D^b_c(X;R)$ for $i=1,2$. Then
$$\cD\left( \cF_1^{\bullet}\overset{L}{\boxtimes}  \cF_2^{\bullet}\right) \simeq
\cD\left(  \cF_1^{\bullet}\right) \overset{L}{\boxtimes} \cD\left( \cF_2^{\bullet}\right) \:.$$
In particular,  $\bD_{X_1\times X_2}^{\bullet}\simeq \bD_{X_1}^{\bullet}\overset{L}{\boxtimes}  \bD_{X_2}^{\bullet}$, so that one gets the classical K\"{u}nneth formula for homology, e.g., for $R$ a field:
$$H_*(X_1\times X_2;R) \simeq H_*(X_1;R) \otimes H_*(X_2;R)\:,$$
 and in the algebraic context also for Borel Moore homology, e.g., for $R$ a field:
$$H^{BM}_*(X_1\times X_2;R) \simeq H^{BM}_*(X_1;R) \otimes H^{BM}_*(X_2;R)\:.$$
\ep

The general calculus of constructible sheaves also includes finiteness results for the cohomology (with compact support) of constructible sheaf complexes.

\bc\label{finc} Assume that $\cF^{\bullet} \in D^b_c(X;R)$ and that either
\ben
\item[(a)] $X$ is a complex algebraic variety, or
\item[(b)] $X$ is an analytic space and ${\rm supp}(\cF^{\bullet})$ is compact. 
\een
 Then the hypercohomology groups $H^i(X;\cF^{\bullet})$ and $H_c^i(X;\cF^{\bullet})$ are finite type $R$-modules for every $i \in \Z$ (which are zero for $|i|$ large enough).
\item[(c)] Assume $X$ is a compact analytic space with a Whitney stratification $\cS$ so that $j\colon  U\hookrightarrow X$ is the inclusion of the open complement of a closed union of strata. If  $\cG^{\bullet} \in D^b_{\cS|_U-c}(U;R)$, then also
the hypercohomology groups $H^i(U;\cG^{\bullet})$ and $H_c^i(U;\cG^{\bullet})$ are finite type $R$-modules for every $i \in \Z$ (which are zero for $|i|$ large enough).
\ec

With $\cF^{\bullet} \in D^b_c(X;R)$ as in the above corollary, we make the following.
\bd Assume $R$ is a field. The {\it (compactly supported) Euler characteristic \index{Euler characteristic} of $\cF^{\bullet} \in D^b_c(X;R)$} is defined as:
$$\chi_{(c)}(X,\cF^{\bullet}):=\chi(H^*_{(c)}(X;\cF^{\bullet})):=\sum_{i \in \Z} (-1)^i \dim_R H_{(c)}^i(X;\cF^{\bullet}).$$
(Here, we use the notation $\chi_{(c)}$ and $H_{(c)}^i$ to indicate that the definition applies to the compactly supported Euler characteristic $\chi_{c}$ by using $H_{c}^i$, as well as to the usual Euler characteristic $\chi$ by using $H^i$.)
\ed
As it will be explained  later on (see  Example \ref{ex:analconstfct}, and also \cite[Section 2.3 and Section 6.0.6]{Sc}), in this complex context we have the equality:
\be
\chi(X,\cF^{\bullet})=\chi_{c}(X,\cF^{\bullet}).
\ee
Moreover, this Euler characteristic depends only on the associated constructible function 
\be
\chi_{stalk}\left(\cF^{\bullet}\right)\in CF(X)
\ee 
given by the stalkwise Euler characteristic $\chi_{stalk}\left(\cF^{\bullet}\right)(x):=\chi\left(\cF^{\bullet}_x\right)$
for $x\in X$, with $CF(X)$ the corresponding abelian group of \index{constructible function} constructible functions given by (locally) finite $\bZ$-linear combinations of indicator functions $1_Z$, for $Z\subset X$ a closed irreducible algebraic (or analytic)
subset of $X$. Similarly for the abelian group $CF_{\cS}(X)$ of $\cS$-constructible functions given by (locally) finite $\bZ$-linear combinations of indicator functions $1_S$ or $1_{\bar{S}}$ for $S\in \cS$ a stratum, i.e., $\bZ$-valued functions which are constant on all strata $S\in\cS$. This
induces a surjective group homomorphism
\be 
\chi_{stalk} \colon \: K_0\left(D^b_{(\cS-)c}(X;R)\right)\to CF_{(\cS)}(X)\:,
 \ee
with 
\be
\chi_{c}(X,\cF^{\bullet})=\sum_{S\in \cS}\: \chi_c(S)\cdot \chi_{stalk}\left(\cF^{\bullet}\right)(S)
\ee
for an $\cS$-constructible complex  $\cF^{\bullet}$ on $X$ in the complex algebraic context (or analytic context with $X$ compact).
Here $\chi_c(S):=\chi_c(R_S)=\chi(H_c^*(S;R))$ is the corresponding Euler characteristic of a stratum $S\in \cS$.\\

Let us also mention here the following (co)stalk calculation.
\bp\label{hil2}
 Let ${\mathscr{F}^{\bullet}}\in D_c^b(X;R)$, $x \in X$, and $i_x\colon \{x\} \hookrightarrow X$ the inclusion. Then
\be\label{lequ1} \mathscr{H}^j({\mathscr{F}^{\bullet}})_x \simeq {H}^j(i_x^*{\mathscr{F}^{\bullet}}) \simeq H^j(\mathring{B}_{\epsilon,x};\mathscr{F}^{\bullet}),\ee
\be\label{lequ2}  {H}^j(i_x^!{\mathscr{F}^{\bullet}}) \simeq {H}_c^j(\mathring{B}_{\epsilon,x};\mathscr{F}^{\bullet}) \simeq {H}^j(\mathring{B}_{\epsilon,x},
\mathring{B}_{\epsilon,x}\setminus {x};\mathscr{F}^{\bullet}),\ee
where $\mathring{B}_{\epsilon,x}$ is the intersection of $X$ with an open small $\epsilon$-ball neighborhood of $x$ in some local embedding
of $X$ in $\mathbb{C}^N$. Here, $i_x^*\cFb$ and $i_x^!\cFb$ are called the {\it stalk} and, respectively, {\it costalk} of $\cFb$ at $x$. \index{stalk} \index{costalk}
\ep

In fact, this will be an easy application of the Morse theoretical results explained later on, since the proper real analytic function $r$ given by the squared distance to $x$
(in these local coordinates) has no stratified  critical values in an interval $\,]0,\epsilon[\,$ for $\epsilon>0$ small enough 
(see Lemma \ref{loc.const} and, e.g., \cite[Lemma 5.1.1]{Sc}).

\subsection{Perverse sheaves}

Perverse sheaves are an important class of constructible complexes, introduced in \cite{BBD} as a formalization of Kashiwara's Riemann--Hilbert correspondence \cite{Ka} (see also \cite{HTT}), which relates the topology of complex algebraic, resp., analytic  varieties (intersection homology) and the algebraic, resp., analytic theory of differential equations (holonomic $D$-modules). We recall their definition below. 

\bd\label{sper} 
\noindent(a) \ The \index{perverse t-structure} {\it perverse t-structure} on $D_c^b(X;R)$ consists of the two strictly full  subcategories ${^pD}^{\leq 0}(X;R)$ and ${^pD}^{\geq 0}(X;R)$ of $D_c^b(X;R)$ defined as:
$${^pD}^{\leq 0}(X;R):=\{ \cFb \in D^b_c(X;R) \mid \dim {\rm supp}^{-j}(\cFb) \leq j, \forall j \in \bZ \},$$
$${^pD}^{\geq 0}(X;R):=\{ \cFb \in D^b_c(X;R) \mid \dim {\rm cosupp}^{j}(\cFb) \leq j, \forall j \in \bZ \},$$
where, for $i_x\colon \{x\} \hookrightarrow X$ denoting the point inclusion, we define the \index{support} {\it $j$-th support} and, respectively, \index{cosupport} the {\it $j$-th cosupport} of $\cFb \in D^b_c(X;R)$ by:
$${\rm supp}^{j}(\cFb) = \overline{\{x \in X \mid H^j(i_x^*\cFb) \neq 0\}},$$
$${\rm cosupp}^{j}(\cFb) = \overline{\{x \in X \mid H^j(i_x^!\cFb) \neq 0\}}.$$
(For a constructible complex $\cFb$, the sets $\supp^{j}(\mathscr{F}^{\bullet})$ and $\cosupp^{j}(\mathscr{F}^{\bullet})$ are closed algebraic (or analytic) subvarieties of $X$, hence their dimensions are well defined.) \\
\noindent(b) For a given Whitney stratification $\cS$ of $X$ this also induces the {\it perverse t-structure} on $D_{\cS-c}^b(X;R)$ with
${^pD}_{\cS}^{\leq 0}(X;R)$, resp., ${^pD}_{\cS}^{\geq 0}(X;R)$ defined by the same (co)support conditions.\\
\noindent(c) A  ($\cS$-)constructible complex $\cFb \in D^b_{(\cS-)c}(X;R)$ is called a \index{perverse sheaf} {\it perverse sheaf} on $X$ if 
$$\cFb \in Perv_{(\cS)}(X;R):={^pD}_{(\cS)}^{\leq 0}(X;R) \cap {^pD}_{(\cS)}^{\geq 0}(X;R).$$
\ed

The category of  ($\cS$-constructible) perverse sheaves is the {\it heart} of the perverse t-structure, hence it is an abelian category, and it is stable by extensions (see, e.g., \cite[Theorem 1.3.6]{BBD}). 

\br (a) The same definition also defines the {\it perverse t-structure} for ($\cS$-)\\ weakly constructible complexes with  heart the abelian category of
 ($\cS$-)weakly constructible perverse sheaves.\\
\noindent(b) An algebraically (weakly) constructible complex is perverse if and only if it is so when viewed as  an analytically (weakly) constructible complex.\\
\noindent(c) 
If $R$ is a {\it field}, the Universal Coefficient Theorem can be used to show that the Verdier duality functor $\cD\colon D^b_c(X;R) \to D^b_c(X;R)$ satisfies:
\be {\rm cosupp}^{j}(\cFb)={\rm supp}^{-j}(\cD\cFb),
\ee
In particular, $\cD$ exchanges  ${^pD}^{\leq 0}(X;R)$ and ${^pD}^{\geq 0}(X;R)$, so that it
preserves ($\cS$-) constructible perverse sheaves with field coefficients.
\er

Recall here that two subcategories ${^pD}^{\leq 0}(X;R)$ and ${^pD}^{\geq 0}(X;R)$ of $D_{(\cS-w)c}^b(X;R)$ define a {\it t-structure} just means (see, e.g., \cite[Definition 1.3.1]{BBD}):
\begin{enumerate}
\item $\Hom_{D^b(X;R)}(\cFb,\cGb[-1]) =0$ for all $\cFb\in {^pD}^{\leq 0}(X;R)$ and $\cGb \in {^pD}^{\geq 0}(X;R)$.
\item  ${^pD}^{\leq 0}(X;R)$ is stable under $[1]$, and ${^pD}^{\geq 0}(X;R)$ is stable under $[-1]$.
\item For any $\cEb\in D_{(\cS-w)c}^b(X;R)$ there is a distinguished triangle
$$ \cFb \lra \cEb \lra \cGb[-1] \overset{[1]}{\lra}$$
for some $\cFb\in {^pD}^{\leq 0}(X;R)$ and $\cGb \in {^pD}^{\geq 0}(X;R)$.
\end{enumerate}

Then it is enough to check these properties for a fixed Whitney stratification $\cS$, where they can be proved by induction on $\dim X$ via the 
{\it gluing of t-structures} as in \cite[Corollary  2.1.4, Proposition 2.1.14]{BBD}. 
Here it is important to note that the conditions  $\cFb \in {^pD}^{\leq 0}(X;R)$ and, resp.,  $\cFb \in {^pD}^{\geq 0}(X;R)$ can also be described in terms of a fixed Whitney stratification of $X$ for which $\cFb$ is $\cS$-(weakly) constructible. Indeed, the perverse t-structure can be characterized as follows:
\bt\label{pervstr}
Assume $\cFb \in D^b_{(w)c}(X;R)$ is (weakly) constructible with respect to a Whitney stratification $\cS$ of $X$. Then:
\begin{itemize}
\item[(i)] stalk vanishing: \index{stalk vanishing}
\begin{center}
$\cFb \in {^pD}^{\leq 0}(X;R)$ $\iff$ $\forall S \in \cS$,  $\forall x \in S$: $H^j(i_x^*\cFb)=0$ for all $j>-\dim S$.
\end{center}
\item[(ii)] $\,$ costalk vanishing: \index{costalk vanishing}
\begin{center}
$\cFb \in {^pD}^{\geq 0}(X;R)$ $\iff$ $\forall S \in \cS$,  $\forall x \in S$: $H^j(i_x^!\cFb)=0$ for all $j<\dim S$.
\end{center}
\end{itemize}
\et

Of course the isomorphism class of the stalk  $i_x^*\cFb$ for $x\in S$ as above does not depend on the choice of the point $x\in S$. In fact, if $i_S\colon  S\hookrightarrow X$ denotes the inclusion of the stratum $S$, the stalk vanishing condition is equivalent to 
\begin{itemize}
\item[(i') ]\ \  $\cFb \in {^pD}^{\leq 0}(X;R)$ $\iff$ $\forall S \in \cS$: $\cH^j(i_S^*\cFb)=0$ for all $j>-\dim S$.
\end{itemize}
Similarly for the costalk $i_x^!\cFb$ for $x\in S$, since this costalk is also isomorphic to $k_x^!i_S^!\cFb$ with $k_x\colon  \{ x \} \hookrightarrow S$ the inclusion of the point $\{x\}$ into $S$. But $i_S^!\cFb$ has locally constant cohomology sheaves so that
$k_x^!i_S^!\cFb[2\dim S]\simeq k_x^*i_S^!\cFb$ for all $x\in S$. And then the costalk condition is equivalent to
\begin{itemize}
\item[(ii')] \ \ 
$\;\cFb \in {^pD}^{\geq 0}(X;R)$ $\iff$ $\forall S \in \cS$: $\cH^j(i_S^!\cFb)=0$ for all $j<-\dim S$.
\end{itemize}

\bex\label{ct1} Assume $X$ is of pure complex dimension with $c\colon X \to pt$ the constant map to a point space. As usual, for an $R$-module $M$ we denote by $M_X=c^*M$ the constant $R$-sheaf on $X$ with stalk $M_x=M$ for all $x\in X$. Then:
\begin{itemize}\label{e47}
\item[(a)] \ ${M}_X[\dim X] \in {^pD}^{\leq 0}(X;R)$. More generally, if $\cL$ is a local system on $X$, then $\cL[\dim X]  \in {^pD}^{\leq 0}(X;R)$
(since this is a local condition).
\item[(b)] \ If $X$ is smooth with the trivial stratification, and $\cL$ is a local system on $X$, then $\cL[\dim X]$ is perverse on $X$
(since in this case  $i_x^!\cL [2\dim X]\simeq i_x^*\cL$ for all $x\in X$).
\item[(c)] \ The \index{intersection complex} {\it intersection ($IC$) complexes} on $X$ of Goresky-MacPherson \cite{GM} (for the middle perversity) are examples of perverse sheaves. For a given Whitney stratification $\cS$ of $X$ these are defined by further imposing, for all strata $S\in \cS$ with $\dim S<\dim X$, the stronger stalk/costalk vanishing conditions in Theorem \ref{pervstr} obtained by replacing  $>$ with $\geq$ and $<$ with $\leq$. Such an $IC$ complex on $X$ is $\cS$-weakly constructible and determined by its restriction to the top dimensional strata, which is (isomorphic to) a shifted local system $\cL[\dim X]$, so we may denote it unambiguously  by $IC_X(\cL)$. It is $\cS$-constructible if $\cL$ has finitely generated stalks. If $\cL$ is the constant sheaf $R$ on the top dimensional strata, we use the notation $IC_X$. In the end, these perverse $IC$-complexes do not depend on the chosen Whitney stratification $\cS$, but only on the 
{\it generically defined}  local system $\cL$ (on the complement of a closed algebraic or analytic subset of $X$ of dimension smaller than 
$\dim X$).
\item[(d)] \ If $X$ is a local complete intersection then ${R}_X[\dim X]$ is a perverse sheaf on $X$ (see Example \ref{purity} and, e.g.,
\cite[Example 6.0.11]{Sc}). More generally, if $\cL$ is a local system on $X$, then $\cL[\dim X]$ is perverse on $X$. By (a) one only has to show that  $\cL[\dim X]  \in {^pD}^{\geq 0}(X;R)$. Since this is a local condition in the classical topology, we can assume $\cL=M_X$ is a constant sheaf, with $i\colon  X\hookrightarrow X'$ the inclusion of a closed analytic subset defined as the zero set of $k$ holomorphic functions on a pure dimensional complex manifold $X'$ with $\dim X'-k= \dim X$. Then $M_X=i^*M_{X'}$ with 
$M_{X'}[\dim X']  \in {^pD}^{\geq 0}(X';R)$ by (b), so that the claim follows from Proposition \ref{prop:subspace} below (see, e.g.,  \cite[Proposition 6.0.2]{Sc}).
\end{itemize}
\eex

To simplify the formulation of some results, let us recall the following.
\bd For $n\in \bZ$, define the following two strictly full  subcategories ${^pD}^{\leq n}_{(\cS)}(X;R)$ and ${^pD}^{\geq n}_{(\cS)}(X;R)$ of 
$D_{(\cS-)c}^b(X;R)$, 
resp., $D_{(\cS-)wc}^b(X;R)$ by
$${^pD}^{\leq n}_{(\cS)}(X;R):={^pD}^{\leq 0}_{(\cS)}(X;R)[-n] \quad \text{and} \quad {^pD}^{\geq n}_{(\cS)}(X;R):={^pD}^{\geq 0}_{(\cS)}(X;R)[-n]$$
so that ${^pD}^{\leq -1}_{(\cS)}(X;R)\subset {^pD}^{\leq 0}_{(\cS)}(X;R)$ and ${^pD}^{\geq 1}_{(\cS)}(X;R)\subset {^pD}^{\geq 0}_{(\cS)}(X;R)$.
\ed

\br The two subcategories ${^pD}^{\leq n}_{(\cS)}(X;R)$ and ${^pD}^{\geq n}_{(\cS)}(X;R)$ define a \index{shifted perverse t-structure}
{\it shifted perverse t-structure} on
$D_{(\cS-)c}^b(X;R)$
and  $D_{(\cS-)wc}^b(X;R)$. But for most of the results of this paper (especially those proved by stratified Morse theory later on), we only need the following obvious properties (as in \cite[Chapter VI]{Sc}):
\begin{enumerate}
\item The zero object belongs to ${^pD}^{\leq n}_{(\cS)}(X;R)$ and ${^pD}^{\geq n}_{(\cS)}(X;R)$.
\item ${^pD}^{\leq n}_{(\cS)}(X;R)$ and ${^pD}^{\geq n}_{(\cS)}(X;R)$ are stable by extensions.
\item  ${^pD}^{\leq n}(X;R)$ is stable under $[1]$, and ${^pD}^{\geq n}(X;R)$ is stable under $[-1]$.
\end{enumerate}
\er

\bex On a point space $X=\{pt\}$ one gets:
\begin{enumerate}
\item ${^pD}^{\leq n}(\{pt\};R)$ is given by the bounded complexes of $R$-modules, whose cohomology is concentrated in degree $\leq n$.
This condition is stable under shifting a complex to the left (i.e., the shift $[1]$).
\item ${^pD}^{\geq n}(\{pt\};R)$ is given by the bounded complexes of $R$-modules, whose cohomology is concentrated in degree $\geq n$.
This condition is stable under shifting a complex to the right (i.e., the shift $[-1]$).
\end{enumerate}
\eex

\bl Let $f\colon  X\to Y$ be a morphism of complex algebraic (or analytic) varieties, whose \index{fiber dimension}
fiber dimension is bounded by $d\in \bN_0$. Then:
\begin{itemize}
\item[(a)] \ $f^!$ maps ${^pD}^{\geq n}(Y;R)$ into ${^pD}^{\geq n-d}(X;R)$.
\item[(b)] \ $f^*$ maps ${^pD}^{\leq n}(Y;R)$ into ${^pD}^{\leq n+d}(X;R)$.
\end{itemize}
\el
This follows  from the definitions using $i_x^!f^!\simeq i_{f(x)}^!$ and $i_x^*f^*\simeq i_{f(x)}^*$ for all $x\in X$.
As an example, one can take for $f$ a locally closed inclusion or an unramified covering map, both of which have fiber dimension $d=0$.

\bex Let $X$ be a complex algebraic (or analytic) variety of dimension $d=\dim X$, with $c\colon  X\to pt$ a constant map. Then
$$\bD_X^{\bullet}=c^!R \in {^pD}^{\geq -d}(X;R) \quad \text{and} \quad R_X=c^*R\in {^pD}^{\leq d}(X;R) \:.$$
\eex

\bc
Let $f\colon  X\to Y$ be a \index{smooth morphism} \index{submersion}
{\it smooth} morphism (i.e., a submersion) of complex algebraic (or analytic) varieties, with constant relative (or fiber) dimension 
$d\in \bN_0$. Then $f^!\simeq f^*[2d]$ so that
\begin{itemize}
\item[(a)] \ $f^*$ maps ${^pD}^{\geq n}(Y;R)$ into ${^pD}^{\geq n+d}(X;R)$.
\item[(b)] \ $f^!$ maps ${^pD}^{\leq n}(Y;R)$ into ${^pD}^{\leq n-d}(X;R)$.
\end{itemize}
In particular  $f^![-d]\simeq f^*[d]$ maps $Perv(Y;R)$ into $Perv(X;R)$.
\ec

\bex Let $f\colon  X\to Y$ be a {\it smooth} morphism (i.e., a submersion) of complex algebraic (or analytic) varieties, with constant relative (or fiber) dimension 
$d$. Assume $f$ is surjective and $Y$ (and then also $X$) is pure dimensional. Then
$$f^*IC_Y(\cL)[d]\simeq IC_X(f^*\cL) $$
for a generically defined local system $\cL$ on $Y$, with $f^*\cL$ the corresponding generically defined local system  on $X$ defined by pullback.
\eex

The following result will be very important for the applications of the {\it stratified Morse theory for constructible sheaves} in the next sections.

\bp\label{prop-tr} Let $Y\hookrightarrow M$ be a closed complex algebraic (or analytic) subvariety of an ambient complex algebraic (or analytic) manifold 
$M$.
Assume $N\hookrightarrow M$ is a closed complex algebraic (or analytic) submanifold of constant codimension $d=\dim M -\dim N$, which is 
\index{transversal intersection}
{\it transversal} to a Whitney stratification $\cS$ of $Y$
(i.e., $N$ is transversal to all strata $S\in \cS$). Then  $X:=Y\cap N$ gets an induced Whitney stratification $\cS'$ with strata $S'$ the connected components of the intersections $S\cap N$ for $S\in \cS$, with $\dim S\cap N= \dim S -d$ for all $S\in \cS$ (and $S\cap N\neq \emptyset$). Let $i\colon  X=Y\cap N\hookrightarrow Y$ be the (stratified) closed inclusion. Then
\begin{itemize}
\item[(a)] \ $i^*$ maps ${^pD}_{\cS}^{\leq n}(Y;R)$ into ${^pD}^{\leq n-d}_{\cS'}(X;R)$.
\item[(b)] \ $i^*$ maps ${^pD}^{\geq n}_{\cS}(Y;R)$ into ${^pD}^{\geq n-d}_{\cS'}(X;R)$.
\end{itemize}
In particular  $i^*[-d]$ maps $Perv_{\cS}(Y;R)$ into $Perv_{\cS'}(X;R)$.
\ep

\begin{proof} Consider the following cartesian diagram of closed inclusions for $S\in \cS$:
$$\begin{CD}
S':=S\cap N @> i_{S'} >> X=Y\cap N \\
@V i' VV @VV i V \\
S @> i_S >> Y \:.
\end{CD}$$
Then (a) follows from $i^*_{S'}i^*\simeq i'^*i^*_S$ for checking the stalk vanishing condition (i').
Similarly (b) follows from the {\it base change isomorphism}
 \be\label{bcs}
i^!_{S'}i^*\cFb\simeq i'^* i^!_S\cFb \quad \text{for $\cFb \in D^b_{\cS-wc}(X;R)$}
\ee
 for checking the costalk vanishing condition (ii').
\end{proof}

\bex\label{IC-transv} Consider the context of Proposition \ref{prop-tr}, with $Y$ (and therefore also $X=Y\cap N$) pure dimensional. Let $\cL$ be a local system on the open subset $U$ of $Y$ given by the top dimensional stratum of $\cS$ (which is then dense in $Y$). Similarly $U\cap N$ is open and dense in $X$. Then
$$i^*IC_Y(\cL)[-d]\simeq IC_X(i'^*\cL) \:,$$
with $i'\colon  U\cap N \to U$ the induced inclusion.
\eex

In the base change isomorphism (\ref{bcs}) used above, we can even assume that $S\hookrightarrow Y$ is a {\it closed} stratum, by restriction to the open complement of $\partial S$. Then it is a special case (with $Y'=S$) of the following more general result 
(see, e.g., \cite[Proposition 4.3.1 and Remark 4.3.6]{Sc}).

\bt[Base change isomorphisms]\label{bci} 
Let $Y\hookrightarrow M$ be a closed complex algebraic (or analytic) subvariety of an ambient complex algebraic (or analytic) manifold $M$.
Assume $N\hookrightarrow M$ is a closed complex algebraic (or analytic) submanifold, which is \it{transversal} to a Whitney stratification $\cS$ of $Y$.
Let $i\colon  Y'\hookrightarrow Y$ be the inclusion of a closed union of strata of $\cS$, with $j\colon  U:=Y\backslash Y'\hookrightarrow Y$ the inclusion of the open complement with its induced stratification $\cS|_U$. Consider the cartesian diagram
$$\begin{CD}
X':=Y'\cap N @> i' >> X:=Y\cap N @< j' << U':=U\cap N\\
@VV k' V @VV k V @VV k'' V \\
Y' @> i >> Y @< j << U \:.
\end{CD}$$
Then one has the following  \index{base change isomorphism} {\it base change isomorphisms}:
\be\label{bc-open}
k^*Rj_*\cFb \simeq Rj'_*k''^*\cFb \quad \text{for $\cFb \in D^b_{\cS|_U-wc}(U;R)$}
\ee
and
\be
k'^*i^!\cFb \simeq i'^!k^*\cFb \quad \text{for $\cFb \in D^b_{\cS-wc}(Y;R)$.}
\ee
\et

Next we study the relation between the {\it perverse t-structure} and {\it external tensor products}.

\bp\label{p1230}  Let $X_i$ be  complex algebraic (or analytic) varieties ($i=1,2$). Then 
\begin{itemize}
\item[(a)] \ The external tensor product $\overset{L}{\boxtimes}$ induces
$$\overset{L}{\boxtimes}:\:{^pD}^{\leq n}(X_1;R)\times {^pD}^{\leq m}(X_2;R)\to
{^pD}^{\leq n+m}(X_1\times X_2;R) \:.$$
\item[(b)] \ Assume $R$ is a field. Thent $\overset{L}{\boxtimes}$ induces
$$\overset{L}{\boxtimes}:\:{^pD}^{\geq n}(X_1;R)\times {^pD}^{\geq m}(X_2;R)\to
{^pD}^{\geq n+m}(X_1\times X_2;R) \:.$$
\end{itemize}
In particular, if $R$ is a field, $\overset{L}{\boxtimes}$ induces
$$\overset{L}{\boxtimes}:\:  Perv(X_1;R)\times Perv(X_2;R) \to Perv(X_1\times X_2;R)\:.$$
\ep
Property (a) follows from $i_{(x_1,x_2)}^*(-\overset{L}{\boxtimes}-)\simeq i_{x_1}^*(-)\overset{L}{\otimes}i_{x_2}^*(-)$ and the right exactness
of the tensor product $\otimes$. Property (b) is a consequence of  the K\"{u}nneth isomorphism $i_{(x_1,x_2)}^!(-\overset{L}{\boxtimes}-)\simeq i_{x_1}^!(-)\overset{L}{\otimes}i_{x_2}^!(-)$ and the exactness
of the tensor product $\otimes$ for $R$ a field.

\bex\label{ex1230} Let $X_i$ be  pure dimensional complex algebraic (or analytic) varieties ($i=1,2$), with $R$ a field. Then 
$$IC_{X_1}(\cL_1)\overset{L}{\boxtimes}IC_{X_2}(\cL_2)\simeq IC_{X_1\times X_2}(\cL_1\boxtimes \cL_2)$$
for a generically defined local system $\cL_i$ on $X_i$ ($i=1,2$).
\eex

\bex\label{ex1232}
Let $X$ be a  complex algebraic (or analytic) variety, with $R$ a field. 
If $\cFb \in Perv(X;R)$ and $\cL$ is a locally constant sheaf on $X$, then 
$\cFb \otimes \cL \in Perv(X;R)$. In particular, if $X$ is also pure dimensional, one gets:
$$IC_X(\cL_1)\otimes \cL_2\simeq IC_X(\cL_1\otimes \cL_2),$$
with $\cL_1$ a generically defined local system on $X$, and $\cL_2$ a local system on all of $X$. As a special case, one also gets:
$$IC_X \otimes \cL \simeq IC_X(\cL)$$
if $X$ is pure dimensional and $\cL$ is a local system on $X$.
These assertions can be checked locally in the analytic topology (so it suffices to assume that $\cL_2$ is a constant sheaf), by using Proposition \ref{p1230}  and Example \ref{ex1230} in which we take $X_2$ to be a point space.
\eex

\bc\label{cor-tri} Let $X_i\hookrightarrow M$ be  closed complex algebraic (or analytic) subvarieties  of the complex algebraic (or analytic) manifold $M$ of pure dimension $d$ ($i=1,2$). Let $\cS_i$ be Whitney stratifications of $X_i$ ($i=1,2$) which are \index{transversal intersection}
{\it transversal} in $M$, i.e., all strata $S_1\in \cS_1$ are transversal to all strata $S_2\in \cS_2$. This is equivalent to the {\it diagonal embedding} $\Delta\colon  M\to M\times M$ being transversal to the product stratification $\cS_1\times \cS_2$ of $X_1\times X_2$. In particular, $X_1\cap X_2\simeq \left(X_1\times X_2\right)\cap \Delta(M)$ gets an induced Whitney stratification $\cS_1\cap \cS_2$ with strata the connected components of the intersections $S_1\cap S_2$. Then one gets for the induced map $\Delta\colon  X_1\cap X_2\to X_1\times X_2$ the following:
\begin{itemize}
\item[(a)] \ The tensor product $\overset{L}{-\otimes-}\simeq \Delta^*(-\overset{L}{\boxtimes -})$ induces
$$\overset{L}{\otimes}:\:{^pD}^{\leq n}_{\cS_1}(X_1;R)\times {^pD}^{\leq m}_{\cS_2}(X_2;R)\to
{^pD}^{\leq n+m-d}_{\cS_1\cap \cS_2}(X_1\cap X_2;R) \:.$$
\item[(b)] \ Assume $R$ is a field. Then the  tensor product $\overset{L}{-\otimes-}\simeq \Delta^*(-\overset{L}{\boxtimes -})$ induces
$$\overset{L}{\otimes}:\:{^pD}^{\geq n}_{\cS_1}(X_1;R)\times {^pD}^{\geq m}_{\cS_2}(X_2;R)\to
{^pD}^{\geq n+m-d}_{\cS_1\cap \cS_2}(X_1\cap X_2;R) \:.$$
\end{itemize}
In particular, if $R$ is a field,  $\left(-\overset{L}{\otimes}-\right)[-d]$ induces 
$$Perv_{\cS_1}(X_1;R)\times Perv_{\cS_2}(X_2;R) \to Perv_{\cS_1\cap \cS_2}(X_1\cap X_2;R) \:.$$
\ec

\bex In the context of Corollary \ref{cor-tri},  let $X_i$ be in addition  pure dimensional ($i=1,2$), with $R$ a field. Then 
$$\left(IC_{X_1}(\cL_1)\overset{L}{\otimes}IC_{X_2}(\cL_2) \right)[-d]\simeq IC_{X_1\cap X_2}(\cL_1\otimes \cL_2)$$
for a local system $\cL_i$ defined on the open dense  subset $U_i$ of $X_i$ given by the top dimensional stratum ($i=1,2$).
Here $X_1\cap X_2$ is also pure dimensional with $\cL_1\otimes \cL_2$ defined on the open dense subset $U_1\cap U_2$ of
$X_1\cap X_2$.
\eex

The existence of the perverse t-structure on $D^b_{(w)c}(X;R)$ implies the existence of \index{perverse truncation functors} {\it perverse truncation functors} 
${^p\tau}_{\leq 0}, {^p\tau}_{\geq 0}$, which are adjoint to the inclusions 
$${^pD}^{\leq 0}(X;R)\hookrightarrow D^b_{(w)c}(X;R) \hookleftarrow {^pD}^{\geq 0}(X;R)\:.$$
 In particular, for every $k \in \bZ$, there are adjunction maps $\mathscr{F}^\bullet \to {}^{p}\tau_{\ge k}\mathscr{F}^\bullet$ and ${}^{p}\tau_{\le k}\mathscr{F}^\bullet \to \mathscr{F}^\bullet$ (see, e.g., \cite[Proposition 1.3.3]{BBD}).
These perverse truncation functors can be used to associate to any constructible complex $\cFb \in D^b_{(w)c}(X;R)$ its \index{perverse cohomology} {\it perverse cohomology} sheaves defined as:
$${^p\cH}^i(\cFb):= {^p\tau}_{\leq 0} {^p\tau}_{\geq 0} (\cFb[i]) \in Perv(X;R).$$
These are $\cS$-(weakly) constructible for $\cFb$ $\cS$-(weakly) constructible.
It then follows that $\cFb \in {^pD}^{\leq 0}(X;R)$ if and only if ${^p\cH}^i(\cFb)=0$ for all $i>0$. Similarly, $\cFb \in {^pD}^{\geq 0}(X;R)$ if and only if ${^p\cH}^i(\cFb)=0$ for all $i<0$. In particular, $\cFb \in Perv(X;R)$ if and only if ${^p\cH}^i(\cFb)=0$ for all $i\neq 0$ and ${^p\cH}^0(\cFb)=\cFb$
(see, e.g., \cite[Proposition 1.3.7]{BBD}). Perverse cohomology sheaves can be used to calculate the (hyper)cohomology groups of any $\mathscr{F}^\bullet \in D_{(w)c}^b(X)$ via the \index{perverse cohomology spectral sequence} {\it perverse cohomology spectral sequence} 
\be\label{perler} E_2^{i,j} = H^i(X;{^p\cH^j}(\cFb)) \Longrightarrow H^{i+j}(X;\cFb).\ee


\bd A functor $F\colon D_1\to D_2$ of triangulated categories with t-structures is \index{left t-exact functor} {\it left t-exact} if $F(D_1^{\geq 0}) \subseteq D_2^{\geq 0}$, \index{right t-exact functor} {\it right t-exact} if $F(D_1^{\leq 0}) \subseteq D_2^{\leq 0}$, and \index{t-exact functor} {\it t-exact} if $F$ is both left and right t-exact. 
\ed

\bex The inclusion of full subcategories
$$D^b_{\cS-(w)c}(X;R)\subset D^b_{(w)c}(X;R) \quad \text{and} \quad
D^b_{(\cS-)wc}(X;R)\subset D^b_{(\cS-)c}(X;R)$$
are t-exact with respect to the perverse t-structures.
\eex

\br If $F$ is a t-exact functor, it restricts to a functor on the corresponding hearts. More generally, if $F\colon D_1\to D_2$ is a functor of triangulated categories with t-structures, and we let  
$\mathcal{C}_1, \mathcal{C}_2$ be the corresponding hearts with $k_i \colon  \mathcal{C}_i \hookrightarrow D_i$, then 
\begin{displaymath}
 ^pF := \ ^tH^0 \circ F \circ k_1 \colon  \mathcal{C}_1 \rightarrow \mathcal{C}_2
\end{displaymath}
is called the \index{perverse functor} {\it perverse functor associated to $F$}. (Here, if $\tau_{\le 0}$ and $\tau_{\ge 0}$ are the truncation functors on a triangulated category $D$ with heart $\mathcal{C}$, we set $^tH^0:=\tau_{\ge 0}\tau_{\le 0}=\tau_{\le 0}\tau_{\ge 0}\colon  D\to \mathcal{C}$.)
In this paper we work only with the perverse t-structure, so a t-exact functor preserves perverse sheaves. \er

\bex\label{ex410} Let $X$ be a complex  analytic (or algebraic) variety, and let $Z \subseteq X$ be a closed subset. Fix a Whitney stratification $\cS$ of the pair $(X,Z)$, i.e., $Z$ is a union of strata of $\cS$. Then $Z$ and $U:=X \setminus Z$ inherit  Whitney stratifications as well, and if we denote by $i\colon Z \hookrightarrow X$ and $j\colon U \hookrightarrow X$  the stratified  inclusion maps, then the functors $j^*=j^!$, $i^!$, $i^*$, $i_*=i_!$, $j_!$ and $Rj_*$ preserve (weak) constructibility with respect to the above fixed stratifications, with $j^*i_!=0$ and $i^!Rj_*=0$. Moreover, the functors $j_!$, $i^*$ are right t-exact, the functors $j^!=j^*$, $i_*=i_!$ are t-exact, and $Rj_*$, $i^!$ are left t-exact. Similarly for the functors   $j^*=j^!$, $i^!$, $i^*$, $i_*=i_!$, $j_!j^*$ and $Rj_*j^*$ (as well as  $j_!$ and $Rj_*$ in the complex algebraic context), if we do not fix a Whitney stratification.
\eex

\bex Let $f\colon X\to Y$ be a {\it smooth} morphism (i.e., a submersion) of complex algebraic (or analytic) varieties, with constant relative (or fiber) dimension 
$d$. Then  $f^![-d]\simeq f^*[d]$ is t-exact. Similarly, the external tensor product $\overset{L}{\boxtimes}$ is t-exact in each variable in case $R$ is a field.
\eex

\bex\label{finitem} Let $f\colon X \to Y$ be a finite map (i.e., proper, with finite fibers). Then $Rf_*=f_!$ is t-exact (see Example \ref{finite}). If, moreover, $X$ is pure dimensional
with $f$ surjective and generically bijective,  then $Rf_*IC_X \simeq f_* IC_X \simeq IC_Y$. The latter fact applies, in particular, to the case when $X$ is the (algebraic) normalization of $Y$.
\eex

The perverse cohomology sheaf construction provides a way to get perverse sheaves out of any (weakly) constructible complex. Another important method for constructing perverse sheaves, the {\it intermediate extension}, will be discussed below. 

Let $j\colon U \hookrightarrow X$ be the inclusion of an open constructible subset of the complex algebraic (or analytic) variety $X$, with $i\colon Z=X\setminus U \hookrightarrow X$ the closed inclusion. We can also work with a  fixed Whitney stratification $\cS$ of $X$ so that $U$ is an open union of strata of $\cS$. A complex $\cFb \in D^b_{(w)c}(X;R)$ is a (weakly) constructible {\it extension} of $\cGb \in D^b_{(w)c}(U;R)$ if $j^*\cFb \simeq \cGb$. In what follows, we are interested to find perverse extensions of $\cGb \in Perv_{(\cS|_U)}(U;R)$.

\bd\label{interm}
The \index{intermediate extension} {\it intermediate extension} of the perverse sheaf $\cGb \in Perv_{(\cS|_U)}(U;R)$ is the image $j_{!*}\cGb$ in the abelian category $Perv_{(\cS)}(X;R)$ of the morphism $${^pj}_!\cGb :={^p\cH}^0\left(j_!\cGb\right)  \to {^p\cH}^0\left(Rj_*\cGb\right) =:{^pj}_*\cGb.$$
Here, ${^pj}_!\cGb \to {^pj}_*\cGb$ is obtained by applying the functor ${^p\cH}^0$ to the natural morphism $j_!\cGb \to Rj_*\cGb$ in 
$D^b_{(w)c}(X;R)$. In the analytic context we explicitly have to assume that $j_!\cGb$ (and then also $Rj_*\cGb\simeq Rj_*j^*j_!\cGb$) is (weakly) constructible, e.g., $\cGb \in Perv(U;R)$ is (weakly) constructible with respect to $\cS|_U$.
\ed

\bex If $X$ is of pure dimension $n$, and $j\colon U \hookrightarrow X$ is the inclusion of a  smooth open subset whose complement is an algebraic (or analytic) subset of dimension $<\dim X$, then for a local system $\cL$ on $U$ one has that $\cL[n] \in Perv(U;R)$ and
$$IC_X(\cL) \simeq j_{!*}(\cL[n]) \:.$$ 
\eex

\bex\label{cur} Let $X$ be smooth of pure dimension one and let $j\colon U \hookrightarrow X$ be the inclusion of a Zariski open and dense subset. If $\cL$ is a local system on $U$, then:
$$IC_X(\cL)\simeq j_{!*}(\cL[1])\simeq(j_*\cL)[1] \:.$$
\eex

The intermediate extension functor plays an important role in describing the simple objects in the abelian category of ($\cS-)$constructible perverse sheaves on $X$. Indeed, we have the following.
\bp\label{suz} Consider the context of  Definition \ref{interm} above.
\begin{itemize}
\item[(a)] \ The intermediate extension $j_{!*}\cGb$ of $\cGb \in Perv_{(\cS|_U)}(U;R)$ has no non-trivial sub-object and no non-trivial quotient object whose supports are contained in $Z=X\setminus U$. 
\item[(b)] \ Moreover, if $\cGb \in Perv_{(\cS|_U)}(U;R)$ is a simple constructible  object then $j_{!*}\cGb \in Perv_{(\cS)}(X;R)$ is also simple constructible.
\end{itemize}
\ep
Moreover, the following important result holds, see \cite[Theorem 4.3.1]{BBD} (where the assumptions are needed for the use of biduality and the stability  of perverse sheaves under duality):
\bt\label{8.4.14} Let $X$ be a complex algebraic variety, assume that the coefficient ring $R$ is a field, and consider (constructible) 
perverse sheaves $Perv(X;R)\subset D^b_{c}(X;R)$.
\ben
\item[(a)] The category of perverse sheaves $Perv(X;R)$ is {\it artinian} and {\it noetherian}, i.e., every perverse sheaf on $X$ admits an increasing finite filtration with quotients {simple perverse sheaves} \index{simple perverse sheaf}.
\item[(b)] The simple $R$-perverse sheaves on $X$ are the twisted intersection complexes $IC_{\overline V}(\cL)$ (regarded as complexes on $X$ via extension by zero), where $V$ runs through the family of smooth connected constructible subvarieties of $X$, $\cL$ is a simple (i.e., irreducible) $R$-local system of finite rank on $V$, and $\overline V$ is the closure of $V$ in $X$.
\een
\et

\br
A similar result as in Theorem \ref{8.4.14} is also true for a compact analytic variety, or if one works in the algebraic or analytic context with $\cS$-constructible perverse sheaves $Perv_{\cS}(X;R)\subset D^b_{\cS-c}(X;R)$ for
a fixed Whitney stratification $\cS$ of $X$ with only {\it finitely} many strata. In the later case one needs to use in (b) for $V$ only the strata $S\in \cS$.
\er

\bex\label{ex221} Let $R$ be a field of coefficients, and let $X=\bC$, $U=\bC^*$ with open inclusion $j\colon U \hookrightarrow X$, and $Z=\{0\}$ with closed inclusion $i\colon Z \hookrightarrow X$. Let $\cL$ be a local system of finite rank on $U$ with stalk $V$ and monodromy automorphism $h\colon V \to V$. Then $\cL[1]\in Perv(U;R)$ and $  IC_X(\cL) \simeq j_*\cL[1]$. Moreover, it is an instructive exercise to show that the following assertions hold:  
\begin{itemize}
\item[(a)] \ $j_!\cL[1]$ and $Rj_*\cL[1]$ are perverse sheaves on $X$.
\item[(b)] \ There is a short exact sequence of perverse sheaves on $X$:
\be\label{compos1}
0 \lra IC_X(\cL) \lra Rj_*\cL[1] \lra i_!IC_Z(V_h) \lra 0 \ ,
\ee
where $V_h=\coker(h-1)$. In particular, if the local system $\cL$ is simple, then the perverse sheaf $Rj_*\cL[1]$ admits the filtration $$ IC_X(\cL) \subset Rj_*\cL[1]$$ with simple quotients  $IC_Z({V_h})$ and $IC_X(\cL)$.
\item[(c)] \ There is a short exact sequence of perverse sheaves on $X$:
\be\label{compos2}
0 \lra i_!IC_Z(V^h)  \lra j_!\cL[1] \lra IC_X(\cL) \lra 0 \ ,
\ee
where $V^h=\ker(h-1)$. Hence, if the local system $\cL$ is simple, then the perverse sheaf $j_!\cL[1]$ admits the filtration $$ IC_Z(V^h) \subset Rj_*\cL[1]$$ with simple quotients  $IC_X(\cL)$ and $IC_Z({V^h})$.
\end{itemize}
\eex

The next result describes the behavior of the intermediate extension with respect to the dualizing functor (see, e.g., \cite[Proposition 8.4.15]{M}).
\bp\label{sed} Consider the context of  Definition \ref{interm}
with coefficient ring $R$ a field, and let $\cGb \in Perv_{(\cS|_U)}(U;R)$ be ($\cS|_U$-)constructible. Then \index{intermediate extension} 
\be
\cD(j_{!*}\cGb) \simeq j_{!*}(\cD \cGb).
\ee
In particular, $\cD IC_X \simeq IC_X$ for $X$ pure dimensional.
\ep


\subsection{Strongly perverse sheaves. dual t-structure and rectified homological depth}
In this section we only consider ($\cS-$)constructible sheaves so that biduality is available, with a  {\it dual perverse t-structure}
$${^{p^+}D}^{\leq 0}(X;R):=\mathcal{D}\left({^pD}^{\geq 0}(X;R)\right) \quad \text{and} \quad  
{^{p^+}D}^{\geq 0}(X;R):=\mathcal{D}\left({^pD}^{\leq 0}(X;R)\right) $$
 on  $D^b_{(\cS-)c}(X;R)$.
Here $\mathcal{D}$ denotes the Verdier duality functor on $X$. But only for $R$ a {\it Dedekind domain} one can give a more explicit description of this \index{dual perverse t-structure}
{\it dual perverse t-structure}.
If $R$ s a field, the Universal Coefficient Theorem yields that $\cFb \in {^pD}^{\geq 0}(X;R)$ if and only if $\mathcal{D}\cFb \in {^pD}^{\leq 0}(X;R)$. More generally, one has the following result
(see also \cite[Section 3.3]{BBD}):
\bp\label{pid}
Assume that the ring $R$ is a Dedekind domain (e.g., a field or a PID). If $\cFb \in D^b_c(X;R)$ is constructible with respect to a Whitney stratification $\cS$ of $X$, then $\mathcal{D}\cFb \in {^pD}^{\leq 0}(X;R)$, or equivalently $\cFb \in {^{p^+}D}^{\geq 0}(X;R)$,
if and only if the following two conditions are satisfied:
\begin{itemize}
\item[(i)] $\cFb \in {^pD}^{\geq 0}(X;R)$;
\item[(ii)] $\,$ for any stratum $S \in \cS$ and any $x \in S$, the costalk cohomology $H^{\dim S}(i_x^!\cFb)$ is torsion-free.
\end{itemize}
\ep
\begin{proof} Let $S \in \cS$ and $x \in S$, with inclusion $i_x\colon \{x\} \hookrightarrow X$.
Properties of the dualizing functor and the Universal Coefficient Theorem yield:
$$H^j(i_x^*\mathcal{D}\cFb)\simeq H^j(\mathcal{D}i_x^!\cFb)
\simeq {\rm Hom}(H^{-j}(i_x^!\cFb), R) \oplus {\rm Ext}(H^{-j+1}(i_x^!\cFb), R).$$
The desired equivalence can now be checked easily.
In particular, $${^{p}D}^{\geq -1}(X;R)\subset {^{p^+}D}^{\geq 0}(X;R)\subset {^{p}D}^{\geq 0}(X;R).$$
\end{proof}

\bd\label{stp} Assume that $R$ is a Dedekind domain.
We say that $\cFb \in D^b_c(X;R)$ is  \index{strongly perverse sheaf} {\it strongly perverse} if $\cFb \in {^pD}^{\leq 0}(X;R) $ and  $\mathcal{D}\cFb \in {^pD}^{\leq 0}(X;R)$. 
Equivalently, $$\cFb \in {^pD}^{\leq 0}(X;R)\cap {^{p^+}D}^{\geq 0}(X;R).$$
\ed

If $R$ is a field, the notions of perverse sheaf and strongly perverse sheaf coincide. 
If $R$ is a  Dedekind domain, then $\cFb$ is strongly perverse (with respect to $\cS$) if and only if $\cFb$ is perverse and property (ii) above holds (i.e., costalks of $\cFb$  in the lowest possible degree are torsion-free on each stratum).

\bex
Let $X$ be a pure dimensional complex algebraic or analytic variety, and let $R$ be a Dedekind domain.
Then $IC_X(\cL)$ is strongly perverse if the generically defined local system $\cL$ has finitely generated and torsion-free stalks
(since then condition (ii) only needs to be checked for the top dimensional strata  $S$, with $H^{\dim S}(i_x^!\cL[\dim X])\simeq \cL_x$ for $x$ in such a stratum $S$).
\eex

Strongly perverse sheaves are related to the notion of {\it rectified homological depth}, which we now define. 

Let $X$ be a complex algebraic or analytic variety. Following \cite[Definition 6.0.4]{Sc}, we make the following.
\bd\label{defrhd} The \index{rectified homological depth} {\it rectified homological depth} $\r(X,R)$ of $X$ with respect to the commutative base ring $R$ is $\geq d$ (for some $d \in \bZ$) if 
\be\label{rhdd}
\mathcal{D}({R}_X [d]) \in {^pD}^{\leq 0}(X;R)\:,
\ee
or, equivalently,
\be
R_X \in {^{p^+}D}^{\geq d}(X;R)\:.
\ee
The {\it rectified homological depth} $\r(X,R)$ of $X\neq \emptyset$ with respect to the commutative base ring $R$ is the maximum of such integers $d$.
\ed
As indicated in \cite[p. 387]{Sc}, the above definition agrees with the notion of {\it rectified homological depth} introduced by Hamm and L\^e \cite{HL} (following an earlier definition of Grothendieck, together with the corresponding homotopical notion) in more geometric terms. In the following, let $\dim_{-} X$ be the minimum of the dimension of the irreducible components of $X$. So for $X\neq \emptyset $ smooth, this is the minimum of the dimension of the connected components of $X$. Similarly, $\dim_{-} X=\dim_{-} X_{reg}$,
for $X_{reg}$ the open dense regular part of $X$.
\begin{example}\label{ex27}
\begin{itemize}
\item[(a)] \ One always has $\r(X,R) \leq \dim_{-} X$, and $\r(X,R) = \dim_{-} X$ if $X$ is smooth and nonempty (since then 
$\bD_X^{\bullet}\simeq R_X[2\dim X]$ with $\dim X$ viewed as a locally constant function). Moreover, $\r(X,R) =\dim X$ implies that $X$ is pure dimensional
(by looking at the regular part).
\item[(b)] \ If $X$ is a pure-dimensional local complete intersection, \index{local complete intersection}
then $\r(X,R)=\dim X$ (see Example \ref{purity} and, e.g.,  \cite[Example 6.0.11]{Sc}).
\end{itemize}
\end{example}
In view of Example \ref{ct1}(a) and Definition \ref{stp}, one has the following equivalence (see also \cite[(6.14)]{Sc}):
\bp\label{p28} Let  $R$ be a Dedekind domain.
For any nonempty complex algebraic or analytic variety $X$ one has:
\begin{center}
$\r(X,R)=\dim X$ $\iff$  $X$ is pure-dimensional and $R_X [\dim X]$ is strongly perverse.
\end{center}
\ep
As a consequence, Proposition \ref{pid} yields the following.
\bc\label{cfield}
Let $X$ be a nonempty pure-dimensional complex algebraic or analytic variety with a Whitney stratification $\cS$.
\begin{itemize}
\item[(a)] \ If $R$ is a field, then:
\begin{center}
$\r(X,R)=\dim X$ $\iff$ $R_X [\dim X]$ is perverse.
\end{center}
\item[(b)] $\,$ If $R$ is a Dedekind domain (e.g., a PID), then $\r(X,R)=\dim X$ if and only if the following two conditions are satisfied:
\begin{itemize}
\item[(i)] $R_X [\dim X]$ is perverse.
\item[(ii)] $\,$ for any stratum $S \in \cS$ and any $x \in S$ with $i_x\colon \{x\} \hookrightarrow X$, the costalk cohomology $H^{\dim S}(i_x^!R_X[\dim X])$ is torsion-free.
\end{itemize}
\end{itemize}
\ec

Let us finish this section with citing the following result from \cite[Section 3.3]{BBD}.

\bp\label{pid2}
Assume that the ring $R$ is a Dedekind domain (e.g., a field or a PID). If $\cFb \in D^b_c(X;R)$ is constructible with respect to a Whitney stratification $\cS$ of $X$, then $\mathcal{D}\cFb \in {^pD}^{\geq 0}(X;R)$, or equivalently $\cFb \in {^{p^+}D}^{\leq 0}(X;R)$,
if and only if the following two conditions are satisfied:
\begin{itemize}
\item[(i)] $\cFb \in {^pD}^{\leq 1}(X;R)$;
\item[(ii)] $\,$ for any stratum $S \in \cS$ and any $x \in S$, the stalk cohomology $H^{-\dim S+1}(i_x^*\cFb)$ is a torsion module.
\end{itemize}
In particular, ${^{p}D}^{\leq 0}(X;R)\subset {^{p^+}D}^{\leq 0}(X;R)\subset {^{p}D}^{\leq 1}(X;R)$.
\ep

Assume $R$ is a Dedekind domain. Let us denote by $D^b_{(\cS-)tc}(X;R)$ the full triangulated subcategory of $(\cS$-)constructible sheaf complexes $\cFb \in D^b_{(\cS-)c}(X;R)$, with all stalk cohomology modules  $H^k(i_x^*\cFb)$ a torsion module for all $x\in X, k\in \bZ$. These are also stable by the usual truncation functors.
Similarly, in the context of Example \ref{ex410}, all considered functors preserve the subcategory $D^b_{\cS-tc}(-;R)$ of these  torsion $\cS$-constructible sheaf complexes. For $i_!=i_*$ and $j^*=j^!$ this is trivial, and for $i^!$ it follows from the corresponding property for $Rj_*$ and the standard distinguished triangle
$$i_*i^! \lra id \lra Rj_*j^* \lra \:.$$
 But the stalk cohomology modules of the
open push-forward $Rj_*$ can be expressed in terms of the compact link  \cite[Remark 4.4.2]{Sc}, and by the K\"{u}nneth formula they vanish after tensoring with the quotient field $Q(R)$ of $R$. Then the usual perverse t-structure and its truncation functors restrict to the perverse t-structure
\be \left(  {^pD}^{\leq 0}_{\cS-t}(X;R),  {^pD}^{\geq 0}_{\cS-t}(X;R)\right) \quad \text{on $D^b_{\cS-tc}(X;R)$,}
\ee
and the corresponding heart of {\it $\cS$-constructible torsion perverse sheaves} \index{torsion perverse sheaves}
$$Perv_{\cS-t}(X;R)\:.$$

\bex Let $X$ be a pure dimensional complex algebraic or analytic variety with a Whitney stratification $\cS$, and let $R$ be a Dedekind domain.
Then $IC_X(\cL)$ is a  $\cS$-constructible torsion perverse sheaf if the  local system $\cL$ defined on the open top dimensional strata  has only finitely generated torsion stalks.
\eex

\bc Assume $R$ is a Dedekind domain. Then one gets for the \index{shifted Verdier duality}
{\it shifted Verdier duality functor} $(\cD(-))[1]$:
\begin{itemize}
\item[(a)] \ $(\cD(-))[1]$ maps ${^pD}^{\leq 0}_{\cS-t}(X;R)$ into ${^pD}^{\geq 0}_{\cS-t}(X;R)$.
\item[(b)] \ $(\cD(-))[1]$ maps ${^pD}^{\geq 0}_{\cS-t}(X;R)$ into ${^pD}^{\leq 0}_{\cS-t}(X;R)$.
\end{itemize}
In particular, the category $Perv_{\cS-t}(X;R)$ of  $\cS$-constructible torsion perverse sheaves  is stable under the 
{\it shifted Verdier duality functor} $(\cD(-))[1]$.
\ec

\bex Assume $R$ is a Dedekind domain. Then on  a point space $X=\{pt\}$ the {\it shifted  duality functor} $(\cD(-))[1]$ is given by
$$ R\cH om^{\bullet}(-,R)[1] \simeq \cH om^{\bullet}(-,Q(R)/R):\: D^b_{tc}(\{pt\};R)\to D^b_{tc}(\{pt\};R)\:,$$
preserving finitely generated torsion modules (viewed as complexes concentrated in degree zero). Here $Q(R)$ is the quotient field of $R$ so that the short exact sequence $$0\lra R \lra Q(R) \lra Q(R)/R\lra 0$$ gives the injective resolution $[Q(R) \lra Q(R)/R]$ of $R$  used for calculating
$R\cH om^{\bullet}(-,R)$.
\eex

In fact this example is also used on a general $X$ for showing that  the  Verdier duality functor $\cD$ preserves $D^b_{(\cS-)tc}(X;R)$, together with
$i_x^*\cD\simeq \cD i_x^!$ for all $x\in X$.
For applications and a discussion of related results see also \cite{CS, GS,  M0, M1, SW}.


\section{Stratified Morse theory for constructible sheaves}\label{smt}

In this section, we explain  the basic results from \cite{Sc} about {\it stratified Morse theory}
in the framework of (weakly)  constructible sheaves in the {\it complex context}, continuing  and extending the recent survey of Goresky \cite {Gor} in this handbook series (as well as the survey of Massey \cite{Ma1}).
We follow the notions of the geometric {\it stratified Morse theory} of
Goresky-Mac\-Pherson so that one can easily compare our results with 
those of \cite{GMs}. Moreover, we also explain some relations to the general {\it micro-local sheaf theory}
of Kashiwara-Schapira \cite{KS}. But unlike most of these  references, we do not need a {\it global embedding} into an ambient complex manifold
(except for the index theorems for characteristic cycles of constructible sheaves and the comparison with the notion of micro-support from \cite{KS}).

\subsection{Morse functions, local and normal Morse data}

We work with a fixed complex algebraic (or analytic) variety $X$ with a given Whitney stratification $\cS$.

\bd
A function $f\colon  X\to \bR$  is a {$C^{k}$-differentiable function} (with $2\leq k\leq \infty$)
 if, for any $x\in X$, there exists a local embedding $(X,x)\hookrightarrow 
(\bC^{n},x)$ such that $f$ is induced by restriction of a
$C^{k}$-function germ $\hat{f}\colon  (\bC^{n},x)\to (\bR,f(x))$.
\ed

Of course, a complex algebraic (or analytic) function or morphism $h\colon  X\to \bC$ is by definition locally given by  restriction of a
complex algebraic (or analytic) function germ $\hat{h}\colon  (\bC^{n},x)\to (\bC,h(x))$. 

\bd A point $x\in X$ is called a \index{stratified critical point} {\it stratified critical point} of $f\colon  X\to \bR$, resp., $h\colon  X\to \bC$ (with respect to $\cS$) if
$f\colon  S\to \bR$, resp., $h\colon  S\to \bC$ is not a submersion at the point $x$ belonging to the stratum $S\in \cS$.
Let $\Sing(f|_S)$, resp., $\Sing(h|_S)$ be the set of critical points of $f$, resp., $h$ in the stratum $S\in\cS$, with
$$\Sing_{\cS}(f):=\bigcup_{S\in \cS} \Sing(f|_S) \quad \text{resp.,} \quad \Sing_{\cS}(h):=\bigcup_{S\in \cS} \Sing(h|_S) \:.$$
Then $\Sing_{\cS}(f)$ (resp.,  $\Sing_{\cS}(h)$) is a closed (complex algebraic or analytic) subset of $X$ by the {\it Whitney a-condition}.
\ed

The following result shows that for $\cS$-weakly constructible sheaf complexes the local change of cohomology with respect to $f$ is located at 
the critical points $x\in \Sing_{\cS}(f)$ (see, e.g., \cite[Example 5.1.1]{Sc}).

\bl \label{loc.const} 
Let $f\colon  X\to \,]c,d[\,\subset \bR$ be a {\it proper} differentiable function, and assume $\cFb \in D^b_{\cS-wc}(X;R)$ is a $\cS$-weakly constructible complex.
\begin{enumerate}
\item If $f$ is a stratified submersion at a point $x\in \{f=e\}$ (with $e \in \,]c,d[\,$, and here $f$ does not need to be proper), then 
$$\bigl(R\Gamma_{\{f\geq e\}}(\cFb)\bigr)_{x} \simeq 0 \simeq \bigl(R\Gamma_{\{f\leq e\}}(\cFb)\bigr)_{x} \:.$$
\item If $f$ is a stratified submersion at all points of $\{f=e\}$ (for $e\in \,]c,d[\,$), then 
$$R\Gamma(\{f=e\},R\Gamma_{\{f\geq e\}}(\cFb))\simeq 0
\simeq R\Gamma(\{f=e\},R\Gamma_{\{f\leq e\}}(\cFb)) \:.$$
\item If $f$ is a  stratified submersion at all points of $\{a<f<b\}$ (with $]a,b[\,\subset
\,]c,d[\,$), then all cohomology sheaves of $(Rf_{*}{\cFb})|_{\,]a,b[\,}$ are locally constant.
\end{enumerate}
\el

These results can be proved by induction on the dimension of $X$ {\it without using the ``first isotopy Lemma of Thom''}. Property 1. above just means that a Whitney b-regular stratification satisfies the {\it ``local stratified acyclicity''} property form \cite[Definition 4.0.3]{Sc}. In the embedded context of a global closed embedding $k\colon  X\hookrightarrow M$ into a complex manifold it also implies the following estimate \cite[Corollary 4.0.3]{Sc} of the 
\index{micro-support} {\it micro-support} of $Rk_*\cFb$ in the sense of Kashiwara-Schapira \cite[Definition 5.1.2]{KS} (with $Rk_*=k_*=k_!$ the extension by zero):
\be\label{emusupp}
\mu supp(Rk_*\cFb)\subset \bigcup_{S\in \cS}\: T^*_SM \hookrightarrow T^*M \quad \text{for any $\cFb\in D^b_{\cS-wc}(X;R)$.}
\ee
Here, $\mu supp(\cGb)\subset T^*M$ is by definition the complement of the largest open subset $U\subset T^*M$ such that 
$R\Gamma_{\{f\geq f(x)\}}(\cGb)_{x} \simeq 0$ for any differentiable function germ $f\colon  (M,x)\to (\bR,f(x))$ with $df_x\in U$ (see \cite[Definition 5.1.2]{KS}).
We also use the notation
$T^{*}_{S}M$ for the {\it conormal bundle} 
of the stratum $S$ in $M$,
i.e., for the kernel of the natural vector bundle epimorphism $T^{*}M|_{S} \to T^{*}S$
dual to the inclusion $TS\hookrightarrow TM|_{S}$. Let us denote by $T^*_{\cS}M$ this union of conormal spaces of strata $S\in \cS$. Then the Whitney a-condition for $\cS$ is equivalent to the fact that $T^*_{\cS}M\hookrightarrow T^*M$ is closed in $T^*M$. Moreover it is {\it conic} in the sense that it is invariant under the natural $\bC^*$-action on the cotangent bundle. The inclusion of the micro-support in (\ref{emusupp}) implies that 
 $\mu supp(Rk_*\cFb)\subset T^*M$ is a closed complex algebraic (or analytic) conic Lagrangian subset of $T^*M$,
and (\ref{emusupp})  is even equivalent to
$\cS$-weak constructibility (see, e.g., \cite[Lemma 4.1]{Sc}).
\bl 
Let $\cGb\in D^b(M;R)$ be given with $\mu supp(\cGb)\subset T^*_{\cS}M$. Then  $supp(\cGb)\subset X$ and $\cGb|_X$ is $\cS$-weakly constructible.
\el

\bex
Assume $X=M$ with $\cGb\in D^b(M;R)$ be given. Then  $\mu supp(\cGb)\subset T^*_{M}M$ the zero section of $T^*M$ if and only if all cohomology sheaves 
$\cH^i(\cGb)$ are locally constant on $M$.
\eex

This  Example can  be  proved directly with the {\it ``non-characteristic deformation Lemma''} of Kashiwara (see, e.g., \cite[Proposition 4.1.1]{Sc}).
Moreover, it can be used as a substitute for the {\it ``first isotopy Lemma of Thom''}.\\

If one does not want to fix a complex Whitney stratification $\cS$, then one gets the following micro-local characterization of weak constructibility
(see, e.g., \cite[Theorem 4.0.1]{Sc} and \cite[Theorem 8.5.5]{KS}).
\bt
Let $M$ be a complex algebraic (or analytic) manifold. Then $\cGb\in D^b(M;R)$ is a weakly constructible complex $\cGb\in D^b_{wc}(M;R)$ in the complex algebraic (or analytic) sense if and only if $\mu supp(\cGb)\subset T^*M$ is a closed complex algebraic (or analytic) conic Lagrangian subset of $T^*M$.
\et

Let us now come back to the general non-embedded context with $\cS$ a fixed Whitney stratification of $X$. Then Lemma \ref{loc.const} implies
(see, e.g., \cite[Corollary 5.1.1]{Sc}).

\bc \label{cor:coarse}
Let $f\colon  X\to \,]c,d[\,\subset \bR$ be a {\it proper} differentiable function, and assume $\cFb \in D^b_{\cS-wc}(X;R)$ is a $\cS$-weakly constructible complex.
If $f$ is a stratified submersion at all points of $\{a\leq f\leq b\}\cap \{f\neq e\}$ (with $[a,b]\subset \,]c,d[\,$ and $e\in \,]a,b[\,$), then one has 
distinguished triangles 
\be \label{eq:triangle1}
R\Gamma(\{f=e\},R\Gamma_{\{f\geq e\}}(\cFb)) \to 
R\Gamma(\{f\leq b\},\cFb) \to R\Gamma(\{f\leq a\},\cFb) \to \:,
\ee
and 
\be\label{eq:triangle2} 
R\Gamma_{c}(\{f<a\},\cFb) \to R\Gamma_{c}(\{f<b\},\cFb) \to
R\Gamma(\{f=e\},R\Gamma_{\{f\leq e\}}(\cFb)) \to \:.
\ee
\ec

Here, $R\Gamma(\{f=e\},R\Gamma_{\{f\geq e\}}(\cFb))$ is the cohomological
counterpart of the \index{coarse Morse datum} {\it ``coarse Morse datum''} as defined in  \cite[p.62, Definition 3.4]{GMs}.
In particular, 
$$R\Gamma(\{f=e\},R\Gamma_{\{f\leq e\}}(\cFb)\simeq
R\Gamma(\{-f=-e\},R\Gamma_{\{-f\geq -e\}}(\cFb))$$
is the ``coarse Morse datum'' of $-f$. So  this ``coarse Morse datum'' of $-f$ corresponds to the
``relative cohomology with compact support'' in the triangle~(\ref{eq:triangle2}). 
But the cohomology with compact support is in some sense
``dual'' to the cohomology with closed support. So this reflects the ``duality''
observed in \cite[p.27, 2.7]{GMs}. In the context of stratified Morse theory
for constructible sheaves, this is a special form of {\it Poincar\'{e}-Verdier duality}
(see Equation (\ref{LMorse-dual})).\\

The next step is to consider a differentiable function $f$ with {\it isolated} stratified critical 
points in $\{f=e\}$, and to localize the ``coarse Morse datum'' at these
critical points. In terms of sheaf theory, this is quite easy (since $R\Gamma_{\{f\geq e\}}(\cFb)|_{\{f=e\}}$ is supported on the stratified critical points).

\bl \label{lem:locmorse}
Let $f\colon  X\to \,]c,d[\,\subset \bR$ be a {\it proper} differentiable function, and assume $\cFb \in D^b_{\cS-wc}(X;R)$ is a $\cS$-weakly constructible complex.
Assume $f$ is a stratified submersion at all points of $\{a\leq f\leq b\}$ except for finitely many $x_{i}\in \{f=e\}$
(with $[a,b]\subset \,]c,d[\,$ and $e\in \,]a,b[\,$). Then
\begin{equation} \label{eq:locmorse}
R\Gamma(\{f=e\},R\Gamma_{\{f\geq e\}}\cFb)) \simeq
\bigoplus_{i} \:\left(R\Gamma_{\{f\geq e\}}(\cFb)\right)_{x_{i}}\: .
\end{equation} 
\el

The stalk complex 
\begin{equation} \label{eq:LMD}
LMD(\cFb,f,x):=
\bigl(R\Gamma_{\{f\geq f(x)\}}(\cFb)\bigr)_{x} 
\end{equation}
is the sheaf theoretic counterpart of the {\it``local Morse datum''} \index{local Morse datum}
of \cite[p.63, Definition 3.5.2]{GMs}.  The ``duality'' between $f$ and $-f$ discussed above is closely related to {\it Verdier duality for the local Morse data} in the following form of \cite[Equation (5.54) on p.314]{Sc}:
\be\label{LMorse-dual}
\cD(LMD(\cFb,-f,x))\simeq LMD(\cD(\cFb),f,x) \quad \text{for $\cFb\in D^b_{\cS-wc}(X;R)$.}
\ee 

Let us now recall two notions from \cite{GMs}. 

\bd \label{def:Morse,nslice}
Fix a point $x$ in the stratum $S\in \cS$.
\begin{enumerate}
\item $x$ is called a \index{stratified Morse critical point}  {\it stratified Morse critical point} 
of the differentiable function germ $f\colon  (X,x)\to \bR$ (with respect to
$\cS$), if $f|_S$ has at $x$ a Morse critical point in the 
classical sense (i.e., its Hessian at $x$ is non-degenerate), and $f=\hat{f}|_X$ is induced 
in some local embedding $(X,x)\hookrightarrow  (\bC^{n},x)$ by a differentiable function germ $\hat{f}\colon  (\bC^{n},x)\to (\bR,f(x))$
such that the
covector $d\hat{f}_{x}$ is {\it non-degenerate} 
in the sense of \cite[p.44, Definition  1.8]{GMs},
i.e., it does not vanish on any generalized tangent space
$$\tau := \lim_{x_{n}\to x} \;T_{x_{n}}S' \quad \text{at $x$,}$$
with $x_{n}$ a sequence in a stratum $S'\in \cS$ with $S\subset \partial S'$. 
\item A \index{normal slice} {\it normal slice} $N$ at $x\in S\hookrightarrow X\hookrightarrow \bC^n$
is a closed complex submanifold germ $N$ of $(\bC^{n},x)$,
with $N\cap S=\{x\}$ such that $N$ intersects $S$ {\it transversally} at $x$ (so that $\codim N = \dim S$).
\end{enumerate}
\ed

Note that a {\it normal slice} $N$ intersects all strata transversally near $x$ by 
{\it Whitney a-regularity}. Therefore $N\cap X$ gets an induced Whitney
stratification $\cS|_N$ near $x$ with strata the connected components of intersections $S'\cap N\neq \emptyset $ ($S'\in \cS$) and $\{x\}=S\cap N$ now a point stratum.

Similarly, a {\it stratified Morse critical point} of $f$ is an
isolated stratified critical point of $f$. Finally, if $f=\hat{f}|_X$ is induced 
in some local embedding $(X,x)\hookrightarrow  (\bC^{n},x)$ by a differentiable function germ $\hat{f}\colon  (\bC^{n},x)\to (\bR,f(x))$
such that the
covector $d\hat{f}_{x}$ is {\it non-degenerate}, then one can find such an $\hat{f}$ in any local embedding  $(X,x)\hookrightarrow  (\bC^{m},x)$.

\bex Let $S$ be an open stratum of $X$. Then $x\in S$ is a stratified Morse critical point of $f$ if and only if $f|_S$ has at $x$ a classical Morse critical point.
If $S=\{x\}$ is a point stratum, then $x$ is by definition a stratified critical point for any differentiable function $f\colon  X\to \bR$. It is 
stratified Morse critical point  if and only if  $f=\hat{f}|_X$ is induced 
in some local embedding $(X,x)\hookrightarrow  (\bC^{n},x)$ by a differentiable function germ $\hat{f}\colon  (\bC^{n},x)\to (\bR,f(x))$
such that the covector $d\hat{f}_{x}$ is non-degenerate.
\eex

\bex
Let $X=M$ be a complex manifold, with $M'\hookrightarrow M$ a closed complex submanifold of positive codimension and the Whitney stratification $\cS$ given by the connected components of $M'$ and $M\backslash M'$. Then a differentiable function $f\colon  M\to \bR$ has a stratified Morse critical point at $x\in M'$  if and only if $x$ is not a critical point of $f\colon  M\to \bR$ (this is equivalent to $df_x$ is nondegenrate) and $x\in M'$ is a classical Morse critical point for $f\colon  M'\to \bR$.
\eex

Let us now state the first main result from stratified Morse theory for constructible sheaves (see, 
e.g., \cite[Theorem 5.0.1]{Sc}).

\bt \label{thm:SMTB} Let $X$ be a complex algebraic (or analytic) variety with a Whitney stratification $\cS$, and let  $\cFb\in D^b_{\cS-wc}(X;R)$ be a given
$\cS$-weakly constructible complex.
Let $f\colon  (X,x)\to \bR$ be a differentiable function germ with a {\it stratified Morse critical
point} at $x\in S$ for some stratum $S\in \cS$. Take a {\it normal slice} $N\subset \bC^n$ at $x$ 
in some local embedding $(X,x)\hookrightarrow  (\bC^{n},x)$.
\begin{enumerate} 
\item Then one has an isomorphism
$$\bigl(R\Gamma_{\{f\geq f(x)\}}(\cFb)\bigr)_{x} \simeq
\bigl(R\Gamma_{\{f|_{N\cap X}\geq f(x)\}}(\cFb|_{N\cap X})\bigr)_{x} [-\tau] \:,$$
with $\tau$ the \index{Morse index} {\it Morse index} of $f|_S$ (i.e., its Hessian at $x$ has exactly
$\tau$ negative eigenvalues).
\item The isomorphism class of 
$$\bigl(R\Gamma_{\{f|_{N\cap X}\geq f(x)\}}(\cFb|_{N\cap X})\bigr)_{x}$$
is {\it independent} of the choice of such a local embedding and normal slice $N$.  Moreover,
this isomorphism class {\it only depends} on the stratum $S\in \cS$, but not on the point $x\in S$ or the function germ $f=\hat{f}|_X$ with
$d\hat{f}_x$ non-degenerate in such a local embedding.
\end{enumerate}
\et

The isomorphism class of \index{normal Morse data}
\begin{equation} \label{eq:NMD} 
NMD(\cFb,S):=
\bigl(R\Gamma_{\{f|_{N\cap X}\geq f(x)\}}(\cFb|_{N\cap X})\bigr)_{x}
\end{equation}
is the sheaf theoretic counterpart of the   {\it ``normal Morse data''} of 
\cite[p.65, Definition 3.6.1]{GMs}. In the (locally) embedded context $k\colon  X\hookrightarrow M$, $NMD(\cFb,S)$ is also isomorphic to the 
{\it ``micro-local type of $Rk_*\cFb$ at the non-degenerate covector $\omega=d\hat{f}_x$''} in the sense of Kashiwara-Schapira \cite[Proposition 6.6.1(ii)]{KS}
(as shown in \cite[Equation (5.52) on p.311]{Sc}).
Part 1. of Theorem \ref{thm:SMTB}  corresponds to
\cite[p.8, Theorem SMT Part B; p.65, Main Theorem 3.7]{GMs}, and part 2. to
\cite[p.93, Theorem 7.5.1]{GMs} (compare also with \cite[p.223, Proposition 6.A.1]{GMs}).\\

Let us consider again the embedded context with $k\colon  X\hookrightarrow M$ a closed embedding into a complex algebraic (or analytic) manifold,
and $\cFb\in D^b_{\cS-wc}(X;R)$.
By using the estimate  (\ref{emusupp}) and the {\it involutivity} of $\mu supp(Rk_{*}\cFb)$
 \cite[Theorem 6.5.4]{KS} one can show that $\mu supp(Rk_{*}\cFb)$ is a union of
closures of conormal bundles $T^*_SM$ of some strata $S\in \cS$. Then the explicit characterization  of the micro-suport $\mu supp(Rk_{*}\cFb)$
is given as follows (see, e.g., \cite[Proposition 5.0.1]{Sc}). \index{micro-support}

\bp \label{pmusupp}
Let $\cFb\in D^b_{\cS-wc}(X;R)$ be given. Then
$$ 
\overline{T^*_SM} \subset \mu supp(Rk_{*}\cFb) \iff NMD(\cFb,S) \not\simeq 0  \:.$$
\ep

For the applications of stratified Morse theory for constructible sheaves, it is important to get information about this {\it normal Morse datum}
$NMD(\cFb,S)$ for a stratum $S\in \cS$. In our complex context this can be obtained in the following way. Since it is a local study,
we can assume that for a given point $x\in S$ we have (locally) a closed embedding  $(X,x)\hookrightarrow (\bC^n,x=0)$.

We already know that the isomorphism class of the ``normal Morse 
datum''  $NMD(\cFb,S)$ only depends  on the non-degenerate covector $\omega=d\hat{f}_{x}\in T^*_S\bC^n|_x$. So one can use 
a {\it holomorphic} function germ $g\colon  (\bC^n,x)\to (\bC,0)$ with $d\hat{f}_{x}=
d(Re(g))_{x}$. Also choose a local {\it normal slice} $N$ at $x\in S\hookrightarrow X\hookrightarrow \bC^n$, with $\{x\}=S\cap N$ a point stratum of the induced Whitney stratification $\cS|_N$ of $X\cap N$.
Then one gets (see   Remark \ref{rNMD} and also, e.g., \cite[Proposition 5.0.3]{Sc}).

\bp\label{prop:NMDtriangleC}
Let $k_{x}\colon  \{x\}\to N\cap X$ be the inclusion of this point stratum, and $\cFb\in D^b_{\cS-wc}(X;R)$ be given. Then there exists two 
distinguished triangles
\be\label{clink-stalk}\begin{CD}
R\Gamma(l_{X},\cFb)[-1] @>>> 
NMD(\cFb,S) @>>>
k_{x}^{*}(\cFb|_{N\cap X}) @>>> \:,
\end{CD} \ee
and
\be \label{clink-costalk}\begin{CD}
k_{x}^{!}(\cFb|_{N\cap X})  @>>> 
NMD(\cFb,S) @>>>
R\Gamma(l_{X},\partial l_{X},\cFb)[-1]  @>>> \:.
\end{CD} \ee
Moreover, the isomorphism classes of $R\Gamma(l_{X},\cFb)$ and
$R\Gamma(l_{X},\partial l_{X},\cFb)$ only depend
on the stratum $S$ of $X$ (but on no other choice like the local embedding $(X,x)\hookrightarrow (\bC^n,x=0)$, or the choice of $N$ or $g$ as above).
\ep

Here we use the following notation (with $0<|w| \ll \delta \ll 1$ and $r(z):=\sum_{1\leq i \leq n} z_i\bar{z}_i$ the square of the local distance to $x=0\in \bC^n$)
\be 
(l_{X},\partial l_{X}) := (X\cap N\cap\{r\leq\delta, g=w\},X\cap N\cap\{r=\delta, g=w\}) 
\ee
for the \index{complex link} {\it complex link} of $X$ with respect to $S\in \cS$
(see \cite[p.161, Definition.2.2]{GMs}). In other words, this is the local \index{Milnor fiber} {\it Milnor fiber (with its boundary)} (see, e.g., \cite{Le4} or \cite[Example 1.1.2]{Sc}) of the holomorphic function germ
$g\colon  (X\cap N,x)\to (\bC,0)$, which has an isolated stratified critical point in $x$ with respect to $\cS|_N$ (since $dg_x$ is a non-degenerate covector).

\bex\label{top-NMD} \index{open stratum}
Assume $S\in\cS$ is an {\it open} stratum in $X$. Then for $x\in S$ and a normal slice $N$ at $x$ one gets $X\cap N =\{x\}$, so that the corresponding complex 
link is empty, $l_{X}=\emptyset$, with
$$NMD(\cFb,S)\simeq k_{x}^{*}(\cFb|_{N\cap X}) \simeq i_x^*\cFb$$
just the stalk of $\cFb\in D^b_{\cS-wc}(X;R)$ at $x$.
\eex

\bex\label{ex:NMDtriangleC}
We state here cohomological counterparts of some results of \cite{GMs}.
\begin{enumerate}
\item If we take $\cFb=\bZ_{X}$, then we get by the first distinguished
triangle  (\ref{clink-stalk}) for the cohomology of the ``normal Morse datum'' $NMD(\cFb,S)$ for $x\in S$:
$$\bigl(R^{k}\Gamma_{\{Re(g)|_{N\cap X}\geq 0\}}(\cFb|_{N\cap X})\bigr)_{x} \simeq
\begin{cases} 
\bZ &\text{for $k=0$ and $l_{X}=\emptyset$,}\\
0 &\text{for $k \neq 0$ and $l_{X}=\emptyset$,}\\
\widetilde{H}^{k-1}(l_{X}; \bZ) &\text{for $l_{X} \neq \emptyset$.}
\end{cases} $$
This corresponds to the ``first part of the fundamental theorem of complex stratified
Morse theory'' \cite[p.16, Theorem CSMT Part A, p.166 Corollary 2.4.1]{GMs}.
\item Consider now the case $\cFb=Rj_{*}\bZ_{U}$ for $j\colon  U\hookrightarrow X$
the inclusion of an open subspace, which is a union of strata. Assume 
$x\not \in U$. 
Then one gets by the second distinguished
triangle  (\ref{clink-costalk})   for the cohomology of the ``normal Morse datum'' $NMD(\cFb,S)$ for $x\in S$:
$$\bigl(R^{k}\Gamma_{\{Re(g)|_{N\cap X}\geq 0\}}(\cFb|_{N\cap X})\bigr)_{x} \simeq
H^{k-1}(l_{X},\partial l_{X}; Rj_{*}\bZ_{U}) \:.$$
Moreover, by the base change isomorphism (\ref{bc-open})   one gets
$$H^{k-1}(l_{X},\partial l_{X}; Rj_{*}\bZ_{U}) \simeq
H^{k-1}(l_{U},\partial l_{U}; \bZ) \:,$$
with $(l_{U},\partial l_{U}) :=
(l_{X},\partial l_{X}) \cap U$
the {\it complex link of $U$}. 
This corresponds to the ``second part of the fundamental theorem of complex stratified
Morse theory'' \cite[p.18, Theorem CSMT Part B, p.169, Corollary 2.6.1]{GMs}.
\end{enumerate} \eex

Proposition~\ref{prop:NMDtriangleC} and the corresponding distinguished triangles for $dg_x\in T^*_S\bC^n$  a non-degenerate covector are related to the theory \index{nearby cycle functor} \index{vanishing cycle functor}
of the {\it nearby and vanishing cycle functors} of Deligne, as will be explained  and discussed in later sections (see Remark  \ref{rNMD}):
\be  \label{eq:NMD=vancycle}
NMD(\cFb,S) \simeq \left(\,\varphi_{g|_{N\cap X}}[-1](\cFb|_{N\cap X})\,\right)_{x} 
\ee
and
\be
R\Gamma(l_{X},\cFb)\simeq \left(\,\psi_{g|_{N\cap X}}(\cFb|_{N\cap X})\,\right)_{x} \:,
\ee
resp.,
\be
R\Gamma(l_{X},\partial l_{X},\cFb) \simeq k_x^! \left(\,\psi_{g|_{N\cap X}}(\cFb|_{N\cap X})\,\right) \:.
\ee

One also has for a holomorphic function germ $g\colon  (X,x)\to (\bC,0)$ and
$\cFb \in D^b_{\cS-wc}(X;R)$ an isomorphism (see Corollary \ref {van-stalk} or, e.g., \cite[Corollary 1.1.1]{Sc})
\be \label{eq:LMD=vancycle}
LMD(\cFb,Re(g),x) = \bigl(\,R\Gamma_{\{Re(g)\geq 0\}}(\cFb)\,\bigr)_{x} \simeq
\left(\,\varphi_{g}[-1](\cFb)\,\right)_{x} \:.
\ee

\bex[Normal Morse data for products]\label{NMDproduct} Let $X_i$ be a complex algebraic (or analytic) variety, with 
 $\cF_i^{\bullet}$ weakly constructible
with respect to a Whitney stratification $\cS_i$ of $X_i$, for $i=1,2$. Then $\cF_1^{\bullet}\overset{L}{\boxtimes}  \cF_2^{\bullet}$ is
  weakly constructible with respect to the product Whitney stratification $\cS_1\times \cS_2$ of $X_1\times X_2$ (with strata $S_1\times S_2$ for
$S_1\in \cS_1$ and $S_2\in \cS_2$). Moreover
$$NMD \left(\cF_1^{\bullet}\overset{L}{\boxtimes}  \cF_2^{\bullet},S_1\times S_2\right) \simeq
NMD(\cFb_1,S_1) \overset{L}{\otimes}  NMD(\cFb_2,S_2)\:.$$
\eex

Let $N_i\subset \bC^{n_i}$ be normal slices at $x_i\in S_i$, with a corresponding non-degenerate covector $dg_{i,x_i}\in T^*_{S_i}\bC^{n_i}$
for a holomorphic function germ $g_i\colon  (\C^{n_i},x_i)\to (\bC,0)$ as above ($i=1,2$). Then $N=N_1\times N_2\subset \bC^{n_1}\times \bC^{n_2}$ is a normal slice at $(x_1,x_2)\in S_1\times S_2$, with $dg_{(x_1,x_2)}=\left(dg_{1,x_1},dg_{2,x_2}\right)$ a corresponding non-degenerate covector 
defined by $g(z_1,z_2):=g_1(z_1)+g_2(z_2)$. Example \ref{NMDproduct} is then by (\ref{eq:NMD=vancycle}) a very special case of the
\index{Thom-Sebastiani Theorem for vanishing cycles}
{\it Thom-Sebastiani Theorem for vanishing cycles} stated in  Theorem \ref{tsvan}.


\subsection{Perverse sheaf description via normal Morse data}

The distinguished triangles of Proposition~\ref{prop:NMDtriangleC} can be used to get
(inductive) information about the ``normal Morse data''.
The advantage in the complex  context comes from the fact
that $X\cap N\cap\{g=w\}$ is a complex analytic space (of lower dimension)  with an induced {\it Whitney stratification}. 
Moreover, the distance function $r$ used for the
stratified Morse theory on $X\cap N\cap\{g=w\}$ (in the induction step)
is now {\it strongly plurisubharmonic} so that one has nice estimates
for the Morse index (of an approximation by a Morse function) of $r$ 
and $-r$.

\bd \label{q-convex function}  \index{q-convex function}
A $C^{k}$-function $f\colon  X\to \bR$ ($k\geq 2$)
is called {\it q-convex in} $x\in X$, 
if one can choose $\hat{f}\colon  (\bC^n,x)\to (\bR,f(x))$ in a local embedding $(X,x)\hookrightarrow (\bC^n,x)$ with $\hat{f}|_X=f$
so that its \index{Levi form} {\it Levi form} $L_{\hat{f}}(x)$ at the point $x$ has at most $q$ 
non-positive eigenvalues. Here the Levi form of $\hat{f}$ at the point $x$
is the Hermitian form given by the matrix of partial derivatives: 
$$L_{\hat{f}}(x) := 
\left(\; \frac{\partial^{2} \hat{f}}{\partial z_{j}\, \partial \bar{z}_{k}} (x) 
\;\right)_{1 \leq j,k \leq n} \:.$$
$f$ is called {\it q-convex} (in
a subset $Z\subset X$) if it is q-convex in
all points (of $Z$).
\ed

For example, a {\it strongly plurisubharmonic} \index{strongly plurisubharmonic}
$C^{k}$-function on an open subset $U$ of $\bC^{n}$ (with $k\geq 2$) is 
$0$-convex in all of $U$. In particular, the distance function 
$$r\colon  \bC^{n}\to \bR^{\geq 0}\,; \: r(z):=\sum_{i=1}^{n} \: 
z_{i}\bar{z}_{i}  
\quad \text{is $0$-convex in $\bC^{n}$.}$$

We next state the second {\it main result} for stratified Morse theory of constructible sheaves in the complex  context
(see \cite[Theorem 6.0.1]{Sc} for a more general formulation). We include here a complete proof based on Theorem \ref{thm:SMTB},  
for explaining the use of stratified Morse theory and  because similar ideas can also be used later on for other important results.
 
\bt\label{thm:plurisub}
Let $\cS$ be a Whitney 
stratification of the complex algebraic (or analytic) variety $X$, and consider a  
$C^{\infty}$-function $f\colon X\to \bR$.
Assume $a<b$ are  stratified regular values of $f$ with respect to $\cS$
such that $\{a\leq f \leq b\}$ is compact. Let $\cFb\in D^b_{\cS-(w)c}(X;R)$ be given.
\begin{enumerate}
\item Suppose $f$ is {\it q-convex} in $\{a\leq f \leq b\}$, and let $q':=\min\{q,\dim X\}$.
$$\cFb \in \;^{p}D^{\leq n}_{\cS}(X;R)\: \Rightarrow \:
R\Gamma(\{a\leq f \leq b\},\{f=a\},\cFb) \in \;^{p}D^{\leq n+q'}(\{pt\};R) \: .$$
If $R$ is a Dedekind domain and  $\cFb$ is $\cS$-constructible, then also
$$\cFb \in \;^{p^+}D^{\leq n}_{\cS}(X;R)\: \Rightarrow \:
R\Gamma(\{a\leq f \leq b\},\{f=a\},\cFb) \in \;^{p^+}D^{\leq n+q'}(\{pt\};R) \: .$$
\item Suppose $-f$ is {\it q-convex} in $\{a\leq f \leq b\}$, and let $q':=\min\{q,\dim X\}$.
$$\cFb \in  \;^{p}D^{\geq n}_{\cS}(X;R) \: \Rightarrow 
\: R\Gamma(\{a\leq f \leq b\},\{f=a\},\cFb) \in \;^{p}D^{\geq n-q'}(\{pt\};R) \: .$$
If $R$ is a Dedekind domain and  $\cFb$ is $\cS$-constructible, then also
$$\cFb \in  \;^{p^+}D^{\geq n}_{\cS}(X;R) \: \Rightarrow 
\: R\Gamma(\{a\leq f \leq b\},\{f=a\},\cFb) \in\;^{p^+}D^{\geq n-q'}(\{pt\};R) \: .$$
\end{enumerate} 
Here $R\Gamma(\{a\leq f \leq b\},\{f=a\},\cFb)\in D^b_{c}(\{pt\};R)$ in case $\cFb$ is $\cS$-constructible.
\et

\begin{proof}
We prove the claim by induction over $\dim X$, where the case
$\dim X =0$ is trivial (since the conditions to be checked contain the zero object and are stable by finite direct sums).

After approximating $f$ by a stratified Morse function (see \cite{Ben}), we can assume that $f$ has only  finitely many  stratified Morse critical points
in $\{a\leq f \leq b\}$.
If the approximation is close enough, then also this new $\pm f$ is q-convex.
This implies the following important estimate for the 
Morse index $\tau(f|_S)$ of $f|_S$ for a stratum $S\in \cS$ in terms of $q':=\min\{q,\dim(X)\}$
(compare with \cite[p.191, Lemma 4.A.2]{GMs}):
\be
\label{eq:est.index}
\begin{cases}
\tau(f|_S) \leq \dim S+q', & \quad \text{if $f$ is q-convex.}\\
\tau(f|_S) = 2\dim S -\tau(-f|_S)  \geq \dim S-q', & \quad \text{if $-f$ is q-convex.}
\end{cases}
\ee

Then we apply  Lemma \ref{lem:locmorse}.
Since the conditions to be checked contain the zero object and are stable by extensions,
 we only have to show that the
``local Morse datum''
$$LMD(\cFb,f,x) = (R\Gamma_{\{f\geq f(x)\}}(\cFb))_{x} $$
of $f$ at a stratified Morse
critical point $x$ in a stratum $S\in \cS$  belongs to $\;^{p^{(+)}}D^{\leq n+q'}(\{pt\};R)$ for  statement 1.
(resp., belongs to $\;^{p^{(+)}}D^{\geq n+q'}(\{pt\};R)$ for  statement 2.).

By Theorem \ref{thm:SMTB} we get
$$LMD(\cFb,f,x) \simeq 
NMD(\cFb,S)[-\tau] ,$$
with $NMD(\cFb,S)$ the ``normal Morse datum'' of $\cFb$ in the stratum $S$
and $\tau$ the Morse index of $f|_S$.  
Here we can use the 
estimate~(\ref{eq:est.index}) for this Morse index.
Since the conditions to be checked contain the zero object and are stable by extensions as well as shifting by $[1]$ (resp., $[-1]$) for the first (resp., second) result, we only have to show that (with $s:=\dim S$):
\be\label{eq:NMD1}
NMD(\cFb,S) \in \begin{cases}
\;^{p^{(+)}}D^{\leq n-s}(\{pt\};R)
\quad & \text{in part 1. of the theorem.}\\
\;^{p^{(+)}}D^{\geq n-s}(\{pt\};R)
\quad & \text{in part 2. of the theorem.}
\end{cases} \ee

By Proposition~\ref{prop:NMDtriangleC} we have two distinguished triangles:
$$\begin{CD}
R\Gamma(l_{X},\cFb)[-1] @>>> 
NMD(\cFb,S) @>>>
k_{x}^{*}(\cFb|_{N\cap X}) @>>> \:,
\end{CD} $$
and
$$\begin{CD}
k_{x}^{!}(\cFb|_{N\cap X})  @>>> 
NMD(\cFb,S) @>>>
R\Gamma(l_{X},\partial l_{X},\cFb)[-1]  @>>> \:,
\end{CD} $$
where we use, as before, the notation (with $0<|w| \ll \delta \ll 1$)
$$ 
(l_{X},\partial l_{X}) := (X\cap N\cap\{r\leq\delta, g=w\},X\cap N\cap\{r=\delta, g=w\}) 
$$
for the \index{complex link} {\it complex link} of $X$ with respect to $S\in \cS$.
Also, $k_{x}\colon  \{x\}\hookrightarrow N\cap X$
is the inclusion, and $r$ is the {\it 0-convex distance function} 
$r(z):=\sum_{i=1}^{n} \: z_{i}\bar{z}_{i}$ in our fixed local embedding $(X,x)\hookrightarrow (\bC^n,x=0)$.
   
But $N\cap\{g=w\}$ is {\it transversal} to $\cS$ near $x$ (for 
 $0<|w|\ll 1$). So 
$$L:=X\cap N\cap\{g=w\}$$
gets an induced Whitney  stratification $\cS|_L$. Moreover,
 for $0<|w|\ll \delta\ll 1$, $\delta $ is a  regular value of $r$ with respect to
$\cS|_L$, since $dg_{x}$ is non-degenerate.
Then we can use the induction hypothesis with $q=0$ for
$$ \begin{cases}
\cFb|_L,\:r \quad &\text{and $[a,b]:=[-1,\delta]$ in part 1. of the theorem.}\\
\cFb|_L,\:-r 
\quad &\text{and $[a,b]:=[-\delta,1]$ in part 2. of the theorem.}
\end{cases} $$
Note that 
$$\cFb|_L[-1] \in \;^{p^{(+)}}D^{\leq n-s}_{\cS|_L}(L;R)
\quad \text{or} \quad  \;^{p^{(+)}}D^{\geq n-s}_{\cS|_L}(L;R)$$
by transversality and Proposition~\ref{prop-tr} (or its proof), with $\codim L=s+1$.

By the induction hypothesis we get 
$$ R\Gamma(l_{X},\cFb)[-1] \in\;^{p^{(+)}}D^{\leq n-s}(\{pt\};R)$$
or
$$ R\Gamma(l_{X},\partial l_{X},\cFb)[-1] \in \;^{p^{(+)}}D^{\geq n-s}(\{pt\};R) \:.$$
But, by definition, $k_{x}^{*}(\cFb|_{N\cap X}) \in \;^{p^{(+)}}D^{\leq n-s}(\{pt\};R)$,
which implies part 1. of the theorem  by the first distinguished triangle above.
Similarly
$$\cFb|_{X\cap N} \in \;^{p^{(+)}}D^{\geq n-s}_{\cS|_N}(X\cap N;R)$$
 by Proposition~\ref{prop-tr} (or its proof), so that $k_{x}^{!}(\cFb|_{N\cap X}) \in \;^{p^{(+)}}D^{\geq n-s}(\{pt\};R)$. This implies part 2. of the theorem  by the second distinguished triangle above.
The last claim of Theorem \ref{thm:plurisub} follows by inspection of the proof.
\end{proof}

\bc\label{cor:dimension}
Let $\bar{X}$ be a  compact complex algebraic (or analytic) variety of dimension $d$,
with $\cFb\in D^b_{wc}(\bar{X};R)$ (resp.,  $\cFb\in D^b_{c}(\bar{X};R)$ and $R$ is a  Dedekind domain 
in case we consider the dual perverse t-structure). Let $j\colon  X\hookrightarrow \bar{X}$ be the inclusion of the open dense complement of a closed algebraic (or analytic) subset $Z\subset \bar{X}$ (so that $d=\dim X=\dim \bar{X}$). Then:
\begin{enumerate}
\item 
$j^*\cFb \in \;^{p^{(+)}}D^{\leq n}(X;R)\: \Rightarrow \: Rj_!j^*\cFb \in \;^{p^{(+)}}D^{\leq n}(\bar{X};R)\:$
$$\Rightarrow \:
R\Gamma_c(X,j^*\cFb)\simeq R\Gamma(\bar{X},Rj_!j^*\cFb)  \in \;^{p^{(+)}}D^{\leq n+d}(\{pt\};R) \: .$$
\item 
$j^*\cFb \in  \;^{p^{(+)}}D^{\geq n}(X;R) \: \Rightarrow \: Rj_*j^*\cFb \in \;^{p^{(+)}}D^{\geq n}(\bar{X};R) \:$
$$\Rightarrow 
\: R\Gamma(X,j^*\cFb) \simeq R\Gamma(\bar{X},Rj_*j^*\cFb) \in \;^{p^{(+)}}D^{\geq n-d}(\{pt\};R) \: .$$
\end{enumerate} 
\ec

Note that by a theorem of Nagata (see, e.g., \cite{Luet, Con}) such a compactification is always available in the complex algebraic context.

The (method of) proof of Theorem \ref{thm:plurisub} also implies the following.
\bp\label{prop:T-stalk}
Let $\cS$ be a Whitney 
stratification of the complex algebraic (or analytic) variety $X$, and consider a  
differentiable function $f\colon  X\to \bR$.
Assume $a<b$ are  stratified regular values of $f$ with respect to $\cS$
such that $\{a\leq f \leq b\}$ is compact. 
Let $T\subset D^b(\{pt\};R)$ be a fixed \index{null system} ``null system'', i.e., a full triangulated subcategory stable by isomorphisms.
\begin{enumerate}
\item Denote by $D^b_{(\cS-)T-stalk}(X;R)$ the induced ``null system'' given by all complexes $\cFb\in D^b_{(\cS-)wc}(X;R)$ with stalks
$i_x^*\cFb\in T$ for all $x\in X$. Then 
$$\cFb \in D^b_{\cS-T-stalk}(X;R)\: \Rightarrow \: R\Gamma(\{a\leq f \leq b\},\{f=a\},\cFb)\in T \:.$$
\item 
Denote by $D^b_{(\cS-)T-costalk}(X;R)$ the induced ``null system'' given by all $\cFb\in D^b_{(\cS-)wc}(X;R)$ with  costalks
$i_x^!\cFb\in T$ for all $x\in X$. Then 
$$ \cFb \in D^b_{\cS-T-costalk}(X;R)\: \Rightarrow \: R\Gamma(\{a\leq f \leq b\},\{f=a\},\cFb)\in T\:.$$
\end{enumerate}
\ep

\bc\label{cor-T-stalk}
Let $\bar{X}$ be a  compact complex algebraic (or analytic) variety,
with $\cFb\in D^b_{wc}(\bar{X};R)$.  Let $j\colon  X\hookrightarrow \bar{X}$ be the inclusion of the open  complement of a closed algebraic (or analytic) 
subset $Z\subset \bar{X}$, and fix a given ``null system'' $T\subset D^b(\{pt\};R)$.
\begin{enumerate}
\item 
$j^*\cFb \in D^b_{T-stalk}(X;R)\: \Rightarrow \: Rj_!j^*\cFb \in D^b_{T-stalk}(\bar{X};R)\:$
$$\Rightarrow \:
R\Gamma_c(X,j^*\cFb)\simeq R\Gamma(\bar{X},Rj_!j^*\cFb)  \in T\: .$$
\item 
$j^*\cFb \in  D^b_{T-costalk}(X;R) \: \Rightarrow \: Rj_*j^*\cFb \in D^b_{T-costalk}(\bar{X};R) \:$
$$\Rightarrow 
\: R\Gamma(X,j^*\cFb) \simeq R\Gamma(\bar{X},Rj_*j^*\cFb) \in T \: .$$
\end{enumerate} 
\ec

\bex Let us mention here some important examples of such a ``null system'' $T\subset D^b(\{pt\};R)$.
\begin{enumerate}
\item $T=D^b_{c}(\{pt\};R)$ are the complexes with finitely generated cohomology.
Then $D^b_{(\cS-)T-stalk}(X;R)$ is the category of ($\cS$-)constructible sheaf complexes.
\item Assume $R$ is a Dedekind domain and $T=D^b_{tc}(\{pt\};R)$ are the complexes with finitely generated  torsion cohomology.
Then $D^b_{(\cS-)T-stalk}(X;R)$ is the category of ($\cS$-)constructible torsion sheaf complexes.
\item Assume $R$ is a field and $T=D^b_{c,\chi=0}(\{pt\};R)$ are the complexes with finitely generated  cohomology whose Euler characteristic is zero.
Then $D^b_{(\cS-)T-stalk}(X;R)$ is the category of ($\cS$-)constructible sheaf complexes $\cFb \in D^b_{(\cS-)c}(X;R)$
with  associated ($\cS$-)constructible function 
$$\chi_{stalk}\left(\cF^{\bullet}\right) = 0\in CF_{(\cS)}(X)\:,$$
i.e., whose stalkwise Euler characteristic $\chi_{stalk}\left(\cF^{\bullet}\right)(x):=\chi\left(\cF^{\bullet}_x\right)$ vanishes for all $x\in X$.
Then $[\cFb] \in K_0\left(D^b_{(\cS-)c}(X;R)\right)$ is in the kernel of 
$$
\chi_{stalk}\colon  \: K_0\left(D^b_{(\cS)-c}(X;R)\right)\to CF_{(\cS)}(X)\:
\iff \cFb \in D^b_{(\cS-)T-stalk}(X;R)\:.$$
\end{enumerate}
\eex

\bex \label{chi-compact}
Let $\bar{X}$ be a  compact complex algebraic (or analytic) variety,
with a given Whitney stratification $\cS$. Let  $\cFb,\cGb\in D^b_{\cS-c}(\bar{X};R)$ and assume $R$ a field.  Let $j\colon  X\hookrightarrow \bar{X}$ be the inclusion of the open  complement of a closed algebraic (or analytic) 
subset $Z\subset \bar{X}$, which is a union of strata $S\in \cS$, with $\cS|_X$ the induced Whitney stratification of $X$.
Then $\chi_{stalk}(j^*\cFb)=\chi_{stalk}(j^*\cGb)$ as a constructible function  implies 
$$\chi_c(X,j^*\cFb)=\chi_c(X,j^*\cGb)\:.$$
So the global Euler characteristic with compact support $\chi_c(X,j^*\cFb)$ only depends on the underlying constructible function
$\chi_{stalk}(j^*\cFb)=\alpha \in CF_{\cS}(X)$, with
\be \chi_c(X,\cFb)= \int_X \alpha d\chi := \sum_{S\in \cS|_X}\: \chi_c(S)\cdot \alpha(S) \:.
\ee
Here $\chi_c(S):=\chi_c(S,R_S)=\chi(H_c^*(S;R))$ is the corresponding \index{Euler characteristic}
Euler characteristic of a stratum $S\in \cS|_X$.

Similar considerations apply if $X$ is given as  a relatively compact open subset $X=\{f<b\}$ (resp., a compact subset $X=\{f\leq b\}$) for a proper differentiable function $f\colon  \bar{X}\to [a,d[\;\subset \bR$,
with $b<d$ a stratified regular value of $f$ with respect to $\cS$. Here $\bar{X}$ does not need to be compact, with strata $S$ of $\cS|_{\{f<b\}}$ given by the connected components of $S'\cap \{f < b\}$ for $S'\in \cS$. Similarly for the decomposition $\cS|_{\{f\leq b\}}$ of $\{f\leq b\}$, but here the parts
$S'\cap \{f \leq  b\}$ for $S'\in \cS$ are complex manifolds with possible non-empty boundary $S'\cap \{f = b\}$.
\eex

The proof of Theorem \ref{thm:plurisub} implies directly  the following characterizations (see, e.g., \cite[Corollary 6.0.2]{Sc}).
\bc \label{cor:NMD} 
Let $\cS$ be a Whitney 
stratification of the complex algebraic (or analytic) variety $X$, and
consider $\cFb\in D^b_{\cS-wc}(X;R)$ (resp.,  $\cFb\in D^b_{c}(X;R)$ and $R$ is a  Dedekind domain 
in case we consider the dual perverse t-structure).
Then we have:
\begin{enumerate}
\item 
\begin{center}
$\cFb \in \, ^{p^{(+)}}D^{\leq 0}_{\cS}(X;R)$ $\iff$ 
$NMD(\cFb,S)[-\dim(S)] \in \, ^{p^{(+)}}D^{\leq 0}(\{pt\};R) $\\
$\quad$ for all strata $S\in \cS$.
\end{center}
\item 
\begin{center}
$\cFb \in \, ^{p^{(+)}}D^{\geq 0}_{\cS}(X;R) $ $\iff$ 
$NMD(\cFb,S)[-\dim(S)] \in \, ^{p^{(+)}}D^{\geq 0}(\{pt\};R) $\\
$\quad$ for all strata $S\in \cS$.
\end{center}
\item Let $T\subset D^b(\{pt\};R)$ be a fixed ``null system''. Then \index{null system} 
\begin{eqnarray*}
\cFb \in D^b_{\cS-T-stalk}(X;R) &\iff &
NMD(\cFb,S)\in T
\quad \text{for all strata $S\in \cS$}\\
&\iff& \cFb \in D^b_{\cS-T-costalk}(X;R)\:.
\end{eqnarray*}
\end{enumerate}
\ec

\begin{proof}
The implications $\Rightarrow$ are already contained in (\ref{eq:NMD1}) from the proof of Theorem \ref{thm:plurisub}.
For the other implications $\Leftarrow$ in case of a point stratum $S=\{x\}$ one applies the  proof of Theorem \ref{thm:plurisub}
for the (co)stalk description
$$ i_x^*\cFb\simeq  R\Gamma(\{r\leq \delta\},\cFb) \simeq
 R\Gamma(\{-1\leq r\leq \delta\},\{r=-1\},\cFb)$$
and
$$i_x^!\cFb\simeq R\Gamma_c(\{r<\delta\},\cFb) \simeq
R\Gamma(\{-\delta \leq -r \leq 1\},\{-r=-\delta\},\cFb)\:,$$
with $0<\delta$ small and $r$ the {\it 0-convex distance function} 
$r(z):=\sum_{i=1}^{n} \: z_{i}\bar{z}_{i}$ in some fixed local embedding $(X,x)\hookrightarrow (\bC^n,x=0)$.
The general case is reduced to a point stratum via intersecting with a normal slice $N$ to $S\in \cS$ at a given point $x\in S$.
Here, the normal Morse data $NMD(\cFb,S')=NMD(\cFb|_N, S'\cap N)$ for $S'\in \cS$ close to $x$ does not change.
\end{proof}

\bex\label{perverse-NMD}
Let $\cS$ be a Whitney 
stratification of the complex algebraic (or analytic) variety $X$, with $\cFb\in D^b_{\cS-wc}(X;R)$. Then $\cFb$ is a \index{perverse sheaf} {\it perverse sheaf}\, if and only if 
$$H^i\left(NMD(\cFb,S)[-\dim(S)]\right)=0 \quad \text{for $i\neq 0$ and all strata $S\in \cS$.}$$
If, moreover, $\cFb$ is constructible and $R$ a Dedekind domain, then $\cFb$ is a \index{strongly perverse sheaf} {\it strongly perverse sheaf}\, if and only if, in addition,
$$H^0\left(NMD(\cFb,S)[-\dim(S)]\right)\quad \text{ is finitely generated and torsion-free }$$
for all strata $S\in \cS$.
\eex

\bex \label{chi-NMD}
Let $\cS$ be a Whitney 
stratification of the complex algebraic (or analytic) variety $X$, with $\cFb, \cGb\in D^b_{\cS-c}(X;R)$ and $R$ a field.
Then $\chi_{stalk}(\cFb)=\chi_{stalk}(\cGb)$ as an $\cS$-constructible function if and only if
$$\chi\left(NMD(\cFb,S)\right) = \chi\left(NMD(\cGb,S)\right) \quad \text{for all strata $S\in \cS$.}$$
So the Euler characteristic $\chi\left(NMD(\cFb,S)\right)=:\chi\left(NMD(\alpha,S)\right)$ for $S\in \cS$ depends only on the underlying $\cS$-constructible function
$\alpha=\chi_{stalk}(\cFb)$.
\eex

\br Corollary \ref{cor:NMD} also shows that the corresponding properties of the normal Morse data do not depend on the choice of the
Whitney stratification $\cS$ of $X$ with $\cFb$ $\cS$-(weakly) constructible.
The given characterization of the perverse t-structure in terms of normal Morse data is equivalent (in an embedded context) to
\cite[Theorem 10.3.12]{KS}, if one identifies the  normal Morse data with the corresponding ``micro-local type''
in the sense of Kashiwara-Schapira \cite[Proposition 6.6.1(ii)]{KS}
(as shown in \cite[Equation (5.52) on p.311]{Sc}).
The characterization of perverse sheaves in terms of normal Morse data  corresponds to 
the ``dimension axiom'' in the Morse theoretic approach to perverse sheaves
as outlined in \cite{Mac}. 
It also implies  the {\it purity result} stated in \cite[p.223, 6.A.3]{GMs}.
\er

Here are the final applications of Theorem \ref{thm:plurisub} for this section (see, e.g., \cite[Corollary 6.0.5]{Sc}),
which for $Y$ a point space reduce to Corollary \ref{cor:dimension}.

\bp \label{prop:fiberdim}
Let $\bar{f}\colon  \bar{X} \to Y$ be a {\it proper} holomorphic map of complex varieties,
with $Z\subset \bar{X}$ a closed complex subvariety. Consider the induced holomorphic
map
$$f:= \bar{f}|_X\colon \, X:=\bar{X}\backslash Z \to Y \:.$$
Suppose the \index{fiber dimension} {\it fiber dimension of $f$ is bounded above by $d$}, and consider
 $\cFb\in D^b_{(w)c}(X;R)$ (resp.,  $\cFb\in D^b_{c}(X;R)$ and $R$ is a  Dedekind domain 
in case we consider the dual perverse t-structure). Then:
\begin{enumerate}
\item  $\cFb|_X \in \, ^{p^{(+)}}D^{\geq n}(X;R) \:\Rightarrow \: Rf_{*}(\cFb|_X) \in \, ^{p^{(+)}}D^{\geq n-d}(Y;R)$.
\item  $\cFb|_X \in \, ^{p^{(+)}}D^{\leq n}(X;R) \:\Rightarrow \: Rf_{!}(\cFb|_X) \in \, ^{p^{(+)}}D^{\leq n+d}(Y;R)$.
\end{enumerate}
\ep

Note that by a theorem of Nagata (see, e.g., \cite{Luet, Con}) such a (partial) compactification $\bar{f}\colon  \bar{X} \to Y$
is always available in the complex algebraic context for a morphism $f\colon  X\to Y$.

\bc \label{cor:fiberdim}
Let $f\colon X\to Y$ be a morphism of complex algebraic varieties.
Suppose the \index{fiber dimension} {\it fiber dimension of $f$ is bounded above by $d$}, and consider
 $\cFb\in D^b_{(w)c}(X;R)$ (resp.,  $\cFb\in D^b_{c}(X;R)$ and $R$ is a  Dedekind domain 
in case we consider the dual perverse t-structure). Then:
\begin{enumerate}
\item  $\cFb \in \, ^{p^{(+)}}D^{\geq n}(X;R) \:\Rightarrow \: Rf_{*}\cFb \in \, ^{p^{(+)}}D^{\geq n-d}(Y;R)$.
\item  $\cFb \in \, ^{p^{(+)}}D^{\leq n}(X;R) \:\Rightarrow \: Rf_{!}\cFb \in \, ^{p^{(+)}}D^{\leq n+d}(Y;R)$.
\end{enumerate}
\ec

\bex[Finite morphism] \label{finite} \index{finite morphism}
Let $f\colon  X\to Y$ be a {\it finite} morphism of complex algebraic (or analytic) varieties (i.e., $f$ is proper with fiber dimension $d=0$).
Then  $Rf_!=Rf_*=f_*$ is {\it t-exact} with respect to the perverse t-structure, and in case $R$ a Dedekind domain also with respect to the dual
t-structure.
\eex


\subsection{Characteristic cycles and index theorems}

In this section, we give an  introduction to the theory
of Lagrangian cycles in the complex analytic and algebraic context, 
using the language of {\it stratified
 Morse theory for constructible functions and sheaves}, as developed in the previous sections.
We explain from this viewpoint the Euler isomorphism between 
constructible functions and Lagrangian cycles, together with some index theorems.  \\

The theory of Lagrangian cycles in the complex analytic context started in $1973$
with Kashiwara's local index formula for holonomic D-modules, as formulated in
\cite{Ka1} and proved in \cite[Chapter 6]{Ka2} (compare also with \cite{Mal}).
Kashiwara introduced for this a local invariant of a singular complex analytic set,
which around the same time was independently introduced by MacPherson 
as the {\it local Euler obstruction} in his celebrated  work \cite{Mch} on
Chern homology
classes for singular complex algebraic varieties (establishing a conjecture
of Grothendieck and Deligne). 
It was Dubson \cite{Du, Du1, BDK} who observed some years later that these two invariants are the same. The definitions
of Kashiwara and MacPherson are of transcendental nature. A purely algebraic 
definition of the "local Euler obstruction" was found by Gonzalez-Sprinberg and
Verdier \cite[Expos\'{e} 1]{CEP} (compare with \cite[Example 4.2.9]{Ful}).\\

In this section we work in a global embedded context, with $k\colon  X\hookrightarrow M$ the closed embedding of a complex algebraic (or analytic) variety $X$ into a complex algebraic (or analytic) manifold $M$ of dimension $\dim M=m$ (with $M$ pure-dimensional or $\dim M$ viewed as a locally constant function).
Let $\cS$ be a given Whitney stratification of $X$, with conormal space
$$T^*_{\cS}M:=\bigcup_{S\in \cS}\: T^*_SM \hookrightarrow T^*M|_X \hookrightarrow T^*M\:.$$
Here $T^*_{\cS}M$ is a closed subset by the {\it Whitney a-condition} of our stratification $\cS$. 
The open subset of {\it non-degenerate} covectors is given by \index{non-degenerate covector}
$$\left(T^*_{\cS}M\right)^{\circ}= \bigcup_{S\in \cS}\: \left(T^*_SM\right)^{\circ}\:,$$
with
$$\left(T^*_SM\right)^{\circ}:=T^*_SM\backslash \bigcup_{S\neq S'\in \cS} \overline{T^*_{S'}M}\:.$$

\bd The abelian group $L(\cS,T^*M)$ of ($\bC^*$-conic) \index{Lagrangian cycles} {\it Lagrangian cycles} in $T^*_{\cS}M$ is given by
$$L(\cS,T^*M):=H^{BM}_{2m}\left(T^*_{\cS}M;\bZ\right)\:,$$
so that a corresponding Lagrangian cycle is uniquely given by a (locally) finite sum
\be\label{Lcycle}
\sum_{S\in \cS} m(S)\cdot \left[\overline{T^*_SM}\right] \quad \text{with $m(S)\in \bZ$,}
\ee
and $\left[\overline{T^*_SM}\right]$ the corresponding {\it fundamental class}.
If we do not want to fix the stratification $\cS$, then we denote by $L(T^*M|_X)$ the abelian group  of ($\bC^*$-conic) {\it Lagrangian cycles} in 
$T^*M|_X$ given by a similar (locally) finite $\bZ$-linear combination of fundamental classes
$$\left[T^*_ZM\right]:=\left[\overline{T^*_{Z_{reg}}M}\right]$$
with $Z\hookrightarrow X$ an irreducible closed algebraic (or analytic) subvariety.
\ed

\br In the complex algebraic context, the Borel-Moore homology group
$$H^{BM}_{2m}\left(T^*_{\cS}M;\bZ\right)\simeq A_m\left(T^*_{\cS}M \right)\simeq Z_m\left(T^*_{\cS}M\right)$$
is also the same as the corresponding {\it Chow- and cycle group} of $T^*_{\cS}M$ as in \cite{Ful}, since $T^*_{\cS}M$ is of  dimension 
$m$.
So many of the following results could also be stated in this language, but this does not apply to our method of proof based on {\it stratified Morse
theory for constructible functions and sheaves}. For this reason, we work in this section only with the homological language, which at the same time
also  applies to the complex analytic context.
\er

For simplicity, in this section we only consider ($\cS$-)constructible sheaf complexes in $D^b_{(\cS-)c}(-;R)$ with $R$ a field and their corresponding (compactly supported)
{\it Euler characteristic}, even though  most of our results and proofs work for more general ``stalk properties and corresponding additive functions''
(see, e.g., \cite[Section 5.0.3]{Sc}).

\bd \label{CC} The \index{characteristic cycle} {\it characteristic cycle} $CC(\cFb)\in L(\cS,T^*M)$ of a constructible complex 
$\cFb\in D^b_{\cS-c}(X;R)$ is defined via the following multiplicities in (\ref{Lcycle}):
\be
m(S):=(-1)^{\dim S}\cdot \chi(NMD(\cFb,S)) \in \bZ\quad \text{for a stratum $S\in \cS$.}
\ee
By Example \ref{chi-NMD} the Euler characteristic 
$$\chi(NMD(\cFb,S))=\chi(NMD(\alpha,S))$$
 only depends on the associated $\cS$-constructible function
$\alpha = \chi_{stalk}(\cFb)\in CF_{\cS}(X)$ so that we can also define the characteristic cycle $CC(\alpha)$ of $\alpha \in CF_{\cS}(X)$ via
\be\label{l-integral}
(-1)^{\dim S}\cdot m(S):=\chi(NMD(\alpha,S)) = \alpha(x)-\int_{l_X} \alpha \,d\chi \in \bZ\:,
\ee
with $x\in S$, $(l_X,\partial l_X)$ the corresponding {\it complex link} of $X$ in $S\in \cS$ and $l_X^{\circ}:=l_X\backslash \partial l_X$.
This can also be reformulated as
\be
(-1)^{\dim S}\cdot m(S) = \alpha(x)- \sum_{S\subset  \partial S'} \:c(S,S')\cdot \alpha(S') \in \bZ\:,
\ee
with $c(S,S'):=\int_{l_X} 1_{S'} d \chi =   \int_{l_X^{\circ}} 1_{S'} d \chi =\chi_c(l_X^{\circ}\cap S',\bQ)$ for $S\subset \partial S'$ and 
$S,S'\in \cS$ a topological invariant of the Whitney stratification $\cS$ of $X$. Note  that 
\be \label{eq:Sullivan}
 \chi_{c}(l_{X} \cap S', \bQ) = 
\chi_{c}(l_X^{\circ}\cap S', \bQ)  \:,
\ee
since $\partial l_{X}$ is a compact real analytic Whitney stratified set with {\it odd-dimensional}
strata so that $\chi_{c}(\partial l_{X}\cap S', \bQ)=0$
(see, e.g., \cite[Lemma 5.0.3]{Sc} and \cite{Su}).
\ed

Here, the last equality in (\ref{l-integral}) follows from the distinguished triangle (\ref{clink-stalk}).
The choice of the sign $(-1)^{\dim S}$ will become clear in a moment.

\bex\label{CC-top} Let $X\hookrightarrow M$ be a {\it closed} complex submanifold of $M$, with strata $S\in \cS$ given by the connected components of $X$.
Let $\cL$ a local system of rank r on $X$. Then
$$CC(\cL)=CC(r\cdot 1_X)=(-1)^{\dim X}\cdot r\cdot \left[T^*_XM\right]\:.$$
Similarly, in the general context of the above Definition with $S\in \cS$ an {\it open} stratum \index{open stratum}
of $X$ and $\cFb\in D^b_{\cS-c}(X;R)$.
Then by  Example \ref{top-NMD},
the multiplicity $m(S)$ of the characteristic cycle $CC(\cFb)$ is given by $m(S)=(-1)^{\dim S}\cdot \chi_{stalk}(\cFb)(S)\in \bZ$.
\eex

\br\label{supp-CC} 
Let $\cFb\in D^b_{\cS-c}(X;R)$ be given, with $k\colon  X\hookrightarrow  M$ the global closed embedding into the complex manifold $M$. Then for a stratum $S\in \cS$ we get by definition \index{support of a characteristic cycle}
$$\overline{T^*_SM} \subset supp(CC(\cFb)) \iff \chi(NMD(\cFb,S)) \neq 0  \ .$$
 And, in general, this is much weaker than the corresponding condition
$$ 
\overline{T^*_SM}
\subset \mu supp(Rk_{*}\cFb) \iff 
NMD(\cFb,S) \not\simeq 0  $$
from Proposition \ref{pmusupp}. So in general the inclusion
\be
supp(CC(\cFb)) \subset  \mu supp(Rk_{*}\cFb) \quad \text{for $\cFb\in D^b_{\cS-c}(X;R)$}
\ee
does not need to be an equality. But Example \ref{perverse-NMD} implies
\be
supp(CC(\cFb)) = \mu supp(Rk_{*}\cFb) \quad \text{for $\cFb\in Perv_{\cS}(X;R)$.}
\ee
Moreover, Example \ref{perverse-NMD} also implies that $CC(\cFb)$ is an {\it effective} Lagrangian cycle for $0\neq \cFb\in Perv_{\cS}(X;R)$, i.e., the multiplicity
$m(S)\geq 0$ for all $S\in \cS$, with $m(S)>0$ for $S$ a top-dimensional stratum in $supp(\cFb)$.
\er 

By induction on $\dim X$ and Example \ref{CC-top}, one easily gets that the induced group homomorphism 
$$CC\colon  K_0(D^b_{\cS-c}(X;R))\to L(\cS,T^*M)=H^{BM}_{2m}\left(T^*_{\cS}M;\bZ\right)$$
is {\it surjective}. Moreover, $CC$ factorizes over $\chi_{stalk}$ and both homomorphisms
$$CC\colon  K_0(D^b_{\cS-c}(X;R))\to L(\cS,T^*M) \quad \text{and} \quad CC\colon  CF_{\cS}(X)\to L(\cS,T^*M)$$
have by Example \ref{chi-NMD} the {\it same kernel}, so that $CC$ induces for a fixed Whitney stratification $\cS$ of $X$ an isomorphism of abelian groups
\be \begin{CD}
CC\colon  CF_{\cS}(X) @> \sim >> L(\cS,T^*M)=H^{BM}_{2m}\left(T^*_{\cS}M;\bZ\right)\:.
\end{CD} \ee

Finally, $CC(\cFb)\in L(T^*M|_X)$ (and then also $CC(\alpha)\in L(T^*M|_X)$) does not depend on the choice of the Whitney stratification $\cS$ with
$\cFb\in D^b_{\cS-c}(X;R)$. In fact, if $\cT$ is another Whitney stratification of $X$ which refines $\cS$, then any stratum $T\in \cT$ is contained in a stratum $S\in \cS$. If $T\subset S$ is {\it open}, then $NMD(\cFb,S)=NMD(\cFb,T)$ just by the definitions. If $T\subset S$ is a proper closed subset, then
$NMD(\cFb,T)=0$. In fact, if $x\in T$ is a stratified Morse critical point of a function germ $f\colon (M,x)\to (\bR,f(x))$ with respect to $\cT$, then
$LMD(\cFb,f,x)\simeq NMD(\cFb,T)[-\tau]$ by Theorem  \ref{thm:SMTB}. But then $x\in S$ is not a stratified critical point of $f\colon (M,x)\to (\bR,f(x))$ with respect to $\cS$ so that $LMD(\cFb,f,x)=0$ by Lemma \ref{loc.const}. So in the limit over all such Whiney stratifications $\cS$ of $X$, we get a surjective group homomorphism 
$$CC\colon  K_0(D^b_{c}(X;R))\to L(T^*M|_X)\:,$$
which factorizes via $\chi_{stalk}$ over an isomorphism
\be \begin{CD}
CC\colon  CF(X)@>\sim >> L(T^*M|_X)\:.
\end{CD} \ee

Let us now explain the choice of the sign $(-1)^{\dim S}$ in the Definition \ref{CC} of the characteristic cycle $CC$.

The differentiable $C^k$-function germ  $f\colon (M,x)\to (\bR,f(x))$ (with $2\leq k\leq \infty$)  has a \index{stratified Morse critical point} {\it stratified Morse critical point} at $x\in S$, for a stratum $S\in \cS$,
if and only if (see, e.g., \cite[p.311]{KS} or \cite[p.286]{Sc}):
\be \begin{cases}
& df_{x}\in \left(T^*_SM\right)^{\circ}\; , \quad \text{i.e., $df_x\in T^*_SM$ is non-degenerate, and}\\
& \text{the graph $df(M)\subset T^{*}M$ of $df$ intersects $T_{S}^{*}M$
{\it transversally} at $df_{x}\,$.}
\end{cases} \ee
Using the {\it complex orientations} of $T^*M$ and $M\simeq df(M)$ one gets in this case for the corresponding
{\it local intersection number}: 
\begin{equation} \label{eq:Mindexintersec} 
(-1)^{\dim S}\cdot \sharp_{df_{x}}\bigl( \, [T_{S}^{*}M] \cap [df(M)] \,\bigr) = (-1)^{\lambda} \:,
\end{equation}
with $\lambda$ the \index{Morse index} {\it Morse index} of $f|_S$ at $x$
(i.e., its Hessian at $x$  has exactly $\lambda$ negative eigenvalues).
Here $[df(M)]\in H^{BM}_{2m}(df(M);\bZ)$ is the fundamental class of the oriented manifold $df(M)\simeq M$.

The local intersection number is defined similarly to the following {\it global
intersection number}
$$\sharp\colon  \: H^{BM}_{2m}(A; \bZ) \times  H^{BM}_{2m}(B; \bZ) \to \bZ$$
for two {\it  closed} subsets $A,B\subset T^{*}M$ with $A\cap B$ {\it compact}:
$$ \begin{CD}
H^{2m}_{A}(T^{*}M; \bZ) \times  H^{2m}_{B}(T^{*}M; \bZ) @> \cup >> H^{4m}_{c}(T^{*}M; \bZ) @> tr >> \bZ \\
@V PD V \wr V  @V PD V \wr V    @| \\
H^{BM}_{2m}(A;\bZ) \times  H^{BM}_{2m}(B; \bZ) @> \cap >> 
 H_{0}(T^{*}M;\bZ) @> deg >> \bZ 
\end{CD} $$
Here, $PD$ is {\it Poincar\'{e} duality} \index{Poincar\'{e} duality}
given by the cap-product with the 
fundamental class $[T^{*}M]$ of the complex manifold $T^{*}M$.

\bex\label{ex1337} Let $g\colon  (M,x)\to (\bC,g(x))$ be a holomorphic function germ so that the graph $dg(M)\subset T^{*}M$ of $dg$ intersects $T_{S}^{*}M$
{\it transversally} at $dg_{x}\,$. Then 
$$\sharp_{dg_{x}}\bigl( \, [T_{S}^{*}M] \cap [dg(M)] \,\bigr) = 1$$
and $g|_S$ has a \index{complex Morse critical point} {\it complex Morse} critical point at  $x\in S$. But if we identify the complex cotangent bundle $T^*M$ with the real cotangent bundle of the underlying real (oriented) manifold, then the graph $dg(M)$ of $g$ gets identified with the graph $df(M)$ of the real part $f:=Re(g)$ of $g$.
So, with the complex orientations, one also has
$$\sharp_{df_{x}}\bigl( \, [T_{S}^{*}M] \cap [df(M)] \,\bigr) = 1\:,$$
but $f|_S$ has at $x\in S$ a classical Morse critical point of index $\lambda=\dim S$.
\eex

The  intersection theory in the ambient cotangent bundle $T^*M$ 
 is used for the formulation of the following beautiful {\it intersection formula} \index{intersection formula}
(see, e.g., \cite[Theorem 5.0.4]{Sc} for a more general version).

\bt \label{thm:intersec}
Let $f\colon  M\to [a,d[\, \subset \bR$ be a $C^{\infty}$-function
($a \leq d \leq \infty$),
with $f|_X$ {\it  proper}.
Suppose that $T^*_{\cS}M \cap df(M)$ is {\it compact},
with $T^*_{\cS}M$ the union of conormal spaces $T_S^*M$ to the strata $S\in \cS$ of the Whitney stratification $\cS$ of $X\hookrightarrow M$.
Then we have for all $\cFb\in D^b_{\cS-c}(X;R)$:
\begin{eqnarray}
\dim_R H^*(X,\cFb) &<& \infty  \:  \text{, with} \nonumber \\
\chi(X,\cFb)  &=&  \sharp(\:CC(\chi_{stalk}(\cFb)) \cap [df(M)]\:) \:,
\end{eqnarray}
and
\begin{eqnarray}
\dim_R H^*_c(X,\cFb) &<& \infty\:  \text{, with} \nonumber \\
\chi_c(X,\cFb)  &=&  \sharp(\:CC(\chi_{stalk}(\cFb)) \cap [-df(M)]\:) . 
\end{eqnarray}
\et

\begin{proof} 
The reader should compare the following proof also with the proof of Theorem \ref{thm:plurisub}.
 By approximation \cite{Ben}, we can assume that $f$ has only {\it stratified Morse critical points}
with respect to $\cS$, since the corresponding intersection number does not change
by a homotopy argument.
By  assumption, $T^*_{\cS}M \cap df(M)$ is compact so that $f$ has only  finitely
many stratified Morse critical points. Choose $b\in [a,d[\,$ with all these
critical points contained in $\{a\leq f < b\}$. Then one gets by Lemma \ref{loc.const} 
$$R\Gamma(X\cap \{f \leq b\},\cFb)\simeq  R\Gamma(X,\cFb) \:,$$
and
$$R\Gamma_{c}(X\cap \{f < b\},\cFb)\simeq  R\Gamma_{c}(X,\cFb) \:.$$
And these complexes belong to $D^b_c(\{pt\};R)$ by Proposition \ref{prop:T-stalk}.
Similarly,  the normal Morse data $NMD(\cFb,S)\in D^b_c(\{pt\};R)$ by Corollary \ref{cor:NMD}.
Since the Euler characteristic $\chi$ is additive,  by Corollary \ref{cor:coarse}
and Lemma \ref{lem:locmorse} it is enough
to show that for such a stratified Morse critical point $x\in S$:
$$\chi(LMD(\cFb,\pm f,x)) =  
\sharp_{\pm \, df_{x}}\,(\:CC(\cFb) \cap [\pm \,df(M)]\:) \:.$$
But this follows from Theorem \ref{thm:SMTB}:
$$LMD(\cFb,\pm f,x)\simeq NMD(\cFb,S)[-\tau] \:,$$
with $\tau$ the Morse index of $\pm f|_S$ at $x$.
Finally, by the definition of $CC(\cFb)$ and  Equation (\ref{eq:Mindexintersec}) we have
\begin{eqnarray*}
\chi(LMD(\cFb,\pm f,x)) &=&  (-1)^{\tau} \cdot \chi(NMD(\cFb,S)) \\
&=& \sharp_{\pm \,df_{x}}\,(\:CC(\cFb) \cap [\pm \,df(M)]\:) \:,
 \end{eqnarray*}
 which completes the proof.
\end{proof}

This  intersection Theorem \ref{thm:intersec} goes back 
to Dubson \cite{Du2}, Sabbah \cite{Sab1}
and Ginsburg \cite{Gi1} (partially  in the context of
holonomic $D$-modules).\\

Note that the condition that $T^*_{\cS}M \cap df(M)$ is compact just means that the stratified critical locus $\Sing_{\cS}(f):=\bigcup_{S\in \cS} \Sing(f|_S)$
is compact. The intersection formula above can also be reformulated in the language of constructible functions.

Assume $f\colon M\to [a,d[\, \subset \bR$ is a {$C^{\infty}$-function
with $f|_X$ {\it proper} and $\pi\left(T^*_{\cS}M \cap df(M)\right)\subset [a,b]$ for some $b\in [a,d[\,$, 
with $\pi$ the projection $T^{*}M\to M$.
In particular, $T^*_{\cS}M \cap df(M)$ is {\it compact}. 

Then one also has the following counterpart
of Theorem~\ref{thm:intersec}
for a constructible function $\alpha \in CF_{\cS}(X)$:
\be \label{eq:intersecCF1}
\int_{X\cap \{f\leq r\}}  \alpha \: d\chi   = \sharp(\:CC(\alpha) \cap [df(M)]\:) 
\quad  \text{for all $r\in \,]b,d[\,$,}
\ee
and
\be \label{eq:intersecCF2}
\int_{X\cap \{f < r\}}  \alpha \: d\chi   = \sharp(\:CC(\alpha) \cap [-df(M)]\:)
\quad  \text{for all $r\in \,]b,d[\,$.}
\ee

Note that $r\in \,]b,d[\,$ is a stratified  regular value of $f|_X$ so that the left Euler characteristics are defined by Example \ref{chi-compact}.
Let us give some examples, where the above conditions on $f$ are satisfied.

\bex[Global index formula and Poincar\'{e}-Hopf theorem for singular spaces] \label{PH-sing} \index{global index formula} \index{Poincar\'{e}-Hopf theorem}
Let $X\hookrightarrow M$ be a  compact complex algebraic (or analytic) subvariety of $M$, with $\cS$ a given Whitney stratification of $X$. 
Then we can take for $f$ a constant function
so that $df(M)$ is the zero-section.
Then we get for any $\cS$-constructible function $\alpha\in CF_{\cS}(X)$:
\be \label{eq:PHopf}
\int_X \, \alpha \: d\chi   = \sharp(\:CC(\alpha) \cap [T^{*}_{M}M]\:) = \sharp(\:CC(\alpha) \cap [\omega(M)]\:)  
\ee
for any differentiable one-form  $\omega$ on $M$, i.e., a section of $T^*M\to M$.
Assume that $T^*_{\cS}M\cap \omega(M)$ is finite, i.e., the set  $\Sing_{\cS}(\omega):=\bigcup_{S\in \cS} \Sing(\omega|_S)$
of critical points of $\omega$ with respect to $\cS$ is finite. If one defines for $x\in \Sing_{\cS}(\omega)$ and $\alpha\in CF_{\cS}(X)$ the \index{local index}
{\it local index} $ind_x(\omega,\alpha)$ of $\omega$ with respect to $\alpha$ by 
\begin{equation} \label{eq:loc-index} 
ind_x(\omega,\alpha):= \sharp_{\omega_{x}}\bigl( \, CC(\alpha)\cap [\omega(M)] \,\bigr) \:,
\end{equation}
then one gets the Poincar\'{e}-Hopf theorem
\be
\int_X \, \alpha \: d\chi   = \sum_{x\in \Sing_{\cS}(\omega)}\:ind_x(\omega,\alpha)\:.
\ee
\eex

See also \cite{Se} for an overview of different indices of vector fields and one forms in the context of singular complex varieties.
The following recent application is due to \cite{AMSS2} (and implicitly already contained in \cite{FK,ST}).

\bex[Effective characteristic cycles on an abelian variety] \label{cc-abv}
Let $M=A$ be a complex abelian variety, \index{abelian variety} so that $A$ is a projective abelian algebraic group with
trivial cotangent bundle $T^*A$. Let $\alpha\in CF(A)$ be a constructible function with \index{effective characteristic cycles}
$CC(\alpha)$ an {\it effective} cycle in $T^*A$.
By Kleiman's transversality theorem, a {\it generic algebraic one-form} $\omega$ intersects $supp(CC(\alpha))$ only in finitely many points
$\omega_x$ (see, e.g., \cite[Proposition 2.8]{ST}). Then the local intersection number  (\ref{eq:loc-index}) as well as the global
Euler characteristic (\ref{eq:loc-index}) \index{signed Euler characteristic property} are non-negative:
\be
ind_x(\omega,\alpha)= \sharp_{\omega_{x}}\bigl( \, CC(\alpha)\cap [\omega(A)] \,\bigr)\geq 0 \quad \text{and} \quad \int_A \alpha \:d\chi\geq 0 \:. 
\ee
For example $\alpha$ could be given by
\begin{enumerate}
\item $\alpha=Eu_Z^{\vee}=(-1)^{\dim Z}\cdot Eu_Z$ is the dual Euler obstruction of a pure dimensional closed subvariety $Z\hookrightarrow A$
(see Equation (\ref{dual-euler}) below).
\item $\alpha=\chi_{stalk}(\cFb)$ for a perverse sheaf $\cFb\in Perv(A;R)$.
\item $\alpha=\chi_{stalk}(Rp_*\cFb)$ for a perverse sheaf $\cFb\in Perv(G;R)$ on a {\it semi-abelian variety} $G$, with $p\colon  G\to A$ the  projection onto
the corresponding abelian variety $A$. Then
$$\chi(G;\cFb)=\chi(A;Rp_*\cFb)=\int_A \alpha \:d\chi\geq 0\:.$$
See \cite[Proposition 8.4, Example 8.5]{AMSS2} for a more general class of morphisms $p\colon  G\to A$ to an abelian variety with this property
for any perverse sheaf $\cFb\in Perv(G;R)$
(for a different approach via {\it generic vanishing theorems} on semi-abelian varieties see (\ref{eq6}) in Section \ref{semi-ab}).
\end{enumerate}
\eex

\bex[Affine varieties and global Euler obstruction] \label{index-affine}
Let $X\hookrightarrow \bC^n$ be an {\it affine} complex algebraic variety endowed with a complex algebraic Whitney stratification $\cS$.
Then the {\it semi-algebraic} distance function $r\colon \bC^n\to [0,\infty[$,  $r(z):=\sum_{i=1}^n z_i\cdot \bar{z}_i$ has only finitely many critical values with respect to $\cS$. So one gets for $\cFb\in D^b_{\cS-c}(X;R)$: 
$$
\chi(X,\cFb)  =  \sharp(\:CC(\chi_{stalk}(\cFb)) \cap [dr(\bC^n)]\:) \:,
$$
and
$$
\chi_c(X,\cFb)  =  \sharp(\:CC(\chi_{stalk}(\cFb)) \cap  [-dr(\bC^n)]\:)\: . 
$$
Similarly, for a constructible function $\alpha \in CF_{\cS}(X)$:
$$
\int_{X\cap \{r\leq \delta\}}  \alpha \: d\chi   = \sharp(\:CC(\alpha) \cap  [dr(\bC^n)]\:) 
\quad  \text{for all $\delta>0$ large enough,}
$$
and
$$
\int_{X\cap \{r< \delta\}}  \alpha \: d\chi   = \sharp(\:CC(\alpha) \cap  [-dr(\bC^n)]\:)
\quad  \text{for all $\delta>0$ large enough.}
$$
If $X$ is pure dimensional and $CC(\alpha)=(-1)^{\dim X}\cdot \left[\overline{T^*_{X_{reg}}M}\right]$, then the intersection number 
$$Eu(X):= (-1)^{\dim X}\cdot \sharp\left(\left[\overline{T^*_{X_{reg}}M}\right] \cap  [dr(\bC^n)]\:\right)$$
is by \cite[Equation (5.6.4)]{Sc} the \index{global Euler obstruction} {\it global Euler obstruction} of $X$ in the sense of 
Seade-Tib\u{a}r-Verjovsky \cite{STV}
(and compare also with \cite[Theorem 1]{Du2}). The sign $(-1)^{\dim X}$ comes from different orientation conventions.
As we will see in a moment, here $\alpha=Eu_X$ is the \index{local Euler obstruction} {\it local Euler obstruction} function of
$X$ as defined by MacPherson \cite{Mch}.
\eex

\bex[Local index formula]
Let $f\colon M\to [0,d[\;\subset \bR $ be a {\it real analytic} function, with $f|_X$ proper.  Then 
$$f\circ \pi(\,T^*_{\cS}M \cap df(M)\,) \subset [0,d[ $$  
is {\it subanalytic}
and, by the {\it curve selection lemma}, it is {\it discrete}. Especially, 
$$\pi(\,T^*_{\cS}M  \cap df(U)\,)\subset X\cap\{f=0\} \:,$$
if we restrict to $f\colon U:=\{f< \epsilon\}\to [0,\epsilon[\,$, with $0<\epsilon$ small enough.
Then one gets for $\alpha\in CF_{\cS}(X)$ and $0<r<\epsilon$: 
\be \label{loc-intersec}
\int_{X\cap\{f=0\}} \,\alpha \: d\chi = 
\int_{X\cap \{f\leq r\}} \alpha \: d\chi  = \sharp(\:CC(\alpha) \cap [df(U)]\:) \:.
\ee
\eex

The most important special case  is the local situation
$(M,x)\simeq (\bC^{n},0)$ with $f(z):=r(z):=
\sum_{i=1}^{n} z_i\cdot \bar{z}_i$. In this case we get for $0<\epsilon \ll 1$:
\begin{equation} \label{eq:Euler} 
\alpha(x) = \sharp_{dr_{x}} (\:CC(\alpha) \cap [dr(\{r<\epsilon\})]\:) \:.
\end{equation}
In particular, this local intersection number is independent of the choice of
the real analytic  function $r$ with $\{r=0\}=\{x\}$.
Moreover, it defines the inverse
of the characteristic cycle map (see, e.g., \cite[Corollary 5.0.1]{Sc} and compare with \cite[Theorem 11.7]{Gi1}).

\bc \label{cor:invCC}
The inverse of the characteristic cycle map \index{inverse of characteristic cycles}
$$CC \colon  CF_{\cS}(X) \to  L(\cS;T^*M)$$
is given by 
$$Eu^{\vee}\colon  L(\cS;T^*M) \to  CF_{\cS}(X);$$
$$ [C] \mapsto \alpha\;, \quad \text{with} \:\:
\alpha(x):=  \sharp_{dr_{x}} (\:[C] \cap [dr(\{r<\epsilon\})]\:) \:,$$
for $r, \epsilon $ as above. In particular, for $CC(\alpha)=\sum_{S\in \cS} m(S)\cdot \left[\overline{T^*_SM}\right] $, we get for $x\in X$:
\be\label{eq:inCC}
\alpha(x)=\sum_{S\in \cS} m(S)\cdot Eu^{\vee}\left(\left[\overline{T^*_SM}\right] \right)\:.
\ee
\ec

Note that this corollary is also a {\it refinement} of \cite[Theorem 11.7]{Gi1},
since we work with a {\it fixed stratification}. Assume now that $X\hookrightarrow M$ is pure dimensional,
so that 
$$C=\left[T^*_XM\right]:=\left[\overline{T^*_{X_{reg}}M}\right] \in L(\cS;T^*M)$$
defines a corresponding Lagrangian cycle. Then one gets by \cite[Equation (5.35) and p.323-324]{Sc} that
\be\label{dual-euler}
Eu^{\vee}\left(\left[T^*_XM\right]\right) = Eu_X^{\vee}:=(-1)^{\dim X}\cdot Eu_X \in CF_{\cS}(X)
\ee
is the \index{dual local Euler obstruction} {\it dual local Euler obstruction} function $Eu_X^{\vee}$, with $Eu_X$ the famous \index{local Euler obstruction}  {\it local Euler obstruction} function of
$X$ as defined by MacPherson \cite{Mch} (compare also with \cite[Theorem 11.7]{Gi1}). 
The sign $(-1)^{\dim X}$ comes again from different orientation conventions.
In particular, we get without any calculation
that the {\it local Euler obstruction} $Eu_{X}$ is constructible with respect
to {\it any Whitney stratification} of $X$.
Since this is a local result, it is then also true for any pure-dimensional
complex  algebraic (or analytic) variety $X$ (without any embedding into a complex manifold).
This result is due to Dubson \cite[Proposition 1, Theorem 3]{Du}
and Brasselet-Schwartz \cite[p.125, Corollary 10.2]{BSchw}. Then the \index{local index formula}  {\it local index formula} (\ref{eq:inCC}) can be reformulated as
\be\label{eq:inCC2}
\alpha(x)=\sum_{S\in \cS} m(S)\cdot (-1)^{\dim S}Eu_{\bar{S}}(x) \quad \text{for $x\in X$}
\ee
and $CC(\alpha)=\sum_{S\in \cS} m(S)\cdot \left[\overline{T^*_SM}\right] $ for a given $\alpha\in CF_{\cS}(X)$.

\br \label{rem:D-CC}
The {\it local index formula}~(\ref{eq:inCC2}) 
goes back to a corresponding
{\it index formula} of Kashiwara \cite{Ka1} (compare with \cite[Theorem 6.3.1]{Ka2})
for the {\it solution complex} \index{solution complex}
$Sol_M(\cM):=Rhom_{\cD_M}(\cM,\cO_{M})$ of a {\it holonomic} $\cD_M$-module \index{holonomic} on the complex 
manifold $M$. Note that this solution complex is a  complex analytically (or algebraically)
constructible complex of sheaves of $\bC$-vector spaces, with finite dimensional
stalks (\cite[Theorem 11.3.7]{KS}). Kashiwara's formula corresponds to~(\ref{eq:inCC2}) 
for the constructible function
$$\alpha :=\chi_{stalk}(\,Rhom_{\cD_M}(\cM,\cO_{M})[\dim M]\,)\:.$$
 Moreover, he works directly with the {\it characteristic cycle}
of a holonomic $\cD_M$-module.
The corresponding multiplicities of characteristic cycles fit by \cite[Example 5.3.4]{Sc} with our conventions.
Similarly, Kashiwara introduced for his index formula some topological invariants,
which are nothing else but the {\it Euler obstructions} of the closures $\bar{S}$ for the strata $S\in \cS$ of a Whitney stratification $\cS$ of $M$ so that
the characteristic variety $char(\cM)\subset T^*_{\cS}M$.
But this fact was only observed later on by Dubson (compare with \cite{BDK}, and also with
\cite[Introduction, p.xiii]{Ka2}).

Other references for these formulae are \cite{Du2}, \cite[p.545]{BMM} and \cite[Theorem 8.2, Theorem 11.7, p.393]{Gi1}
(but with an incorrect sign in \cite[Corollary 6.19(b) and  Theorem 8.2]{Gi1}).
\er

So one gets in the complex analytic (or algebraic) context a commutative diagram (see also \cite{Gi1}):
\be \label{CC-D-mod} \begin{CD}
K_0(\text{holonomic $\cD_M$-modules}) @> CC >> L(T^*M) \\
@V DR_M  V Sol_M[\dim M] V @A \wr A CC A \\
K_0(D^b_c(M;\bC)) @> \chi_{stalk} >> CF(M) \:.
\end{CD}\ee

Here $DR_M(\cM)\simeq \,Rhom_{\cD_M}(\cO_{M},\cM)[\dim M]$ is the \index{De Rham complex} {\it De Rham complex} of the holonomic $\cD_M$-module $\cM$.
In the algebraic context one has of course to use the De Rham and Solution complex of the associated analytic $\cD_M$-module.
That $DR_M(\cM) $ and $Sol_M(\cM)$ for a holonomic $\cD_M$-module $\cM$ have the same characteristic cycle follows, e.g., from the fact that they are exchanged by the 
{\it duality} of holonomic $\cD_M$-modules, resp., {\it Verdier duality} for constructible sheaf complexes
(see, e.g., \cite[Corollary 4.6.5]{HTT}). Finally, also the left vertical arrows become isomorphisms by the famous {\it Riemann-Hilbert correspondence} \index{Riemann-Hilbert correspondence}
(see, e.g., \cite[Theorem 7.2.1]{HTT}), if we restrict ourselves to \index{regular holonomic}
{\it regular holonomic} $\cD_M$-modules.\\


\subsection{Functorial calculus of characteristic cycles}

We continue to work in an embedded context and discuss the functorial calculus of characteristic cycles (see also, e.g., \cite{Gi1, Sab1} and compare with \cite[Chapter IX]{KS} for counterparts in real geometry).
We explain the translation
into the context of Lagrangian cycles of the following operations for
constructible functions and sheaves: external product, proper direct image, non-characteristic pullback and  specialization (i.e., nearby cycles),
together with an intersection formula for
vanishing cycles. \\

 Note that by Corollary \ref{cor-T-stalk} and Corollary \ref{cor:NMD}, the constructible complexes 
$$\cFb\in D^b_{(\cS-)c}(X;R)\quad \text{ with}\quad 
\chi_{stalk}(\cFb)=0\in CF_{(\cS)}(X)$$
 are {\it preserved} by all Grothendieck functors like pullback $f^*$, direct image $Rf_!=Rf_*$ for a proper morphism, and the (external) tensor products 
$\otimes$ resp., $\boxtimes$ (which are exact in the case $R$ a field as considered here). So they induce similar transformations of constructible functions $CF_{(\cS)}(X)$ via the surjection
$$\chi_{stalk}\colon  \: K_0\left(D^b_{(\cS)-c}(X;R)\right)\to CF_{(\cS)}(X)\: .$$
Similarly for the nearby and vanishing cycles functors $\psi_f, \varphi_f$ of Deligne, as studied in Section \ref{sec:nvc}. So it is natural to ask for a direct description of the corresponding transformations of characteristic cycles in the embedded context, compatible with the isomorphisms
$$\begin{CD}
CC\colon  CF_{\cS}(X) @> \sim >> L(\cS,T^*M) \quad \text{and} \quad CC\colon  CF(X)@> \sim >> L(T^*M|_X) \:.
\end{CD}$$
Whenever possible, we formulate the corresponding results for the refined context of fixed Whitney stratifications.

\bex[Characteristic cycle of external products] \label{CC-external} \index{external product}
Let $X_i\hookrightarrow M_i$ be a closed embedding of the complex algebraic (or analytic) variety $X_i$ into a complex algebraic (or analytic) manifold $M_i$  (for $i=1,2$). Assume $\cS_i$ is a  Whitney stratification of $X_i$, with $\alpha_i\in CF_{\cS_i}(X_i)$ (for $i=1,2$).
Then $\alpha_1\boxtimes \alpha_2\in CF_{\cS_1\times \cS_2}(X_1\times X_2)$ is defined by 
$$\alpha_1\boxtimes \alpha_2 (x_1,x_2):=\alpha_1(x_1)\cdot \alpha_2(x_2) \quad \text{for $(x_1,x_2) \in X_1\times X_2$,}$$
with
\be
CC(\alpha_1\boxtimes \alpha_2)= CC(\alpha_1)\boxtimes CC(\alpha_2) \in L(\cS_1\times \cS_2,T^*M_1\times T^*M_2)\:.
\ee
Here, we use the identification $T^*(M_1\times M_2)=T^*M_1\times T^*M_2$, and this result follows directly from 
the {\it product formula for normal Morse data} as in Example \ref{NMDproduct}.
Alternatively, it can be deduced from the \index{multiplicativity}
{\it multiplicativity of the local Euler obstruction functions} (as stated in \cite[p. 426]{Mch}):
$$Eu_{X_1\times X_2}=Eu_{X_1}\boxtimes Eu_{X_2}$$
in case the $X_i$ are pure-dimensional ($i=1,2$).
\eex

To formulate the results for suitable  pullbacks $f^*$ and direct images $f_*$ for a  morphism $f\colon  M\to N$ of complex (algebraic) manifolds,
we need to compare the corresponding cotangent bundles as in the following commutative diagram (whose right square is cartesian):

\be \label{diag-coT} \begin{CD}
T^{*}M @< ^tf' << f^{*}(T^{*}N) @> f_{\pi}>> T^{*}N \\
@VV \pi_{M} V  @VV \pi V  @VV \pi_{N} V \\
M @= M @> f>> N \:.
\end{CD} \ee
Here $^tf'$ is the dual of the differential of $f$, with $f_{\pi}$ induced by base change.
We first consider the case of direct images. Let $X\hookrightarrow M$, resp., $Y\hookrightarrow N$ be closed complex algebraic (or analytic) subvarieties endowed with Whitney stratifications $\cS$ of $X$, resp., $\cT$ of $Y$ in such a way that $f\colon  X\to Y$ is a \index{stratified submersion}
{\it proper stratified submersion}.
Then also $f_{\pi}\colon  f^{*}(T^{*}N)|_X\to T^{*}N|_Y$ is {\it proper} by base change. Moreover
$$f_{\pi}\left( ^{t}f'^{-1}\left(T^*_{\cS}M\right)\right) \subset T^*_{\cT}Y \:,$$
since $f\colon  X\to Y$ is a stratified submersion. So on the level of {\it Lagrangian cycles} one gets an induced group homomorphism
\be
{f_{\pi *}}^{t}f'^{*}\colon  L(\cS,T^*M) \to L(\cT,T^*N) \:,
\ee
with the pullback $^tf'^{*}$ defined via Poincar\'{e} duality:
$$\begin{CD}
 H^{2m}_{T^*_{\cS}M}(T^*M) @> ^tf'^{*}  >>
H^{2m}_{ ^{t}f'^{-1}(T^*_{\cS}M)}(f^{*}(T^{*}N),\bZ)\\
@V PD V \wr V @V PD V \wr V \\
H^{BM}_{2m}(T^*_{\cS}M,\bZ) @> ^tf'^{*}  >> H^{BM}_{2n}(^{t}f'^{-1}(T^*_{\cS}M),\bZ) \:,
\end{CD}$$
with $m=\dim M, n=\dim N$ and $T^*M, T^*N, f^*(T^*N)$ oriented as complex manifolds.\\

Similarly, the pushforward of constructible functions fits into a commutative diagram
$$\begin{CD}
K_0\left(D^b_{\cS-c}(X;R)\right) @> f_* >> K_0\left(D^b_{\cT-c}(Y;R)\right)\\
@V \chi_{stalk} VV  @ VV \chi_{stalk} V \\
CF_{\cS}(X) @> f_* >> CF_{\cT}(Y) \:,
\end{CD}$$
with 
\be
f_*(\alpha)(y):= \int_{X\cap \{f=y\}} \alpha d \chi \quad \text{for $y\in Y$ and $\alpha\in CF_{\cS}(X)$.}
\ee

And these pushforwards of Lagrangian cycles and constructible functions are compatible with the characteristic cycle map
(see, e.g., \cite[Section 4.6]{Scf} for the following proof).

\bp\label{CC-pushdown} Let $\alpha\in CF_{\cS}(X)$ be given. Then
$$CC(f_*\alpha) = {f_{\pi *}}^{t}f'^{*}CC(\alpha) \in L(\cT,T^*N) \:.$$
\ep

\begin{proof} Let $\beta\in CF_{\cT}(Y)$ be defined by $CC(\beta)={f_{\pi *}}^{t}f'^{*}CC(\alpha) \in L(\cT,T^*N)$.
Then we need to show $\beta=f_*\alpha$. For this we calculate $\beta(y)$ 
for $y\in Y$ as in Corollary \ref{cor:invCC} via
$$\beta(y):=  \sharp_{dr_{y}} (\: {f_{\pi *}}^{t}f'^{*}CC(\alpha)\cap [dr(\{r<\epsilon\})]\:) \:,$$
with $r\colon  (N,y)\simeq (\bC^n,0)\to [0,\infty[\;$ given by $r(z)=\sum_{i=1}^n z_i\cdot \bar{z}_i$ the distance function to $y$ (in local coordinates) and
$0<\epsilon$ small enough. By the {\it projection formula} one gets
$$\beta(y)= \sharp (\: CC(\alpha) \cap [d(r\circ f)(\{r\circ f<\epsilon\})]\:) \:.$$
And the last intersection number is by the {\it local index formula} (\ref{loc-intersec}) given by
$$\sharp (\: CC(\alpha) \cap [d(r\circ f)(\{r\circ f<\epsilon\})]\:) = \int_{X\cap\{ f=y\}} \,\alpha \: d\chi = f_*\alpha(y) \:,$$
thus completing the proof.
\end{proof}

\bex Let $Y=N=\{pt\}$ be a point with $f$ a constant map, so that $f^*(T^*N)=T^*_MM$ is the zero section of $T^*M$.
Then we recover the \index{global index formula}{\it global index formula} (\ref{eq:PHopf}):
$$
\int_X \, \alpha \: d\chi   = \sharp(\:CC(\alpha) \cap [T^{*}_{M}M]\:) \:.
$$
\eex

To define similarly a pullback of Lagrangian cycles, we have to go in diagram (\ref{diag-coT}) into the opposite direction, and need the following  properness condition for $^{t}f'\colon   f^{*}(T^{*}N) \to T^{*}M$ (see, e.g., \cite[Lemma 4.3.1]{Sc}).

\bd Let $C\subset T^*N$ be a closed $\C^*$-conic subset. Then $f\colon M\to N$ is \index{non-characteristic} {\it non-characteristic} with respect to $C$, if one of the following two equivalent conditions holds:
\begin{enumerate}
\item $f_{\pi}^{-1}(C) \cap kern(^{t}f') \subset f^*(T^{*}_{N}N)$,
with $kern(^{t}f')$ the corresponding kernel bundle, and $f^*(T^{*}_{N}N)$ 
the zero-section of $f^{*}(T^{*}N)$.
\item The map $\,^{t}f'\colon  f_{\pi}^{-1}(C)\to T^{*}M$ is {\it proper} and therefore finite.
\end{enumerate}
\ed

\bex If $f\colon M\to N$ is a {\it submersion}, then $f$ is {\it non-characteristic} with respect to any closed $\C^*$-conic subset $C\subset T^*N$,
since then $kern(^{t}f') = f^*(T^{*}_{N}N)$ is the zero-section of $f^{*}(T^{*}N)$.
\eex

More generally, {\it transversality} with respect to a Whitney stratification $\cT$ of a closed complex algebraic (or analytic) subset
$Y\hookrightarrow N$ can be characterized as follows (see, e.g., \cite[Example 4.3.2]{Sc} and \cite[Definition 4.1.5]{KS}).

\bex\label{ex.transv}
 The morphism $f\colon M\to N$ is \index{transversal morphism} {\it transversal} to $\cT$, i.e., $f$ is transversal to all strata $T\in \cT$, if and only if $f$ is {\it non-characteristic} with respect to
the closed $\C^*$-conic subset $T^*_{\cT}N\subset T^*N$. In this case, $X:=f^{-1}(Y)$ gets an induced Whitney stratification $\cS=f^{-1}\cT$ with strata $S$ given by the connected components of the locally closed complex submanifolds $f^{-1}(T)\subset M$ (for $T\in \cT$).
Here, the codimension $\codim f^{-1}(T)=\codim T$ is preserved, i.e., 
$$\dim f^{-1}(T)=  \dim T + \dim M - \dim N \quad \text{for all $T\in \cT$.}$$
Moreover, $\,^{t}f'\colon  f_{\pi}^{-1}(T^*_{T}N)\to T^*M$ is injective with image $T^*_{f^{-1}(T)}M$. By Poincar\'{e} duality one gets an induced
pullback map of Lagrangian cycles
\be
\,^{t}f'_* f_{\pi}^{*}\colon  L(\cT,T^*N) \to L(\cS,T^*M), 
\ee
with 
$$\,^{t}f'_* f_{\pi}^{*}\left(\left[\overline{T^*_TN}\right]\right)=\left[\overline{T^*_{f^{-1}(T)}M}\right] \quad \text{ for all $T\in \cT$.}$$
Moreover, normal slices to $f^{-1}(T)\subset M $ and $T\subset N$ get identified via $f$, so that the corresponding Euler characteristics of 
normal Morse data do not change
under pullback of constructible functions given by $f^*\alpha :=\alpha\circ f$ for $\alpha\in CF_{\cT}(Y)$.
In particular, the characteristic cycle map $CC$ commutes with this pullback only up to a sign:
$$(-1)^{\dim M -\dim N}\cdot CC(f^*\alpha)= \,^{t}f'_* f_{\pi}^{*} CC(\alpha)  \quad \text{for all $\alpha\in CF_{\cT}(Y)$.}$$
For example, if $Y$ (and hence also $X=f^{-1}(Y)$) is pure-dimensional, then one gets for such a transversal map
$$f^*(Eu_Y)=Eu_X\:.$$
\eex

More generally one has the following result (see, e.g., \cite[Theorem 3.1]{Sct}).

\bt\label{non-cc}
Let  $f\colon  M\to N$ be a morphism of complex algebraic (or analytic)  manifolds of dimension $m=\dim M, n=\dim N$,
with  $Y\subset N$ a closed complex algebraic (or analytic) subvariety and  $X:=f^{-1}(Y)\subset M$.
Assume that $f$ is \index{non-characteristic}
{\it non-characteristic} with respect to the {\it support}
$C:=supp(CC(\alpha))\subset T^*N|_Y$ of the characteristic cycle $CC(\alpha)$ of a constructible
function $\alpha\in CF(Y)$. Then $C':=\,^{t}f'\left( f_{\pi}^{-1}(C)\right)$ is pure $m$-dimensional,
with
\begin{equation}
\,^{t}f'_* f_{\pi}^{*}(CC(\alpha)) = (-1)^{m-n}\cdot CC(f^*(\alpha))\:.
\end{equation}
In particular,  the left hand side is a Lagrangian cycle in $T^*M|_X$.
\et

The following Example \ref{1350} has nice applications in {\it geometric representation theory} for the \index{Weyl group} {\it Weyl group} $\bW$
of a connected semisimple complex Lie group $G$  (see, e.g., \cite{Gi2}).
Here, $M_1=M_2=G/B$ is the \index{flag manifold} {\it Flag manifold} of $G$ (with $B\subset G$ a Borel subgroup), and the Whitney stratification $\cS$ of $G/B\times G/B$ 
is given by the finitely many $G$-orbits $S_w$ of the diagonal $G$-action (indexed by $w\in \bW$). Then any $\alpha\in CF_{\cS}(G/B\times G/B)$ satisfies the following assumption for 
$$C:=\bigcup_{w\in \bW} T^*_{S_w}(G/B\times G/B) \subset T^*(G/B\times G/B)$$ 
the corresponding \index{Steinberg variety}  {\it Steinberg variety}.

\bex\label{1350} Let $M_i$ be three complex algebraic (or analytic) complex manifolds of dimension $m_i=\dim M_i$ ($i=1,2,3$). 
Assume $\alpha\in CF(M_1\times M_2)$ satisfies 
one of the following two equivalent conditions for a closed $\C^*$-conic subset $C \subset T^*(M_1\times M_2)$
with $supp(CC(\alpha))\subset C$:
\begin{enumerate}
\item  The projection $T^*(M_1\times M_2)=T^*M_1\times T^*M_2\to T^*M_1\times M_2$ restricted to $C$ is proper
and therefore finite,
\item $C \cap M_1\times T^*M_2$ is contained in the zero section of $T^*(M_1\times M_2)$.
\end{enumerate}
Then the embedding $d\colon  M_1\times M_2\times M_3 \to (M_1\times M_2)\times (M_2\times M_3)$ induced by the diagonal embedding
$M_2\to M_2\times M_2$ is \index{non-characteristic}
{\it non-characteristic} with respect to $supp(CC(\alpha\boxtimes \beta))$ for {\it any} $\beta \in CF(M_2\times M_3)$,
so that
$$\,^{t}d'_* d_{\pi}^{*}(CC(\alpha\boxtimes \beta)) = (-1)^{m_2}\cdot CC(d^*(\alpha\boxtimes \beta))\:.$$
\eex

Another application of Theorem \ref{non-cc} is the following \index{intersection formula} {\it intersection formula}  (see, e.g., \cite[Corollary 3.1]{Sct}).
\bc\label{cor-int}
Let $M$ a complex algebraic (or analytic) manifold of dimension $m=\dim M$, with $\alpha, \beta\in CF(M)$ given constructible functions.
Assume that the diagonal embedding $d\colon  M\to M\times M$ is {\it non-characteristic} with respect
to $supp(CC(\alpha\boxtimes \beta))$, with  $supp(\alpha\cdot \beta)$ compact.

Then also  $supp(CC(\alpha)\cap CC(\beta))\subset T^*M$ is compact, with
\be\label{mic-int}
\int_M \alpha\cdot \beta \: d \chi = (-1)^m\cdot deg( CC(\alpha)\cap CC(\beta))\:.
\ee
\ec

\bex The assumption $d\colon  M\to M\times M$ is {\it non-characteristic} with respect
to $supp(CC(\alpha\boxtimes \beta))$ for $\alpha, \beta\in CF(M)$  holds in the following cases:
\begin{enumerate}
\item $\alpha$, resp., $\beta$ is constructible with respect to a Whitney stratification $\cS$, resp., $\cT$ of $M$, with $\cS$ and $\cT$ intersecting transversally (i.e., all strata $S\in \cS$ and $T\in \cT$ intersect transversally in $M$),
so that $d\colon  M\to M\times M$ is transversal to the product stratification $\cS\times \cT$ of $M\times M$.
\item  $\alpha$ and $\beta$ are {\it splayed} \index{splayed}
 in the sense of \cite[Definition 2.1]{Sc}, i.e., for any $p\in M$ there are  locally analytic isomorphisms $M=V_1\times V_2$
 of analytic manifolds so that $\alpha=\pi_1^*(\alpha')$ and $\beta=\pi_2^*(\beta')$ for some $\alpha'\in CF(V_1)$ and 
 $\beta'\in C F(V_2)$, with $\pi_i\colon  V_1\times V_2\to V_i$ the projection ($i=1,2$). For $\alpha=1_X$ and $\beta=1 _Y$, with $X,Y\subset M$ closed complex algebraic (or analytic) subvarieties, this just means that $X$ and $Y$ are {\it splayed} in the sense of Aluffi-Faber \cite{AF, AF2}.
\end{enumerate}
\eex

We next explain the relation between \index{specialization}  {\it specialization} of Lagrangian cycles and \index{nearby cycle} {\it nearby cycles} for constructible functions
(see, e.g., \cite{BMM, Gi1, Sab1}).

Let $f\colon N\to \bC$ be a submersion of complex algebraic (or analytic) manifolds, 
with $k\colon M:=\{f=0\}\hookrightarrow N$ the inclusion of a smooth hypersurface with open complement $U=\{f\neq0\}$. Consider the exact sequence of vector bundles on $N$
$$0\to \langle df \rangle= kern(p)\to T^*N \to T_f^*\to 0\:,$$
with the projection $p\colon  T^*N\to T^*_f$ dual to the inclusion $T_f\to TN$ of the subvector bundle of tangents to the fibers of $f$.
Let $X\hookrightarrow N$ be a closed complex algebraic (or analytic) subvariety endowed with a Whitney stratification $\cS$ so that
$X\cap\{f=0\}=:X_0$ is a union of strata, with $\cS|_{X_0}$, resp., $\cS|_U$ the induced Whitney stratification of $X_0$, resp., $X\cap U$.
After shrinking to an open neighborhood, we can and will assume by the {\it curve selection lemma} that $f|_S$ is a submersion for all strata
$S\in \cS|_U$. Finally, we assume that $\cS$ satisfies the following \index{$a_f$-condition} {\it $a_f$-condition of Thom}:
\be\label{af}
\begin{cases} &\text{If $x_n\in S'$, for $S'\in \cS|_U$, is a sequence converging to $x\in S\in \cS|_{X_0}$,}\\
& \text{such that $kern(df_{x_i}|_{T_{x_i}S'}) = T_{x_i}(S'\cap \{f=f(x_i)\})$} \\
& \text{converges to some limiting plane $\tau$, then $T_xS\subset \tau$.}
\end{cases}\ee

 By a classical  theorem of Hironaka \cite[Corollary 1 to Theorem 2, p.248] {Hir} (and see also \cite[Corollary 1.3.5.1]{LeTe}),
one can always refine a given stratification so that it satisfies the 
$a_f$-condition. In fact, by a more recent result \cite[Theorem  4.2.1]{BMM}, this $a_f$-condition is true for any Whitney stratification $\cS$ as above
(with $X_0$ a union of strata).\\

Then $T^*_{S'}N\cap kern(p)$ is contained in the zero-section $T_N^*N|_U$ of $T^*N|_U$,
since $f\colon S'\to \bC$ is a submersion for $S'\in \cS|_U$.
Therefore $p'=p\colon T^*_{S'}N\to T^*_f|_S$ is a proper injection with image
$p'(T^*_{S'}N)=:T^*_{S'}f$, 
 the {\it relative conormal space} of $f|_{S'}$.
Here we consider the commutative diagram 

\be\begin{CD}
T^*_{S'}N @>>> T^*N |_U @> j >> T^*N\\
@V p' VV @V p' VV  @VV p V \\
T^*_{S'}f @>>> T_f^*|_U @> j' >> T^*_f @< i << T^*M=T^*_f|_M\:.
\end{CD}\ee

Note that by the $a_f$-condition the closure $\overline{T^*_{S'}f}\subset T^*_f$ for a stratum $S'\in \cS|_U$ is contained in the closed subset
$$T^*_{\cS}f:=\bigcup_{S\in \cS} T^*_Sf \: \hookrightarrow T^*_f\:,$$
with $T^*_Sf:=T^*_SM$ for $S\in \cS|_{X_0}$ and $p'\colon  T^*_{\cS}N|_U\to T^*_f|_U$ proper.
Then the \index{specialization} specialization of Lagrangian cycles $$sp\colon  L(\cS,T^*N)\to L(\cS|_{X_0},T^*M)$$ is defined as the composition of the following homomorphisms
(with $n=\dim N =\dim M +1$):
\be\begin{CD}
H^{BM}_{2n}(T^*_{\cS}N,\bZ) @> j^* >> H^{BM}_{2n}(T^*_{\cS}N|_U,\bZ )@> p'_* >> H^{BM}_{2n}(T^*_{\cS}f|_U,\bZ )\\
@V sp VV @. @A j'^* A \wr A \\
L(\cS|_{X_0},T^*M) @= H^{BM}_{2n-2}(T^*_{\cS|_{X_0}}M,\bZ ) @< i^* <<H^{BM}_{2n}(T^*_{\cS}f,\bZ )\:,
\end{CD}\ee
with the Gysin map $i^*$ again defined by Poincar\'{e} duality as the intersection with $[T^*M]$.
So $sp([T^*_SN])=0$ for $S\in \cS|_{X_0}$ and
$$sp([\overline{T^*_{S'}N}])= \left[ \overline{T^*_{S'}f}\right] \cap [T^*M] \quad \text{for $S'\in \cS|_U$.}$$

The {\it nearby cycles} \index{nearby cycle} of constructible functions used in the next Theorem are induced from Deligne's nearby cycle functor $\psi_f$
for constructible sheaf complexes as discussed later on  in Section \ref{sec:nvc}, fitting into a commutative diagram
\be\begin{CD}
K_0(D^b_{\cS-c}(X;R)) @> \psi_f >> K_0(D^b_{\cS|_{X_0}-c}(X_0;R))\\
@V \chi_{stalk} VV @VV \chi_{stalk} V \\
CF_{\cS}(X) @ > \psi_f >> CF_{\cS|_{X_0}}(X_0) \:.
\end{CD}\ee

With these notations, we can now show the following result (see, e.g., \cite[Theorem 4.3]{Sab1} and compare with \cite{BMM, Gi1}),
which is also used in \cite{Sct} in the proof of Theorem \ref{non-cc} above.
\bt \label{special} Let $\alpha\in CF_{\cS}(X)$ be given. Then
\be
sp(CC(\alpha))=CC(-\psi_f(\alpha)) \in L(\cS|_{X_0},T^*M)\:,
\ee
with $\psi_f(\alpha)(x):= \int_{M_{f|_X,x}} \alpha d \chi$, for $M_{f|_X,x}$ a local Milnor fiber \index{Milnor fiber} of $f|_X$ at $x\in X_0$.
\et 

\begin{proof} Let $\beta\in CF_{\cS|_{X_0}}(X_0)$ be defined by $CC(\beta)= sp(CC(\alpha)) \in L(\cS|_{X_0},T^*M)$.
Then we need to show that $\beta=-\psi_f(\alpha)$. For this, we calculate $\beta(x)$ 
for $x\in X_0$ as in Corollary \ref{cor:invCC} via
$$\beta(x):=  \sharp_{dr'_{x}} (\: sp(CC(\alpha))\cap [dr'(\{r'<\epsilon\})]\:) \:,$$
with $r':=r|_M$, for $r\colon  (N,x)\simeq (\bC^n,0)\to [0,\infty[\;$ given by $r(z)=\sum_{i=1}^n z_i\cdot \bar{z}_i$ the distance function to $0$ (in local coordinates of $N$), and
$0<\epsilon$ small enough. 
But this intersection number is locally constant in the family $T^*_f|_{\{f=w\}}$ for $|w|\ll 1$ small compared to $\epsilon$,
since $\{r=\epsilon\}$ is transversal to all $S\cap \{f=w\}$, for $S\in \cS$ and  $|w|\ll 1$ small, by the $a_f$-condition as well as the Whitney condition
for $\cS|_{X_0}$.
So instead of specializing at $w=0$, we can do this at a small  stratified regular value $w\neq 0$ and use Example \ref{ex.transv}:
$$\beta(x)= - \:\sharp (\: CC(\alpha|_{\{f=w\}}) \cap  [dr''(\{r''<\epsilon\})]\:) \:,$$
with $r'':=r|_{\{f=w\}}$. Note that  the sign comes from $\dim N -\dim \{f=w\}=1$, i.e., we intersect transversally with a submanifold of codimension $1$.
And the last intersection number is by Equation  (\ref{eq:intersecCF1}) given by
$$\: \sharp (\: CC(\alpha|_{\{f=w\}}) \cap  [dr''(\{r''<\epsilon\})]\:)= \:\int_{M_{f|_X,x}} \,\alpha \: d\chi =  \psi_f(\alpha)(x) \:,$$
with $M_{f|_X,x}:=X\cap \{f=w\}\cap \{r\leq \epsilon'\}$, for $0<|w|\ll \epsilon'<\epsilon$,  a local Milnor fiber of $f|_X$ at $x\in X_0$.
\end{proof}

Note that the $a_f$-condition is also needed to have such a {\it local Milnor fibration} with \index{Milnor fiber}
Milnor fiber $M_{f,x}$ for a general holomorphic function germ
$f\colon  (X,x)\to (\bC,0)$ on a singular complex analytic variety $X$  (see, e.g., \cite[Example 1.3.3]{Sc} and \cite{Le2}).
The definition of {\it nearby cycles} of constructible functions in terms of (weighted) Euler characteristics of these {\it local Milnor fibers} goes back to Verdier 
\cite{Ve3}.\\

As an example, let us explain how Theorem \ref{special} implies Theorem \ref{non-cc} for the inclusion $k\colon M=\{f=0\}\hookrightarrow N$ of a global smooth hypersurface (i.e., of codimension one), with $f\colon N\to \bC$ a submersion and the notations as before. Consider the following cartesian diagram:

\be\begin{CD}
T^*N|_M @>  k_{\pi} >> T^*N  @< j << T^*N|_U\\
@V \,^{t}k' VV @V p VV  @VV p' V \\
T^*M @> i  >> T_f^* @< j' << T^*_f|_U \:.
\end{CD}\ee

Assume  that $k\colon  M\hookrightarrow N$ is {\it non-characteristic} with respect to $supp(CC(\alpha))\subset T^*N$ for a given $\alpha\in CF(N)$, so that
$ \,^{t}k' \colon   k_{\pi}^{-1}(supp(CC(\alpha))) \to T^*M$ is proper. After shrinking $N$ we can then assume that also
$p\colon  supp(CC(\alpha))\to T_f^*$ is proper.
Then one gets by base change
$$j'^*p_*CC(\alpha)= p'_*j^*CC(\alpha)\:,$$
and
\begin{eqnarray*}
\,^{t}k'_*k_{\pi}^*CC(\alpha)  &=& i^*p_*CC(\alpha)= sp(CC(\alpha))\\
&=&  CC(-\psi_f(\alpha)) = - CC(k^*\alpha)\:.
\end{eqnarray*}
Here, the fact that  $\psi_f(\alpha)=k^*\alpha$ follows from (\ref{no-phi}) below, since we have 
$0=\psi_f(\alpha)(x)-k^*\alpha(x)= \psi_f(\alpha)(x)-\alpha(x)$ for all $x\in M$ by the
 non-characteristic assumption, which gives $df_x\not\in supp(CC(\alpha))$ for all $x\in M$.\\

 We finish this section with citing from \cite[Corollary 0.3]{ScL} (and compare with \cite[Theorem 4.5]{Sab1})
 a nice intersection formula related to the vanishing cycle functor $\varphi_f$ introduced in Section \ref{sec:nvc}.
For a description of the {\it characteristic cycle of vanishing cycles} we refer to \cite[Theorem 2.10]{Ma}.

\bt[Global intersection formula for vanishing cycles] \label{thm:intvan}
Let $M$ be a complex algebraic (or analytic) manifold and 
$f\colon M\to \bC$ an algebraic (or  holomorphic) function, with $\cFb\in D^b_c(M;R)$, resp., $\alpha \in CF(M)$  be given.
Suppose that the intersection
of $df(M) $ and the support 
of the characteristic cycle of $\cFb$, resp., $\alpha$ is 
contained in a
{\it compact} complex algebraic (or analytic) subset $I\subset T^{*}M$, with $K:=\pi(I)\subset \{f=0\}$
for $\pi\colon  T^*M\to M$ the projection.
Then one has
\be \label{eq:intvan}
\chi\bigl(R\Gamma(K,\varphi_{f}[-1]\;\cFb)\bigr) = 
\sharp\bigl(\;  [CC(\cFb)] \cap [df(M)]\; \bigr) \:,
\ee
resp.,
\be \label{eq:intvan2}
-\int_K \varphi_{f}(\alpha)\:d \chi=
- \chi(K,\varphi_{f}(\alpha))  = 
\sharp\bigl(\; [CC(\alpha)] \cap [df(M)] \; \bigr) \:.
\ee
\et

\bex Let  $f\colon  M\to \{0\}$  be the constant zero function so that $df(M)=T^*_MM$ is the zero section of $T^*M$,
with $-\varphi_{f}(\alpha)=\alpha$ for all $\alpha\in CF(M)$. Assume $M=:K$ is compact.
Then we recover the \index{global index formula} {\it global index formula} (\ref{eq:PHopf}):
$$
\int_X \, \alpha \: d\chi   = \sharp(\:CC(\alpha) \cap [T^{*}_{M}M]\:) \:.
$$
\eex

In the special case $I=\{\omega\}$ given by a point 
$\omega\in T^{*}M$,  we 
get back by Theorem \ref{thm:intvan}
a formula conjectured by Deligne (with $x:=\pi(\omega)$ and $K:=\{x\}$):
\be \label{eq:int5}
\chi\bigl((\phi_{f}[-1]\;\cFb)_{x}\bigr) =
\sharp_{df_{x}}\bigl(\; [CC(\cFb] \cap [df(M)]\; \bigr)
\ee
and
\be \label{eq:int5CF}
- \varphi_{f}(\alpha)(x)=
\alpha(x) - \int_{M_{f;x}}\, \alpha \,d\chi = 
\sharp_{df_{x}}\bigl(\;  [CC(\alpha)] \cap [df(M)]\; \bigr) \:,
\ee
with $M_{f,x}$ the \index{Milnor fiber}
{\it local Milnor fiber} of $f$ in $x$. In particular,  for $df_x\not\in supp(CC(\alpha))$, we get
\be\label{no-phi}
- \varphi_{f}(\alpha)(x)= \alpha(x)-\psi_f(\alpha)(x)=
\alpha(x) - \int_{M_{f;x}}\, \alpha \,d\chi = 0 \:.
\ee

\bex
Let $f\colon  (M,x)=(\bC^{n+1},0)\to (\bC,0)$ be a holomorphic function germ with an \index{isolated singularity}
 {\it isolated singularity} in $0\in \bC^{n+1}$.
Then the local Milnor fiber $M_{f,0}$ has the homotopy type of a wedge of a finite number of $n$-spheres (for $n\geq 1$). The \index{Milnor number} {\it Milnor number} $\mu(f)$ is the number of these $n$-spheres, so that
$$\chi(M_{f,x},\bQ)=1 + (-1)^n\cdot \mu(f)\:.$$ 
For $n=0$, the local Milnor fiber $M_{f,0}$ consists of $1+\mu(f)$ points.
If we apply the above formula to the perverse sheaf $\cFb=\bQ_M[n+1]$, resp., $\alpha=(-1)^{n+1}\cdot 1_M$,
with $CC(\bQ_M[n+1])=CC((-1)^{n+1}\cdot 1_M)=T^*_MM$ the zero section of $T^*M$, we recover the formula
$$\mu(f)= \sharp_{df_{x}}\bigl(\; [T^*_MM] \cap [df(M)] \; \bigr) $$
for the Milnor number, i.e., 
$$\mu(f)=\dim _{\bC} \frac{\left(\cO_{\bC^{n+1}}\right)_0}{\left(\nabla f\right)_0}\:$$
for the Milnor number in terms of the \index{Jacobian ideal}
Jacobian ideal $\left(\nabla f\right)_0\subset \left(\cO_{\bC^{n+1}}\right)_0 $.
\eex

The {\it local intersection formula for the vanishing cycle functor} (\ref{eq:intvan})  is due 
to Dubson \cite[Theorem 1]{Du3}, Ginsburg \cite[Proposition 7.7.1]{Gi1}, 
L\^{e} \cite[Theorem 4.1.2]{Le5} and Sabbah \cite[Theorem 4.5]{Sab1}. 
For a discussion of the history of this {\it local intersection formula for vanishing cycles}
we recommend the paper \cite{Le5}.

But most of these
references are in the language of  {\it holonomic D-modules} or {\it perverse sheaves}. So the assumption on the 
intersection for a  constructible complex of sheaves $\cFb$
corresponds to an assumption on the {\it micro-support} $\mu supp(\cFb)$.
To be able to state the result also for constructible functions, it is important to work with the weaker assumption about the {\it support of 
the characteristic cycle} $supp(CC(\cFb)) \subset \mu supp(\cFb)$.\\

Let us finally point out that MacPherson's theory \cite {Mch} of \index{Chern classes}
{\it Chern classes $c_*$ of singular varieties and constructible functions} in the embedded context can easily be recovered and improved form the {\it functorial theory of characteristic cycles} as presented above, e.g., like  the following results about $c_*$:
\begin{enumerate}
\item Functoriality of $c_*$ for proper morphism via Proposition \ref{CC-pushdown}.
\item Multiplicativity of $c_*$ with respect to external products via Example \ref{CC-external}.
\item Specialization of $c_*$ in one parameter families via Theorem \ref{special}.
\item A  Verdier-Riemann-Roch Theorem for $c_*$ with respect to non-characteristic pullbacks (e.g., submersions or transversal pullbacks)
via Theorem \ref{non-cc}.
\item An intersection formula for $c_*$ via Theorem \ref{non-cc} in the context of Example \ref{mic-int}.
\end{enumerate}
In fact, the most recent approach of \cite{AMSS} deduces $c_*(\alpha)$ for $\alpha\in CF(X)$ from  $CC(\alpha)$, viewed as a $\bC^*$-invariant
cycle, via intersection with the zero-section $T^*_MM$ in $\C^*$-equivariant Borel-Moore homology.
For further reading also on relations between {\it MacPherson Chern classes} $c_*$ and {\it characteristic cycles} we refer to, e.g., 
\cite{AF, AF2, AMSS, AMSS2, CMSS, Gi1, Gi2, Sab1, Scf, Sct, ST, Ve3}.


\subsection{Vanishing results}
Let us now explain other applications of the {\it stratified Morse theory for constructible sheaves} and of Theorem \ref{thm:plurisub},
to {\it vanishing and weak Lefschetz theorems} in the complex algebraic and analytic context.
In the following we are back to the general context that the base ring $R$ is commutative and noetherian, of finite global dimension, and we work with (weakly) constructible complexes of sheaves of $R$-modules.\\

First we consider  the complex algebraic context with an {\it affine} complex algebraic variety $X\hookrightarrow \bC^n$ and
the {\it strongly plurisubharmonic and semi-algebraic} distance function 
$$r\colon  \bC^{n}\to \bR^{\geq 0}\,; \: r(z):=\sum \nolimits_{i=1}^{n} \: 
z_{i}\bar{z}_{i} \:.$$
If  $\cFb\in D^b_{\cS-wc}(X;R)$ is $\cS$-weakly constructible with respect to a {\it complex algebraic} Whitney stratification $\cS$ of $X$,
then the proper semi-algebraic distance function $r$ has only {\it finitely} many stratified critical values so that one gets 
as in the proof of Theorem \ref{thm:intersec} and Example \ref{index-affine}:
$$R\Gamma(X\cap \{r \leq b\},\cFb)\simeq  R\Gamma(X,\cFb) \quad \text{for $b>0$ large enough,}$$
and
$$R\Gamma_{c}(X\cap \{r < b\},\cFb)\simeq  R\Gamma_{c}(X,\cFb) \quad \text{for $b>0$ large enough.}$$

Then Theorem \ref{thm:plurisub} implies for $f:=r$ (with $q=0$, since $r$ is strongly plurisubharmonic)
directly  the following  important {\it Artin-Grothendieck type result} (see also \cite[Corollary 6.0.4]{Sc} for a more general version):

\bt[Artin vanishing Theorem] \label{Artinvan} \index{Artin vanishing}
Let $X$ be an {\it affine} complex algebraic variety,
with  $\cFb\in D^b_{(w)c}(X;R)$ a {\it complex algebraically} (weakly) constructible sheaf complex
(respectively, complex algebraically constructible in case of a Dedekind domain $R$, if we want to use the dual perverse t-structure).
Then:
\begin{enumerate}
\item $\cFb\in \,^{p^{(+)}}D^{\leq n}(X;R)  \Rightarrow 
R\Gamma(X,\cFb) \in  \,^{p^{(+)}}D^{\leq n}(\{pt\};R) \subset D^b_{(w)c}(\{pt\};R)   $,
\item $\cFb\in \,^{p^{(+)}}D^{\geq n}(X;R)  \Rightarrow 
R\Gamma_c(X,\cFb) \in  \,^{p^{(+)}}D^{\geq n}(\{pt\};R)\subset D^b_{(w)c}(\{pt\};R)  $.
\end{enumerate}
\et

\br\label{opensemia}
The same proof as indicated above gives the {\it Artin vanishing Theorem} more generally for an {\it open semi-algebraic} subset $j\colon  X\hookrightarrow X'$
of a complex algebraic variety $X'$, if it has a proper {\it strongly plurisubharmonic and semi-algebraic}  function 
$r\colon  X\to \bR^{\geq 0}$, and we consider only (weakly) constructible sheaves of the form $\cFb=j^* \cGb$ for some algebraically (weakly) constructible
$\cGb\in D^b_{(w)c}(X',R)$.
\er

Compare  with \cite[Section 4]{BBD}, especially 
\cite[Theorem 4.1.1]{BBD}),
for the corresponding relative counterpart for an affine morphism 
in the context
of the {\it perverse t-structure}  in $l$-adic
cohomology. The relative version of Theorem~\ref{Artinvan} 
for an {\it affine complex algebraic morphism} is also true, as we will see later on.
For now, let us only illustrate this Theorem~\ref{Artinvan}  by the following example
(see, e.g., \cite[Example 6.0.4]{Sc}).

\bex[Weak Lefschetz theorem for singular spaces] \label{ex:weaklef} \index{weak Lefschetz theorem}
Let $X$ be a closed algebraic subvariety of the complex projective space, with $H$ a hyperplane.
Consider the open inclusion $j\colon  U:=X\backslash H \hookrightarrow X$, with $i\colon  X\cap H\hookrightarrow X$ the corresponding closed inclusion.
Assume $\cFb\in D^b_{(w)c}(X;R)$ is {\it complex algebraically} (weakly) constructible
(respectively, complex algebraically constructible in case of a Dedekind domain $R$, if we want to use the dual perverse t-structure).
Then
$$j^{*}\cFb\in \,^{p^{(+)}}D^{\leq n}(U;R)  \Rightarrow$$
$$ R\Gamma(X,Rj_{*}j^{*}\cFb)
\simeq R\Gamma(U,j^{*}\cFb) \in  \,^{p^{(+)}}D^{\leq n}(\{pt\};R) \:,$$
and
$$j^{*}\cFb\in \,^{p^{(+)}}D^{\geq n}(U;R)  \Rightarrow$$
$$ R\Gamma(X,X\cap H,\cFb):=R\Gamma(X,Rj_{!}j^{*}\cFb)
\simeq R\Gamma_c(U,j^{*}\cFb) \in  \,^{p^{(+)}}D^{\geq n}(\{pt\};R) \:.$$
In particular:
\begin{enumerate}
\item The {\it relative homology} $$H_{k}(X,X\cap H; R)=H^{-k}(U; j^{*}\bD_X^{\bullet})=0 \quad \text{for $k< \r(U,R)$.}$$
\item If $R$ is a Dedekind domain, then the {\it relative cohomology}
$$H^{k}(X,X\cap H;R) \simeq H^k_c(U; R_U)= \begin{cases}
0 & \text{for $k< \r(U,R)$,}\\
\text{torsion-free} & \text{for $k= \r(U,R)$.}
\end{cases}$$
\item Suppose $X$ is purely $n$-dimensional, and $H$ is a {\it generic} hyperplane,
i.e., $H$ is transversal to a Whitney  stratification $\cS$ of $X$,  so that $i^*IC_X[-1]\simeq IC_{X\cap H}$ by Example \ref{IC-transv}.
Then the {\it relative intersection cohomology} \index{relative intersection cohomology}
$$IH^k(X, X\cap H; R) := H^{k-n}_c(U; j^{*}IC_X)=0 \quad \text{ for $k<n$.}$$ Moreover, $IH^n(X, X\cap H; R)$ is torsion-free in case $R$ is a Dedekind domain.
\end{enumerate}
\eex

Here, statements 1. and 2. show the role of the {\it rectified homological depth} \index{rectified homological depth} for
the ``weak Lefschetz theorem for singular spaces'' as conjectured by Grothendieck.
Compare also with \cite[Theorem 3.4.1]{HL} for the
corresponding homotopy result.

The ``weak Lefschetz theorem for intersection homology''
3. is due to \cite[Theorem 7.1]{GM}, and compare also with \cite[Theorem 6.10]{GMs}.
Here we only cite the following remark of Goresky-MacPherson from the beginning of
\cite[Section 6.10]{GMs}: ``\dots The following Lefschetz hyperplane theorem was
our original motivation for developing Morse theory on singular spaces \dots''.\\

The  {\it weak Lefschetz theorem} of Example~\ref{ex:weaklef} can be generalized
in many different directions. Especially, it is enough to assume that $X$ is
a {\it quasi-projective} \index{quasi-projective} 
algebraic subvariety $X=X'\backslash A$, with $A\subset X'\subset \bC P^{N}$
closed subvarieties, if $H$ is a {\it generic} hyperplane.
Just take a Whitney  stratification of $X'$ such that $A$ is a union of strata, and
$H$ is {\it transversal} to all strata. Then one can apply the following \index{base change isomorphism} {\it base change isomorphisms} (see, e.g., \cite[Lemma 6.0.5]{Sc}).

\bl \label{lem:basechange6}
Let $A\subset X'$ be closed analytic subvarieties of the complex manifold $M$.
Fix a Whitney stratification $\cS$ of $X'$, with $A$ a union of strata.
Suppose $H$ is closed complex submanifold of $M$, which is {\it transversal} to all strata,
and consider the following cartesian diagram of inclusions:
\begin{displaymath} \begin{CD}
X:= X'\backslash A @< j << U:=X\backslash H \\
@VV k' V @VV k V \\
X' @< j' << U':=X'\backslash H \:.
\end{CD} \end{displaymath}
Then  for $\cFb\in D^b_{\cS|_X-wc}(X;R)$ a weakly constructible complex with respect to the induced stratification $\cS|_X$ of $X$, one has natural isomorphisms:
\be \label{basechange6}
j'_{!}Rk_{*}j^{*}\cFb\simeq Rk'_{*}j_{!}j^{*}\cFb
\quad \text{and} \quad
Rj'_{*}k_{!}j^{*}\cFb \simeq k'_{!}Rj_{*}j^{*}\cFb \:.
\ee
If $X'$ is also {\it compact}, then one gets for $\cFb\in D^b_{\cS|_X-wc}(X;R)$:
$$R\Gamma_{c}(X,Rj_{*}j^{*}\cFb
\simeq R\Gamma(X',k'_{!}Rj_{*}j^{*}\cFb)
\simeq R\Gamma(U',k_{!}j^{*}\cFb) \:,$$
and
$$R\Gamma(X,X\cap H,\cFb) \simeq R\Gamma(X',Rk'_{*}j_{!}j^{*}\cFb)
\simeq 
R\Gamma_{c}(U',Rk_{*}j^{*}\cFb) \:. $$
\el

 Moreover,
$$j^{*}\cFb\in \,^{p^{(+)}}D^{\leq n}(U;R) \: \Rightarrow \: k_{!}j^{*}\cFb\in \,^{p^{(+)}}D^{\leq n}(U';R) \:,$$
and
$$j^{*}\cFb\in \,^{p^{(+)}}D^{\geq n}(U;R) \: \Rightarrow \: Rk_{*}j^{*}\cFb\in \,^{p^{(+)}}D^{\geq n}(U';R) \:.$$

Applying the Artin vanishing Theorem \ref{Artinvan} to these complexes on the right hand side, one gets that the  
{\it weak Lefschetz theorem} of Example~\ref{ex:weaklef} 
remains true for $X=X'\backslash A$  a quasi-projective variety, with $A\subset X'\subset \bC P^{N}$ closed subvarieties, if $H$ is a {\it generic} hyperplane as discussed above (with $U'\hookrightarrow \bC P^{N}\backslash H=\bC^N$ affine).\\

We next discuss the counterpart of the Artin vanishing Theorem \ref{Artinvan} in the complex analytic context for {\it q-complete} varieties in the following sense.

\bd\label{q-complete} \index{q-complete variety}
A complex variety $X$ is called {\it q-complete} ($q\in \bN_{0}$), 
if there exists a
{\it proper q-convex C$^{\infty}$-function} $f\colon  X \to \bR^{\geq 0}$ (with q-convexity in the sense of Definition \ref{q-convex function}).
In particular, the {\it $0$-complete} complex varieties are just the complex analytic \index{Stein space} {\it Stein spaces}
(see, e.g., \cite{Na}), and a closed complex subvariety of a  q-complete variety is again
 q-complete.
\ed

Then we need to consider a  Morse approximation of $f$ (as in \cite{Ben}) with possibly infinitely many critical points, 
and  get by Theorem~\ref{thm:plurisub} the following result (see, e.g., \cite[Corollary 6.1.2]{Sc} and compare also with \cite[Theorem 10.3.8]{KS}).

\bt[Vanishing theorem for q-complete varieties] \label{van-qc} \
Let $X$ be a complex analytic variety with a Whitney stratification $\cS$,
$f\colon  X\to \bR^{\geq 0}$ a {\it proper q-convex} function, and
$$K_{n}:=\{f\leq r_n\},\: U_{n}:=\{f<r_n\}, \quad \text{  for \ $r_{n}\nearrow \infty$}$$
 a sequence  of regular values of $f$ with respect to $\cS$ ($n\in \bN$).
  Then one gets for  $\cFb\in D^b_{\cS-wc}(X;R)$  a $\cS$-weakly constructible complex and
$q':=\min\,\{q, \dim X\}$:
\begin{enumerate}
\item If $\cFb\in \,^{p}D^{\leq m}_{\cS}(X;R)$, then
$$H^{k}(K_{n+1},K_{n};\cFb) \simeq  H^{k}(\{r_{n} \leq f \leq r_{n+1}\},\{f=r_{n}\};\cFb)=0
\quad \text{for $k>m+q'$},$$
so that $H^{k}(K_{n+1};\cFb)\to H^{k}(K_{n};\cFb)$ is surjective
for $k\geq m+q'$. In particular, the projective system 
$H^{k}(K_{n};\cFb)$, $n\in \bN$, 
satisfies the {\it Mittag-Leffler condition} for $k\geq m+q'$ 
so that 
\be
H^{k}(X;\cFb) \simeq \lim_{\leftarrow} \: H^{k}(K_{n};\cFb) = 0 \quad \text{for $k>m+q'$.}
\ee
\item If $\cFb\in \,^{p}D^{\geq m}_{\cS}(X;R)$, then
\be
H^{k}_{c}(X;\cFb) \simeq \lim_{\rightarrow} \: H^{k}_{c}(U_{n};\cFb) = 0 \quad  \text{for $k<m-q'$.}
\ee
\item If $\cFb\in \,^{p^{+}}D^{\geq m}_{\cS}(X;R)$ is constructible and $R$ is a Dedekind domain, then
\be
H^{m-q'}_{c}(X;\cFb) \simeq \lim_{\rightarrow} \: H^{m-q'}_{c}(U_{n};\cFb)\quad  \text{is torsion-free.}
\ee
\end{enumerate}
\et

Note that in the general complex analytic context one cannot expect any {\it finiteness or torsion-properties} for the corresponding cohomology groups.
If $X$ is the Stein space given by an infinite discrete set $X=\bN\subset \bC$ and $\cF$ is a sheaf on $X$, then
$$H^0(X;\cF)=\prod_{x\in X} \cF_x \quad \text{and} \quad H^0_c(X;\cF)=\bigoplus_{x\in X} \cF_x\:.$$
So even if all stalks $\cF_x$ ($x\in X$) are finitely generated (torsion) $R$-modules, this does not need to be the case for
$H^k(X;\cF)$ or $H^k_c(X;\cF)$. See \cite{Sce} for similar  finiteness results on {\it q-convex and q-concave complex varieties}
(using the same method of proof as in Theorem \ref{van-qc}).

\bex
Let $X$ be a q-complete variety of dimension $n=\dim X$ (e.g., a Stein space for $q=0$), with $q':=\min\,\{q, n\}$.  Then
\begin{enumerate}
\item The {\it homology} $H_{k}(X;R)=H^{-k}_c(X;\bD_X^{\bullet})=0 $ for $k> n+q'$.
 Moreover, $H_{n+q'}(X;R)$ is torsion-free in case $R$ is a Dedekind domain.
\item The {\it cohomology}
$H^{k}(X;R) = 0$ for $k> n+q'$.
\item Suppose $X$ is purely n-dimensional.
Then the {\it  intersection cohomology} \index{intersection cohomology}
$$IH^k(X;R) := H^{k-n}(X;IC_X)=0 \quad \text{ for $k>n+q'$, }$$
and \index{intersection homology} the {\it intersection homology} $IH_k(X;R) := H^{n-k}_c(X; IC_X)=0 $ for $k> n+q'$.
 Moreover, $IH_{n+q'}(X;R)$ is torsion-free in case $R$ is a Dedekind domain.
\end{enumerate}
\eex

The results in 1. and 2. are due to Hamm \cite{Ha1, Ha2}, whereas 3. is due to Goresky-MacPherson \cite[Section 6.9]{GMs}
(at least for Stein spaces). See \cite[Chapter 6]{Sc} for more general results and further references.\\

Before stating the relative versions of these vanishing results, let us recall the following.

\bd A morphism $f\colon  X\to Y$ of complex algebraic varieties is called \index{affine morphism} {\it affine}, if any point $y\in Y$ has an open affine neighborhood $U\subset Y$ such that
$f^{-1}(U)$ is affine. Similarly, a morphism $f\colon  X\to Y$ of complex analytic  varieties is called \index{q-complete morphism} \index{Stein morphism} {\it q-complete}, resp., {\it Stein}, if any point $y\in Y$ has an open neighborhood $U\subset Y$ such that
$f^{-1}(U)$ is q-complete, resp., Stein ($q\in \bN_0$). So, by definition, a {\it Stein morphism} corresponds to a {\it 0-complete morphism}.
\ed 

\bex  A finite morphism (e.g., a closed embedding) is an affine, resp., Stein morphism.
Similarly, if the closed embedding $i\colon  X\hookrightarrow Y$ is locally given by one algebraic, resp., analytic  equation $X=\{f=0\}$,
then the inclusion of the open complement $j\colon  U:=Y\backslash X\hookrightarrow Y$ is an affine, resp., Stein morphism.
Also, an {\it affine algebraic} morphism is {\it Stein} when viewed as a morphism of the underlying analytic varieties.
\eex

\br Let us explain the crucial property of a q-complete (resp., affine) morphism $f\colon  X\to Y$ used for (co)stalk calculations in the following Theorems.
In the analytic context consider a local embedding $(Y,y)\hookrightarrow (\bC^n,0)$ and  an open neighborhood $U\subset Y$ of $y$ such that
$f^{-1}(U)$ is q-complete, with $g\colon  f^{-1}(U)\to \bR^{\geq 0}$ a proper q-convex function. Consider the strongly plurisubharmonic distance function
 $r\colon  \bC^n\to \bR^{\geq 0}$ with $r(z):=\sum_{i=1}^n z_i\cdot \bar{z}_i$. Then for any small $r_0>0$ the function
\be
g':=g+\frac{1}{r-r_0}\circ f: \: f^{-1}(U\cap \{r<r_0\}) \to \bR^{\geq 0}\quad \text{is also proper and q-convex}
\ee
(see, e.g., \cite[p.429]{Sc} and \cite[Proposition 2.2, Proposition 2.4]{SV}). 
Therefore one can apply the 
vanishing Theorem \ref{van-qc}  to $f^{-1}(U\cap \{r<r_0\})$.

For an affine morphism $f\colon  X\to Y$, one takes a global affine embedding $(U,y) \hookrightarrow (\bC^n,0)$ of  an open affine neighborhood $U\subset Y$ of $y$ such that
$f^{-1}(U)$ is also affine with  $g\colon  f^{-1}(U)\to \bR^{\geq 0}$ a proper strongly plurisubharmonic and semialgebraic  function. 
Then $g'$ as before is a {\it proper strongly plurisubharmonic and semialgebraic  function} on the {\it open semi-algebraic} subset 
$f^{-1}(U\cap \{r<r_0\})\subset X$, so that one can use the Artin vanishing Theorem in the version of Remark \ref{opensemia}.
\er 

Let us now state the following {\it Artin-Gothendieck type} result in the complex algebraic context (see, e.g., \cite[Theorem 6.0.4]{Sc}).

\bt[Algebraic Artin-Grothendieck type Theorem] \label{algAG} \index{Artin-Gothendieck type theorem, algebraic}
Let $f\colon  X \to Y$ be an {\it affine} morphism of complex algebraic varieties,
with  $\cFb\in D^b_{(w)c}(X;R)$ a {\it complex algebraically} (weakly) constructible  sheaf complex
(respectively, complex algebraically constructible sheaf complex in case of a Dedekind domain $R$, if we want to use the dual perverse t-structure).
Then:
\begin{enumerate}
\item $\cFb\in \,^{p^{(+)}}D^{\leq n}(X;R)  \Rightarrow 
Rf_*\cFb\in  \,^{p^{(+)}}D^{\leq n}(Y;R)  $,
\item $\cFb\in \,^{p^{(+)}}D^{\geq n}(X;R)  \Rightarrow 
Rf_!\cFb \in  \,^{p^{(+)}}D^{\geq n}(Y;R) $.
\end{enumerate}
\et

Compare  with \cite[Section 4]{BBD}, especially 
\cite[Theorem 4.1.1]{BBD},
for the corresponding relative counterpart for an affine morphism 
in the context
of the {\it perverse t-structure}  in $l$-adic
cohomology. In the case $Y=\{pt\}$ a point space, Theorem \ref{algAG} just reduces to Theorem \ref{Artinvan}.\\

In the complex analytic context one has in addition to assume that the corresponding direct image complexes are again (weakly) constructible
(see, e.g., \cite[Corollary 6.0.8]{Sc} and \cite[Proposition 10.3.17]{KS} for the case of a Stein map).

\bt[Analytic Artin-Grothendieck type Theorem] \label{anaAG} \index{Artin-Gothendieck type theorem, analytic}
Let $f\colon  X \to Y$ be a {\it q-complete } morphism of complex analytic varieties (e.g., a Stein map for $q=0$),
with  $\cFb\in D^b_{(w)c}(X;R)$ a (weakly) constructible  sheaf complex.
Then:
\begin{enumerate}
\item $\cFb\in \,^{p}D^{\leq n}(X;R)$ and $Rf_*\cFb\in  D^b_{(w)c}(Y;R) \Rightarrow 
Rf_*\cFb\in  \,^{p}D^{\leq n+q}(Y;R)  $,
\item $\cFb\in \,^{p}D^{\geq n}(X;R)$ and $Rf_!\cFb\in  D^b_{(w)c}(Y;R)  \Rightarrow 
Rf_!\cFb \in  \,^{p}D^{\geq n-q}(Y;R) $.
\item Assume $\cFb\in D^b_{c}(X;R)$ is constructible with $R$ a Dedekind domain. then\\
$\cFb\in \,^{p^{+}}D^{\geq n}(X;R)$ and $Rf_!\cFb\in  D^b_{c}(Y;R)  \Rightarrow 
Rf_!\cFb \in  \,^{p^{+}}D^{\geq n-q}(Y;R) $.
\end{enumerate}
\et

In the case $Y=\{pt\}$ a point space, Theorem \ref{anaAG} corresponds  to Theorem \ref{van-qc}.

Let us finish this section with the following:

\bex[Relative weak Lefschetz theorem for singular spaces] \label{relweaklef} \index{Relative weak Lefschetz theorem}
Let $V\to Y$ be a complex algebraic (or analytic) vector bundle, and $W\hookrightarrow V$ be a subvector bundle with $rank\; V=rank\; W +1$.
Let $i\colon  \bP(W)\hookrightarrow \bP(V)$ the  closed inclusion of the associated projective bundles, with open complement 
$j\colon  U:=\bP(V)\backslash \bP(W) \hookrightarrow \bP(V)$. Then the projection $\pi\colon  U\to Y$ is an affine (resp., Stein) morphism.

Assume $\cFb\in D^b_{(w)c}(\bP(V);R)$ is (weakly) constructible
(respectively,  constructible in case of a Dedekind domain $R$, if we want to use the dual perverse t-structure), so that
$$R\pi_!j^*\cFb,\: R\pi_*j^*\cFb\in D^b_{(w)c}(Y;R)\:.$$
Then
$$j^{*}\cFb\in \,^{p}D^{\leq n}(U;R)  \Rightarrow
R\pi_*j^*\cFb\in  \,^{p}D^{\leq n}(Y;R)  \:,$$
and
$$j^{*}\cFb\in \,^{p^{(+)}}D^{\geq n}(U;R)  \Rightarrow
R\pi_!j^*\cFb\in  \,^{p^{(+)}}D^{\geq n}(Y;R)  \:.$$
\eex

\section{Nearby and vanishing cycles, applications}\label{sec:nvc}

In this section we recall the construction of Deligne's nearby and vanishing cycle functors (\cite{Gr2}), and indicate their relation with perverse sheaves. We continue to assume that the base ring $R$ is commutative and noetherian, of finite global dimension, and we work with (weakly) constructible complexes of sheaves of $R$-modules in the complex algebraic (or analytic) context.

\subsection{Construction}
Let $f\colon X \to \bC$ be a morphism from a complex algebraic (or analytic)  variety $X$ to $\bC$. Let $X_0=f^{-1}(0)$ be the central fiber, with inclusion map $i\colon X_0 \hookrightarrow X$. Let $X^*:=X \setminus X_0$ and let $f^*\colon X^* \to \bC^*$ be the induced morphism  to the punctured affine line. Consider the following cartesian diagram:
\be\label{diagr-psi}
\xymatrix{
X_0 \ar[d] \ar@{^(->}[r]^i & X \ar[d]_f & \ar@{^(->}[l]_j \ar[d]_{f^*} X^* & \ar[l]_{\widehat{\pi}} \wti{X^*} \ar[d] \\
\{ 0 \} \ar@{^(->}[r] & \bC & \ar@{^(->}[l] \bC^* & \ar[l]^\pi \widetilde{\bC^*}\:,
}
\ee
where $\pi\colon \widetilde{\bC^*} \to \bC^*$ is the infinite cyclic (universal) cover of $\bC^*$ given by $z \mapsto \exp(2\pi i z)$. Then $\widehat{\pi}\colon \wti{X^*} \to X^*$ is an infinite cyclic cover with deck group $\pi_1(\bC^*)\simeq \bZ$. Note that $\pi$ is only a holomorphic but not an algebraic map, so that this fiber square only exists in the complex analytic category, even if we start with a morphism $f\colon  X\to \bC$ in the complex algebraic context.

\bd\label{ncd} 
The \index{nearby cycle functor} {\it nearby cycle functor of $f$} assigns to a bounded complex $\cF^\bullet \in D^b(X;R)$ the complex on $X_0$ defined by
\begin{equation} \psi_f\cFb:=i^*R(j \circ \widehat{\pi})_*(j \circ \widehat{\pi})^* \cFb \simeq i^*Rj_*R\widehat{\pi}_*\widehat{\pi}^*j^*\cFb
\in D^b(X_0;R)\:.\end{equation}  
\ed

\br
By definition,  $\psi_f\cFb$ depends only on the restriction $j^*\cFb$ of $\cFb$ to $X^*$. In the complex analytic context, it would have been enough for the  definition of the nearby cycles  to start with a holomorphic map $f\colon  X\to D\subset \bC$ to a small open disc $D$ around zero in the complex plane, with $\pi\colon \widetilde{D^*} \to D^*$ the induced infinite cyclic (universal) cover of $D^*$.
\er
Note that $\pi_1(\bC^*)\simeq \bZ$ acts naturally on $R\widehat{\pi}_*\widehat{\pi}^*$, with an induced action on the nearby cycles $\psi_f\cFb$.
In our complex context we choose the generator of $\pi_1(\bC^*)\simeq \bZ$ fitting with the complex orientation of $\bC^*$ and call the induced action
$$h=h_f\colon  \psi_f\cFb\to \psi_f\cFb$$ the \index{monodromy} {\it monodromy automorphism} of the nearby cycles. The adjunction morphism 
$$id\to R(j \circ \widehat{\pi})_*(j \circ \widehat{\pi})^*$$ induces the {\it specialization  map} $sp\colon  i^*\to \psi_f$ commuting with the monodromy automorphism $h$ acting trivially on $i^*$ (i.e., acting as the identity on $i^*$). Let us now explain the use of the base change induced by the universal cover $\pi\colon  \widetilde{\bC^*}\to \bC^*$
in the definitions above in a simple but important example.

\bex[Nearby cycle functor for the identity map] \label{psi-id} Let $f=id$ be the identity map of $\bC$, with $i\colon  \{0\}\hookrightarrow \bC$ the inclusion of the point zero. Similarly, let $i_x\colon  \{x\}\hookrightarrow 
D_r:=\{|z|<r\}\subset \bC$, for $0<|x|<r\leq \infty$, be the inclusion of a (nearby) point $x\neq 0$ in a (small) open disc of radius $r$ around zero.
Finally consider a complex $\cFb\in D^b_{\cS-wc}(\bC;R)$ which is weakly constructible with respect  to the Whitney stratification $\cS$ of $\bC$ given by the two strata $S=\{0\}$ and $S'=\bC^*$. 

Then the cohomology sheaves $\cH^k(j^*\cFb)$ are only {\it locally constant} on $\bC^*$, but their pullbacks
$\cH^k(\widehat{\pi}^*j^*\cFb)$ to $\bC\simeq  \widetilde{\bC^*}$ (or their restrictions to $ \widetilde{D_r^*}$) are locally constant and therefore 
{\it constant} (for all $k\in \bZ$)
since $\bC\simeq  \widetilde{\bC^*}$ (or $ \widetilde{D_r^*}$) is convex and {\it contractible}. Then
$$\begin{CD}
\psi_{id}\cFb @< \sim << R\Gamma(D^*_r; R\widehat{\pi}_*\widehat{\pi}^*j^*\cFb)\simeq 
R\Gamma( \widetilde{D_r^*};\widehat{\pi}^*j^*\cFb) \\
@.  @VV \wr V \\
@.  i^*_{\tilde{x}}\widehat{\pi}^*j^*\cFb\simeq i^*_x\cFb
\end{CD}$$
for $i_{\tilde{x}}\colon  \{\tilde{x}\}\hookrightarrow \widetilde{D_r^*}$ the inclusion of a point $\tilde{x}$ with $\pi(\tilde{x})=x$. And the {\it monodromy} $h$ on
$H^k(\psi_{id}(\cFb))$ gets identified with the {\it monodromy of the local system} $\cH^k(j^*\cFb)$  acting on  $i^*_x\cH^k(j^*\cFb)$ (for all $k\in \bZ$). Finally the {\it specialization} map is given by
$$\begin{CD}
sp\colon  \;i^*\cFb @< \sim << R\Gamma(D_r;\cFb) @>>> i^*_x\cFb \simeq \psi_{id}\cFb\:.
\end{CD}$$

\eex

\bd In the context of the diagram (\ref{diagr-psi}), \index{vanishing cycle functor}
the  {\it vanishing cycles} of $\cF^\bullet\in D^b(X;R)$, denoted  $\varphi_f \cF^\bullet$, is the bounded complex on $X_0$ defined by taking ``the'' cone of 
 the \index{specialization  map} {\it specialization  map} $sp\colon  i^*\cFb\to \psi_f\cFb$.
 In particular, one gets a unique distinguished triangle
\be\label{canmo}
i^*\cF^\bullet \overset{{sp}}{\lra} \psi_f \cF^\bullet \overset{can}{\lra} \varphi_f \cF^\bullet \overset{[1]}\lra
\ee
in $D^b(X_0;R)$.
\ed

\br
Note that cones are not functorial, but in the above construction one can for example work with (a suitable truncation of) the canonical flabby resolution  to get $\varphi_f$ as a functor (see \cite[Chapter 8]{KS} or \cite[pp. 25-26]{Sc} for more details). The vanishing cycle functor also comes equipped with a
\index{monodromy} monodromy automorphism, denoted also by $h$, so that $h$ induces an automorphism of the  triangle (\ref{canmo}).
\er

\bex Assume $\cFb\in D^b(X;R)$ is supported on $X_0$, i.e., $j^*\cFb\simeq 0$. Then $\psi_f\cFb\simeq 0$ and
$(\varphi_f\cFb)[-1] \simeq i^*\cFb$, with a trivial monodromy action $h$.
\eex

\bex[Vanishing  cycle functor for the identity map] \label{phi-id}
Consider the context of the Example \ref{psi-id}. Then for any point $0\neq x\in D_r$ in a small open disc $D_r$ around zero:
$$(\varphi_{id}\cFb)[-1] \simeq R\Gamma(D_r,\{x\};\cFb) \simeq 
\bigl(R\Gamma_{\{l\geq 0\}}(\cFb)\bigr)_{0} =: LMD(\cFb,l,0)$$
for any $\bR$-linear map $l\colon  \bC\to\bR$ with $l(x)<0$. In particular the following properties are equivalent:
\begin{enumerate}
\item $sp\colon  i^*\cFb\to \psi_{id}\cFb$ is an isomorphism,
\item  $\varphi_{id}\cFb\simeq 0$,
\item all cohomology sheaves $\cH^k(\cFb)$ are locally constant on $\bC$ or $D_r$ ($k\in \bZ$),
\item the normal  Morse datum  $NMD(\cFb,\{0\})\simeq LMD(\cFb,l,0)\simeq 0$.
\end{enumerate}
\eex

Next we explain why the nearby and vanishing cycle functors preserve (weak) constructibility.
Consider a complex algebraic (or analytic) morphism $f\colon  X\to \bC$ as in diagram (\ref{diagr-psi}).
Assume $X$ is endowed with a Whitney stratification $\cS$ such that  $X_0=\{f=0\}$ is a union of strata, with induced stratification $\cS|_{X_0}$.
Similarly for the induced stratification $\cS|_{X^*}$ of the open complement  $X^*=\{f\neq 0\}$.
But  $\widehat{\pi}\colon \wti{X^*} \to X^*$ is an infinite cyclic covering  with deck group $\pi_1(\bC^*)\simeq \bZ$, so that
$\wti{X^*}$ gets an induced complex analytic Whitney stratification $\widetilde{\cS}$, making $\widehat{\pi}\colon \wti{X^*} \to X^*$ a {\it stratified map}
such that, for any stratum $S\in \cS$, $\widehat{\pi}\colon \widehat{\pi}^{-1}(S)\to S$ is also an infinite cyclic covering  with deck group 
$\pi_1(\bC^*)\simeq \bZ$, i.e., a locally trivial fibration with fiber $\bZ$. But this implies by induction on $\dim X^*$ the following (see, e.g. \cite[Corollary 4.2.1(4)]{Sc}).

\bl $R\widehat{\pi}_*\widehat{\pi}^*$ maps $D^b_{\cS|_{X^*}-wc}(X^*;R)$ to itself.
\el

Together with Example \ref{ex410} this implies (also in the complex algebraic context)  the following important fact.

\bc\label{nv-weak} The nearby and vanishing cycle functors $\psi_f, \varphi_f$ induce 
$$ \psi_f, \; \varphi_f\colon  \:D^b_{(\cS-)wc}(X;R) \to D^b_{(\cS|_{X_0}-)wc}(X_0;R)\:,$$
i.e., they preserve weak constructibility (with respect to $\cS$ and $\cS|_{X_0}$).
\ec

But note that $R\widehat{\pi}_*\widehat{\pi}^*$ {\it does not preserve} $D^b_{(\cS|_{X^*}-)c}(X^*;R)$, i.e., constructibility.
Let $g$ be a generator of $\pi_{1}(\bC^{*}) \simeq \bZ$, given by the complex 
orientation of $\bC^{*}$, which we interpret now as an {\it automorphism} of 
$\widetilde{X^{*}}$.
Let $\cG$ be a sheaf on $X^*$. We claim that the following sequence of sheaves
on $X^*$ is {\it exact}:
$$ \begin{CD}
0 @>>>\cG @>ad_{\widehat{\pi}}>>\widehat{\pi}_*\widehat{\pi}^*\cG @>g^*- id >> 
\widehat{\pi}_*\widehat{\pi}^*\cG @>>> 0 \:.
\end{CD} $$
Take a ball $B$ in $\bC^{*}$ so that the restriction of the
universal covering map of $\bC^{*}$ to $B$ is
isomorphic to the projection $B\times \bZ\to B$, with $g$ corresponding
to the translation $\bZ\to \bZ, i\mapsto i+1$.
If we take an open subset $V\subset f^{-1}(B)$, then $\widehat{\pi}\colon  \widehat{\pi}^{-1}(V)\to V$
is isomorphic to the projection $V\times \bZ\to V$, with $g$ acting as 
before on the second factor. Then 
\[\Gamma(V,\widehat{\pi}_*\widehat{\pi}^*\cG) \simeq \prod \nolimits_{i\in \bZ} \;
\Gamma(V,\cG) \]
such that $g^*$ acts by the {\it permutation} $i\mapsto i+1$.
Moreover the adjunction map corresponds to the diagonal embedding
$\Gamma(V, \cG) \to \prod_{i\in \bZ}\; \Gamma(V,\cG)$.  
Then the sequence
$$\begin{CD}
0 @>>> \Gamma(V,\cG) @>diag>>\prod_{i\in \bZ}\; \Gamma(V,\cG)
 @>g^*-id >> \prod_{i\in \bZ}\; \Gamma(V,\cG)  @>>> 0
\end{CD}$$
is exact. Since each point of $X^*=\{f\neq 0\}$ has a  fundamental system
of open neighborhoods $V$ as before, this implies our claim.\\

By using a flabby resolution we therefore get the distinguished 
triangles (see, e.g., \cite[(5.88) on p.369]{Sc})
$$\begin{CD}
\cG @>ad_{\widehat{\pi}}>> R\widehat{\pi}_*\widehat{\pi}^*\cG @>g^*-id >> 
R\widehat{\pi}_*\widehat{\pi}^*\cG @>[1]>> 
\end{CD} $$
for any $\cGb\in D^b(X^*;R)$, and 
\begin{equation} \label{eq:dtrwang} \begin{CD}
i^{*}Rj_{*}j^{*}\cFb @>ad_{\widehat{\pi}} >> \psi_{f}\cFb @> h_{f} -id  >> \psi_{f}\cFb @>[1]>> 
\end{CD} \end{equation}
for any $\cFb\in D^b(X;R)$.\\

Consider also the distinguished triangle
\be \label{eq:troct} \begin{CD}
i^{*}\cFb @>ad_{j}>> i^{*}Rj_{*}j^{*}\cFb @>>> i^{!}\cFb [1] @>[1]>> \: .
\end{CD} \ee

Since the map $can\colon  i^*\cFb\to \psi_f\cFb$ factorizes as 
$$ \begin{CD}
can\colon  i^{*}\cFb @>ad_{j}>> i^{*}Rj_{*}j^{*}\cFb    @> ad_{\widehat{\pi}} >>\psi_{f}\cFb\:,
\end{CD} $$
we get from the distinguished triangles~(\ref{eq:troct}),
~(\ref{eq:dtrwang}),~(\ref{canmo}) and the 
{\it octahedral axiom} a distinguished triangle

\be \label{eq:dtrvar} \begin{CD}
i^{!}\cFb[1] @>>> \varphi_{f}\cFb @> var >> \psi_{f}\cFb @>[1]>> 
\end{CD} \ee
for $\cFb\in D^b(X;R)$.
Here the \index{variation morphism} variation morphism $$var\colon \varphi_f \cF^\bullet  \to\psi_f \cF^\bullet$$  can be defined by the cone of the pair of morphisms
(applied to a flabby resolution):
$$(0, h - id)\colon [i^*\cF^\bullet\to \psi_f \cF^\bullet] \lra [0 \to \psi_f \cF^\bullet]\:,$$
with $h=h_f$ the monodromy automorphism, so that 
 \be can \circ var=h - id \quad \text{and}  \quad var \circ can=h - id\:.\ee

\br Note that the {\it variation morphism} $var$ depends of the choice of a generator $g\in \pi_1(\bC^*)\simeq \bZ$, i.e., on the choice of an orientation for 
$\bC^*$. Moreover there are then two choices for $var$ in the literature, using a different sign convention so that
$$ can \circ var=\pm(h - id)=\mp(id-h) \quad \text{and}  \quad var \circ can=\pm(h - id)=\mp(id-h)\:.$$
For example, our choice here only fits with $-var$ as used in \cite[Equation (8.6.8)]{KS} and \cite[Equation (5.90)]{Sc}.
\er

Let us now explain an important description of the {\it (co)stalks of the nearby cycles} in terms of {\it local Milnor fibers}.
Consider a complex algebraic (or analytic) morphism $f\colon  X\to \bC$ as in diagram (\ref{diagr-psi}).
Assume $X$ is endowed with a Whitney stratification $\cS$ such that  $X_0=\{f=0\}$ is a union of strata, with induced stratification $\cS|_{X_0}$.
By a classical  theorem of Hironaka \cite[Corollary 1 to Theorem 2, p.248] {Hir} (and see also \cite[Corollary 1.3.5.1]{LeTe}),
one can always refine a given stratification so that it satisfies the \index{$a_f$-condition}
$a_f$-condition (\ref{af}). In fact, by a more recent result \cite[Theorem  4.2.1]{BMM}, this $a_f$-condition is true for the given  Whitney stratification $\cS$, with $X_0$ a union of strata. \\

Note that the $a_f$-condition is needed to have  a {\it local Milnor fibration} at a given point $x\in X_0$, with {\it  Milnor fiber} \index{Milnor fiber}
$$(M_{f,x},\partial M_{f,x}) :=(X\cap \{f=w\}\cap \{r\leq \epsilon\}, X\cap \{f=w\}\cap \{r= \epsilon\}) $$
$$ \text{and} \quad \mathring{M}_{f,x}:=M_{f,x}\backslash \partial M_{f,x}$$
for $0<|w|\ll \epsilon$
and a  general holomorphic function germ
$f\colon  (X,x)\to (\bC,0)$ on a singular complex analytic variety $X$  (see, e.g., \cite[Example 1.3.3]{Sc} and \cite{Le2}).
Here we consider a local embedding $(X,x)\hookrightarrow (\bC^n,0)$ with $r(z):=\sum_{i=1}^n z_i\cdot\bar{z}_i$ the distance to $x=0$.
By the {\it curve selection lemma} one can assume that (locally near $x\in X_0$) $w$ is a stratified regular value, i.e., $\{f=w\}$ is transversal to $\cS$ for
$0<|w|$ small enough. Then by the $a_f$-condition, and for $0<|w| \ll \epsilon$ small enough, $r=\epsilon$ is a stratified regular value of $r$ with respect to 
the induced Whitney stratification of $X\cap \{f=w\}$ (see, e.g., \cite[Example 1.1.3]{Sc}).

Then the  following local calculation is a direct consequence of the definition of the {\it nearby cycle functor} and the existence of such a {\it local Milnor fibration}
(see, e.g., \cite[Example 5.4.2]{Sc}).
\bp\label{Dec} For every $x \in X_0$, with $i_x\colon  \{x\}\hookrightarrow X_0=\{f=0\}$ the inclusion, there are isomorphisms:
\be\label{de3} i_x^*(\psi_f\cF^\bullet) \simeq  R\Gamma(M_{f,x},\cFb) \simeq  R\Gamma(\mathring{M}_{f,x},\cFb) 
\ee 
compatible with the corresponding {\it monodromy} actions,
and
\be \label{de4}  i_x^!(\psi_f\cF^\bullet) \simeq  R\Gamma(M_{f,x},\partial M_{f,x},\cFb) \simeq  R\Gamma_c(\mathring{M}_{f,x},\cFb) 
\ee
for any $\cFb\in D^b_{wc}(X;R)$.
\ep

By Proposition \ref{prop:T-stalk} and the triangle (\ref{canmo}) we get the following.
\bc\label{nvc}
Let $T\subset D^b(\{pt\};R)$ be a fixed \index{null system} ``null system'', i.e., a full triangulated subcategory stable by isomorphisms.
If $\cFb\in D^b_{(\cS-)T-stalk}(X;R)$, then also $$\psi_f\cFb, \;\varphi_f\cFb \in D^b_{(\cS|_{X_0}-)T-stalk}(X;R)\:.$$
In particular, $\psi_f\cFb$ and $\varphi_f\cFb$ are constructible for $\cFb$ constructible.
\ec

\br By taking $T= D^b_{c,\chi=0}(\{pt\};R)$, we get the nearby and vanishing cycles for constructible functions:
$$\psi_f,\; \varphi_f\colon  \: CF_{(\cS)}(X) \to CF_{(\cS|_{X_0})}(X_0)\:,$$
with 
$$\psi_f(\alpha)(x):=\int_{M_{f,x}} \alpha d \chi \quad \text{and} \quad
\varphi_f(\alpha)(x):=\int_{M_{f,x}} \alpha d \chi -\alpha(x)$$
for $\alpha\in CF_{(\cS)}(X)$ and $x\in X_0=\{f=0\}$.
\er

Using Proposition \ref{Dec} and the distinguished triangle (\ref{canmo}), one gets the following (see, e.g., \cite[Lemma 5.4.1, Example 5.4.1]{Sc}).
\bc\label{van-stalk}
For every $x \in X_0=\{f=0\}$, with $i_x\colon  \{x\}\hookrightarrow X_0=\{f=0\}$ the inclusion, there are isomorphisms:
\begin{eqnarray}\label{de5}
 i_x^*(\varphi_f \cF^\bullet)[-1]  \simeq R\Gamma(\mathring{B}_{\epsilon, x}, \mathring{B}_{\epsilon, x} \cap \{f=w\}, \cF^\bullet)  \nonumber \\
 \simeq 
\bigl(R\Gamma_{\{Re(f)\geq 0\}}(\cFb)\bigr)_{x} =: LMD(\cFb,Re(f),x) 
\end{eqnarray}
for $\cFb\in D^b_{wc}(X;R)$ and  $0<|w|\ll \epsilon$ small enough.
Here  $\mathring{B}_{\epsilon, x}=X\cap \{r<\epsilon\}$ is the intersection of $X$ with a small open $\epsilon$-ball
around $x\in X$.
\ec

\bex\label{stcv} 
As a special case of (\ref{de5}), let $\cF^\bullet = R_X$ be the constant sheaf on $X$. Since $\mathring{B}_{\epsilon, x}\cap X_0$ is contractible,  one  gets \index{reduced cohomology}
\begin{align*}
\cH^k(\varphi_f R_X)_x \simeq  \widetilde{H}^k(M_{f,x}; R)
\end{align*} 
is the {\it reduced cohomology} of the Milnor fiber  $M_{f,x}$ of $f$ at $x$. 
If, moreover, $X$ is smooth, then Milnor fibers at smooth points of $X_0$ are contractible, so the above calculation yields the inclusion: 
$$\supp(\varphi_f R_X):=\bigcup_{k}  \supp  \cH^k(\varphi_f \cFb)  \subseteq \Sing(X_0).$$ 
\eex

A more general estimation of the \index{support} support of vanishing cycles is provided by the following result (see, e.g., \cite{Ma} or \cite[Remark 4.2.4]{Sc}).
\bp\label{410} Let $X$ be a complex algebraic (or analytic) variety with a given Whitney stratification $\cS$, and let $f\colon X \to \bC$ be a morphism
with $X_0=\{f=0\}$. For every $\cS$-weakly constructible complex $\cFb$ on $X$ and every integer $k$, one has the inclusion
\be \supp  \cH^k(\varphi_f \cFb)\subseteq X_0 \cap \Sing_{\cS}(f),\ee
where  $$\Sing_{\cS}(f):=\bigcup_{S \in \cS} \Sing(f|_S)$$
is the stratified singular set of $f$ with respect to the stratification $\cS$.
\ep

\bex[Isolated stratified critical point] \label{isostr}\index{isolated stratified critical point}
In the context of Proposition \ref{410}, assume that $x\in X_0=\{f=0\}$ is an {\it isolated} stratified critical  point of $f$.
Then 
$$i_x^!(\varphi \cFb)[-1]\simeq  i_x^*(\varphi \cFb)[-1] \simeq \bigl(R\Gamma_{\{Re(f)\geq 0\}}(\cFb)\bigr)_{x} =: LMD(\cFb,Re(f),x) $$
for $i_x\colon  \{x\}\to X_0$ the point inclusion. Applying $i_x^*$ to the distinguished triangle (\ref{canmo}), one gets
$$\begin{CD}
i_x^*(\psi_f\cFb)[-1] @> can >> i_x^*(\varphi_f\cFb)[-1] @>>> i_x^*(i^{*}\cFb)  @>[1]>> \: ,
\end{CD}$$
with $i_x^*(\psi_f\cFb)\simeq   R\Gamma(M_{f,x},\cFb)$ for a \index{Milnor fiber} local Milnor fiber $M_{f,x}$ as in (\ref{de3}).\\
Applying $i_x^!$ to the distinguished triangle (\ref{eq:dtrvar}), one gets
$$\begin{CD}
i_x^!(i^{!}\cFb)  @>>> i_x^*(\varphi_f\cFb)[-1] @> var >>   i_x^!(\psi_f\cFb)[-1]      @>[1]>> \: ,
\end{CD}$$
with $ i_x^!(\psi_f\cF^\bullet) \simeq  R\Gamma(M_{f,x},\partial M_{f,x},\cFb) $ for a local Milnor fiber $(M_{f,x},\partial M_{f,x})$ as in (\ref{de4}).
\eex

\br\label{rNMD} Consider a local embedding $(X,x)\hookrightarrow (\bC^n,x)$, with $N$ a {\it normal slice} to $x\in S$ for a stratum $S\in \cS$.
Assume $g\colon  (\bC^N,x)\to (\bC,0)$ is a holomorphic function germ such that the covector $dg_x$ is {\it non-degenerate} with respect to $\cS$.
Then $x$ is an isolated stratified critical point of $g|_N$ with respect to the induced Whitney stratification $\cS|_N$ of $X\cap N$, so that 
Example \ref{isostr} gives for this case the distinguished triangles from Proposition \ref{prop:NMDtriangleC}.
\er

The nearby and vanishing cycle functors have the following {\it base change properties}.
Consider a cartesian diagram of morphism 
\be\begin{CD}
Y_0 @> k' >> X_0 @>>> \{0\}\\
@V i' VV @VV i V  @VVV\\
Y @> k >> X @> f >> \bC\:,
\end{CD}\ee
with $f':=f\circ k$. Then one has the following {\it base change isomorphisms} (see, e.g., \cite[Remark 4.3.7, Lemma 4.3.4]{Sc}).

\bp[Base change isomorphisms for nearby and vanishing cycles] \label{bc-nv} \index{base change isomorphism}
The following base change isomorphisms commute with the  maps $can$ and $var$:
\begin{enumerate}
\item Assume $k\colon  Y\to X$ is {\it proper}. Then
$$Rk'_*(\psi_{f'}\cFb )\simeq \psi_{f}(Rk_*\cFb) \quad \text{and} \quad Rk'_*(\varphi_{f'}\cFb) \simeq \varphi_{f}(Rk_*\cFb)$$
for all $\cFb\in D^b(Y;R)$.
\item Assume $k\colon  Y\to X$ is {\it smooth}. Then
$$k'^*(\psi_f\cFb )\simeq \psi_{f'}(k^*\cFb) \quad \text{and} \quad k'^*(\varphi_f\cFb) \simeq \varphi_{f'}(k^*\cFb)$$
for all $\cFb\in D^b(X;R)$.
\item Assume $X\hookrightarrow M$ is a closed subvariety of the complex algebraic (or analytic) manifold $M$, with $\cS$ a Whitney stratification of $X$.
Let $N\hookrightarrow M$ be a closed complex algebraic (or analytic) submanifold which is {\it transversal} to $\cS$ (i.e., transversal to all strata $S\in \cS$),
with $k\colon  Y:=X\cap N\hookrightarrow X$ the induced inclusion. Then
$$k'^*(\psi_f\cFb )\simeq \psi_{f'}(k^*\cFb) \quad \text{and} \quad k'^*(\varphi_f\cFb) \simeq \varphi_{f'}(k^*\cFb)$$
for all $\cFb\in D^b_{\cS-wc}(X;R)$.
\end{enumerate}
\ep

Then one gets by {\it proper base change} and the Examples
\ref{psi-id} and \ref{phi-id} the following (see, e.g., \cite[Example 1.1.1]{Sc}).

\bex\label{De} Let $f\colon  X\to \bC$ be a {\it proper} morphism, with $X_0=\{f=0\}$.  Then 
\be\label{de8} R\Gamma(X_0, \psi_f \cFb) \simeq \psi_{id}(Rf_* \cFb) \simeq i_x^*(Rf_*\cFb) \simeq R\Gamma(\{f=x\}, \cFb)
\ee
and
\begin{eqnarray} R\Gamma(X_0, \varphi_f \cFb)[-1] \simeq \varphi_{id}(Rf_* \cFb)[-1] \quad \quad  \nonumber \\
 \simeq R\Gamma(D_r,\{x\}, Rf_*\cFb) \simeq R\Gamma(\{|f|<r,\{f=x\}, \cFb)
\end{eqnarray}
for $\cFb\in D^b_{wc}(X;R)$ and
a point $0\neq x\in D_r$ in a small open disc $D_r$ around zero. In particular, it follows from \eqref{canmo} and \eqref{de8} that for any $0\neq x\in D_r$ as above, one has for $f$ proper the following \index{specialization sequence} {\it specialization sequence}: 
\be
\cdots \to H^k(X_0;\cFb) \to H^k(\{f=x\};\cFb) \to H^k(X_0;\varphi_f\cFb) \to H^{k+1}(X_0;\cFb) \cdots
\ee
In the case $\cFb=R_X$ and by analogy with the local situation, the groups $H^*(X_0;\varphi_fR_X)$ are usually referred to as the \index{vanishing cohomology} the {\it vanishing cohomology of $f$} (see, e.g., \cite{MP}).
\eex


\subsection{Relation with perverse sheaves and duality}
Let $f\colon X \to \bC$ be a morphism of complex algebraic (or analytic) varieties.
The behavior of the nearby and vanishing cycle functors with regard to Verdier duality is described by the following result (for instance, see \cite[Theorem 3.1, Corollary 3.2]{Ma2}).
\bt\label{nvpp}
The shifted functors  $\psi_f[-1]$ and $\varphi_f[-1]$ commute with the \index{Verdier dual}
Verdier duality functor $\cD$ up to natural isomorphisms. 
\et

Note that this duality result also fits with the following behavior of the nearby and vanishing cycle functors with regard to the (dual) perverse t-structure
(but this is not used in its proof, see, e.g., \cite[Theorem 6.0.2]{Sc}).

\bt\label{thm:nearcycle}
Let $f\colon X \to \bC$ be a morphism of complex algebraic (or analytic) varieties, with $X^*:=\{f\neq 0\}$.
Assume $\cFb\in D^b(X;R))$ is weakly constructible (resp., $\cFb$ is constructible with $R$ a Dedekind domain, in case we want to use the dual perverse t-structure).  Then we have: \index{perverse t-structure} \index{dual perverse t-structure}
\begin{enumerate}
\item $j^{*}\cFb\in\: ^{p^{(+)}}D^{\leq n}(X^*;R) \: \Rightarrow \:
(\psi_{f}\cFb)[-1] \in\:  ^{p^{(+)}}D^{\leq n}(X_{0};R)$.
\item $j^{*}\cFb\in\: ^{p^{(+)}}D^{\geq n}(X^*;R) \: \Rightarrow \:
(\psi_{f}\cFb)[-1] \in\:  ^{p^{(+)}}D^{\geq n}(X_{0};R)$.
\item $\cFb\in\: ^{p^{(+)}}D^{\leq n}(X;R) \: \Rightarrow \:
(\varphi_{f}\cFb)[-1] \in\:  ^{p^{(+)}}D^{\leq n}(X_{0};R)$.
\item $\cFb\in\: ^{p^{(+)}}D^{\geq n}(X;R) \: \Rightarrow \:
(\varphi_{f}\cFb)[-1] \in\:  ^{p^{(+)}}D^{\geq n}(X_{0};R)$.
\end{enumerate} \et

\begin{proof}
Note that the result (3.) resp.,  (4.) for the {\it vanishing cycle functor} follows directly form
the corresponding result (1.) resp., (2.) for the nearby cycle functor, if one uses the distinguished
triangle~(\ref{canmo}) or~(\ref{eq:dtrvar}).

The argument for the nearby cycles is similar
to the proof of Theorem~\ref{thm:plurisub}. But this time we use the
corresponding description of Proposition~\ref{Dec}  for the {\it (co)stalk
of the nearby cycle functor}. 
Choose a complex algebraic  (or analytic) Whitney stratification $\cS$
of $X$ with $X_0=\{f=0\}$ and $X^*=\{f\neq 0\}$
 a union of strata, which
satisfies the \index{$a_f$-condition}  {\it $a_{f}$-condition} of Thom (\ref{af}).
Then we already know, by Corollary \ref{nv-weak}, that $\psi_{f}\cFb$ is (weakly) constructible with respect to the induced Whitney stratification
$\cS|_{X_0}$ of $X_0$.

Consider a point $x\in S$ for a stratum $S\subset \{f=0\}$ of dimension $s$. 
First we assume that $S=\{x\}$ is a point stratum. By Proposition~\ref{Dec}  we get

$$
i_{x}^{*}(\psi_{f}\cFb)\simeq
R\Gamma(X\cap\{r\leq \delta, f=w\},\cFb) \:,
$$
and
$$
i_{x}^{!}(\psi_{f}\cFb)\simeq
R\Gamma(X\cap\{r\leq \delta, f=w\}, X\cap\{r=\delta, f=w\},\cFb) 
$$
for $0<|w|\ll \delta\ll 1$, with $i_{x}\colon \{x\}\to \{f=0\}$ the inclusion
and 
$$r(z):=\sum \nolimits_{i=1}^{n}\: z_{i}\bar{z}_{i}  \quad
\text{in a local embedding} \quad
(X,x)\hookrightarrow  (\bC^{n},0) \:.$$
Then $L:=\{f=w\}$ is {\it transversal} to $\cS$ near $x$
for $0<|w|\ll 1$ (by the curve selection lemma). So  we get by Proposition \ref{prop-tr}:
$$\cFb|_L[-1] \in \;^{p^{(+)}}D^{\leq n}(L;R)  \quad \text{or} \quad
\cFb|_L[-1] \in \;^{p^{(+)}}D^{\geq n}(L;R) \:.$$
Moreover,
 for $0<|w|\ll \delta\ll 1$, $\delta$ is a  regular value of $r$ with respect to
$\cS|_L$. This follows from the {\it $a_{f}$-condition}. 
If we apply Theorem~\ref{thm:plurisub} with $q=0$ to $r$ (or $-r$),
then we get by the (co)stalk formulae above:
$$i_{x}^{*}(\psi_{f}\cFb)[-1] \in \;^{p^{(+)}}D^{\leq n}(\{pt\};R)
\quad \text{or} \quad
i_{x}^{!}(\psi_{f}\cFb)[-1] \in \;^{p^{(+)}}D^{\geq n}(\{pt\};R) \:.$$
This proves our claim for a point stratum.

We reduce the general case to the first case by 
taking a {\it complex analytic normal slice} $N$ at $x$ (in some local embedding), with $\codim N=\dim S=s$.
Consider the cartesian diagram
\begin{displaymath} \begin{CD}
\{x\}=N\cap S @>\kappa_{x} >> N\cap \{f=0\} @>i' >> N\cap X \\
@VV k_{x} V  @VV k' V  @VV k V \\
S @> i_{S}  >> \{f=0\} @>i >> X \:.  
\end{CD} \end{displaymath}
Then we have 
$$k_{x}^{*}i_{S}^{*}(\psi_{f}\cFb)\simeq
\kappa_{x}^{*}k'^{*}(\psi_{f}\cFb) \:.$$
Since $\psi_{f}\cFb$ is weakly constructible with respect to
the induced stratification of $\{f=0\}$, we get by the {\it base change property} (\ref{bcs}):

$$k_{x}^{*}i_{S}^{!}(\psi_{f}\cFb)\simeq
\kappa_{x}^{!}k'^{*}(\psi_{f}\cFb) \:.$$
But we also have by Proposition \ref{bc-nv}  the {\it base change isomorphism}

$$k'^{*}(\psi_{f}\cFb)\simeq \psi_{f'}(k^{*}\cFb) \:,
\quad \text{with $f':=f\circ k$.} $$
Then the claim follows from the first case 
for $f',\, X':=N\cap X$ and $k^{*}\cFb$, since 
$N$ is transversal to $\cS$ near $x$, with $\codim N=\dim S=s$,
so that by Proposition \ref{prop-tr}:
$$k^{*}\cFb \in \;^{p^{(+)}}D^{\leq n-s}(X';R) \quad \text{or} \quad 
k^{*}\cFb\in \;^{p^{(+)}}D^{\geq n-s}(X';R) \:.$$
This completes the proof.
\end{proof}

\br We get in particular 
that the (shifted) {\it nearby and vanishing cycle functors}
${^p}\psi_f:=\psi_{f}[-1]$ and ${^p}\varphi_f:=\varphi_{f}[-1]$ are {\it t-exact} functors with respect to the {\it perverse t-structure}.

This is a result of Gabber in the algebraic
context for $l$-adic cohomology (unpublished, but compare with 
\cite[Proposition 4.4.2]{BBD}, \cite[Theorem 1.2]{Bry} and \cite[Corollary 4.5, 4.6]{Illu}). In the complex analytic
context it is due to Goresky-MacPherson \cite[Theorem 6.5]{GM3}, and 
\cite[p.222, Corollary 6.13.6, p.224 6.A.5]{GMs}. 
For another proof of
this classical case see \cite[Corollary 10.3.11, 10.3.13]{KS}.

But it does not seem to be well known that the result also applies for $R$ a Dedekind domain to the {\it dual t-structure}.
In particular, the functors ${^p}\psi_f=\psi_{f}[-1]$ and ${^p}\varphi_f=\varphi_{f}[-1]$  preserve then \index{strongly perverse sheaf}
{\it strongly perverse} sheaves.
\er

\bex
If $X$ is a pure-dimensional complex algebraic (or analytic) variety satisfying $\r(X,R)=\dim X$ for $R$ a Dedekind domain, then Proposition \ref{p28} and the above Remark yield that ${^p}\psi_f R_X[\dim X]$ and ${^p}\varphi_f R_X[\dim X]$ are strongly perverse sheaves on $X_0$. Therefore  
these perverse sheaves have torsion-free costalks in the lowest possible degree (cf. Corollary \ref{cfield}).
\eex

\bex
Let $f\colon (X,0) \to (\bC,0)$ be a nonconstant holomorphic function germ defined on a pure $(n+1)$-dimensional complex singularity germ contained in some ambient $(\bC^N,0)$. 
Denote  by $M_{f,0}$ the \index{Milnor fiber}
Milnor fiber of the singularity at the origin in $X_0=f^{-1}(0)$. Let $\Sigma:=\Sing_{\cS}(f)$ be the stratified singular locus of $f$ with respect to a fixed Whitney stratification $\cS$ of $X$, and set $r:=\dim_0 \Sigma$.
Let $R$ be a Dedekind domain and assume that $\r(X,R)=n+1$. The support condition for the perverse sheaf ${^p}\varphi_f R_X[n+1]$ (which is supported on $\Sigma$) yields that the only possibly non-trivial \index{reduced cohomology}
reduced cohomology $\widetilde{H}^k(M_{f,0};R)$ of $M_{f,0}$ is concentrated in degrees $n-r \leq k \leq n$. Moreover, as shown, e.g.,  in \cite[Theorem 3.4(d)]{MPT} or \cite[Example 6.0.12]{Sc}, the lowest (possibly) non-trivial module 
$\widetilde{H}^{n-r}(M_{f,0};R)$ is torsion-free.
\eex

For more applications of the fact that vanishing and nearby cycle functors preserve strongly perverse sheaves see \cite{Sc}, e.g., \cite[Example 6.0.14]{Sc} for applications to {\it local Lefschetz and vanishing theorems}, as well as  \cite{MP} (for the study of vanishing cohomology of complex projective hypersurfaces) and \cite{MPT} (for understanding the Milnor fiber cohomology).\\

Using Theorem \ref{thm:nearcycle} and Corollary \ref{cor:NMD}, together with Example \ref{isostr} and Remark \ref{rNMD}, one gets the following characterization (see, e.g., \cite[Corollary 6.0.7]{Sc}).

\bc\label{cor:vancycle3} 
Let $X$ be a complex algebraic (or analytic) variety.
Assume $\cFb\in D^b(X;R)$ is weakly constructible (resp., $\cFb$ is constructible with $R$ a Dedekind domain, in case we want to use the dual perverse t-structure). 
Then we have for\\ ``$?$ given by $\,\leq$'' or ``$?$ given by $\,\geq$'':
$$\cFb\in\: ^{p^{(+)}}D^{\, ?\, n}(X;R) \quad \Leftrightarrow 
\quad
\bigl((\varphi_{f}\cFb)[-1]\bigr)_{x} \in \: ^{p^{(+)}}D^{\, ?\, n}(\{pt\};R) $$
for  all holomorphic function germs $f\colon  (X,x)\to (\bC,0)$ such that $x$ is 
{\it isolated in the support of} $\varphi_{f}(\cFb)$.
 \ec

As an application we get once more the \index{effective characteristic cycles} {\it effectivity of characteristic cycles of perverse sheaves}, this time via {\it vanishing cycles}.
Assume $X$ is a complex manifold and $R$ is a field. Let $\cFb \in D^b_{\cS-c}(X;R)$ be a constructible complex with respect to some Whitney stratification 
$\cS$. Consider the characteristic cycle of $\cFb$, i.e., 
$$CC(\cFb)=\sum_{S \in \cS} m(S) \cdot \overline{T_{S}^*X},$$
with multiplicities $m(S)$ given as in Definition \ref{CC} in terms of the normal Morse data of $\cFb$. 
Here, we recall the calculation of $CC(\cFb)$ in terms of vanishing cycles.

For a stratum $S \in \cS$, let $x \in S$ be a point, and let $g\colon (X,x)\to (\bC,0)$ be a holomorphic function germ at $x$ such that $dg_x$ is a non-degenerate covector and $x$ is a complex Morse critical point of $g\vert_S$ (i.e., $dg(X)$ intersects $T_{S}^*X$ transversally at $dg_x$). Then, as in Example \ref{ex1337},  $Re(g)\vert_{S}$ has a classical Morse critical point at $x$ of Morse index equal to $\dim S$, and the multiplicities $m(S)$ of $CC(\cFb)$ can be computed from Theorem \ref{thm:SMTB} and Corollary \ref{van-stalk}:
\be\label{mul}
m(S)=\chi(({^p\varphi}_g \cFb)_x).
\ee
Then we get by Example  \ref{isostr} and Corollary \ref{cor:vancycle3}:
\bc
If, in the above notations, $\cFb$ is a {\it perverse sheaf}, then its characteristic cycle is  {\it effective}, i.e., $m(S) \geq 0$ for all strata $S$ in $\cS$.
\ec

Let us come back to the general context for finishing this section with  the following applications of Theorem \ref{thm:nearcycle}
(see, e.g., \cite[Proposition 6.0.2]{Sc}).

\bp
\label{prop:subspace}
Let $i\colon  Y\hookrightarrow X$ be the inclusion
of a closed complex algebraic (or analytic) subset, with $j\colon  U:=X\backslash Y\hookrightarrow X$ the
inclusion of the open complement. Suppose that $Y$ can locally be 
described (at each point $x\in Y$) as the common zero-set of at most
$k$ algebraic (or holomorphic) functions on $X$ ($k\geq 1$). 

Assume $\cFb\in D^b(X;R)$ is weakly constructible (resp., $\cFb$ is constructible with $R$ a Dedekind domain, 
in case we want to use the dual perverse t-structure).  
Then one has:
\begin{enumerate}
\item $\cFb \in\: ^{p^{(+)}}D^{\leq n}(X;R) \Rightarrow 
i^{!}\cFb [k] \in\: ^{p^{(+)}}D^{\leq n}(Y;R)$.
\item $ \cFb \in\: ^{p^{(+)}}D^{\geq n}(X;R) \Rightarrow 
i^{*}\cFb [-k]\in\: ^{p^{(+)}}D^{\geq n}(Y;R)$.
\item 
$j^{*}\cFb \in\: ^{p^{(+)}}D^{\leq n}(U;R) \Rightarrow 
i^{*}Rj_{*}j^{*}\cFb [k-1] \in\: ^{p^{(+)}}D^{\leq n}(Y;R)$ and\\
$Rj_{*}j^{*}\cFb [k-1] \in\: ^{p^{(+)}}D^{\leq n}(X;R)\:.$
\item $j^{*}\cFb  \in\: ^{p^{(+)}}D^{\geq n}(U;R) \Rightarrow i^{!}Rj_{!}j^{*}\cFb [-(k-1)]
 \in\:^{p^{(+)}}D^{\geq n}(Y;R)$ and\\
$Rj_{!}j^{*}\cFb [-(k-1)] \in\: ^{p^{(+)}}D^{\geq n}(X;R) \:.$
\end{enumerate} \ep

These are local results, 
so that we can assume that $Y$ is the common zero-set of $k'$ algebraic (or holomorphic)
functions (with $1\leq k'\leq k$). Then the claim for (1.) or (2.) follows by induction
from the case $k'=1$, which is a direct application of Theorem~\ref{thm:nearcycle}
and the distinguished triangle~(\ref{canmo}) or~(\ref{eq:dtrvar}).
Finally (3.) or (4.) is, by $i^{*}Rj_{*}j^{*}\simeq i^{!}Rj_{!}j^{*}[1]$ (see, e.g., \cite[(1.4.6.4)]{BBD}), 
a special case of (1.) or (2.).

\bex[Purity] \label{purity} \index{purity}
Let $i\colon Y\hookrightarrow X$ be the inclusion
of a closed complex algebraic (or analytic) subset,
which can locally be described as the common zero-set of at most
$k$ algebraic (or holomorphic) functions on $X$. Then $i^{!}\bD_{X}\simeq \bD_{Y}$ implies by 
Proposition~\ref{prop:subspace} the following estimate of the \index{rectified homological depth}
{\it rectified homological depth}:
\begin{equation} \label{eq:purity}
\r(X,R) \geq n \:\Rightarrow
\: \r(Y,R \geq n-k  \:.
\end{equation}
Especially for $X$ smooth of pure dimension $n$, we get (using also $i^*R_X=R_Y$)
$$ \r(Y,R) \geq n-k = \dim Y \geq \r(Y,R) \quad \text{and} \quad R_Y[\dim Y]\in  \: ^{p}D^{\geq 0}(Y;R)\:,$$
if $Y$ is locally a {\it set-theoretical complete intersection} of codimension $k$
in $X$. In particular
$$\r(Y,R)=\dim Y \quad \text{and} \quad \text{$R_Y[\dim Y]$ is a perverse sheaf}$$
for $Y$ a {\it pure-dimensional local complete intersection}. \index{local complete intersection}
\eex

Compare also with \cite[Theorem 3.2.1, Corollary.3.2.2]{HL} for the corresponding
homotopy results.


\subsection{Thom-Sebastiani for vanishing cycles}
In this section, we state a Thom-Sebastiani result for vanishing cycles, generalizing \cite{TS} to functions defined on singular ambient spaces, with arbitrary critical loci, and with arbitrary weakly constructible sheaf coefficients. For complete details, see \cite{Ma3} and also \cite[Corollary 1.3.4]{Sc}. 

Let $f\colon X \to \bC$ and $g \colon Y \to \bC$ be complex algebraic (or analytic) functions. Let $pr_1$ and $pr_2$ denote the projections of $X\times Y$ onto $X$ and $Y$, respectively. Consider the function
$$f \boxtimes g:=f \circ pr_1 + g \circ pr_2\colon X \times Y \to \bC.$$
The goal is to express the vanishing cycle functor $\varphi_{f \boxtimes g}$ in terms of the corresponding functors $\varphi_{f}$ and $\varphi_{g}$ for $f$ and, resp., $g$. 

We let $V(f)=\{f=0\}$, and similarly for $V(g)$ and $V(f \boxtimes g)$. Denote by 
$\ell$  
the inclusion of $V(f) \times V(g)$ into $V(f \boxtimes g)$. With these notations, one has the following result.
\bt\label{tsvan} \index{Thom-Sebastiani Theorem for vanishing cycles}
For $\cFb \in D^b_{wc}(X;R)$ and $\cGb \in D^b_{wc}(Y;R)$, there is a natural isomorphism
\be\label{tsf1}
\ell^*{{^p}\varphi}_{f \boxtimes g} ( \cFb \overset{L}{\boxtimes} \cGb ) \simeq {{^p}\varphi}_{f} \cFb \overset{L}{\boxtimes} {{^p}\varphi}_{g} \cGb
\ee
commuting with the corresponding \index{monodromy} monodromies. 

Moreover, if $p=(x,y) \in X \times Y$ is such that $f(x)=0$ and $g(y)=0$, then, in an open neighborhood of $p$, the complex ${{^p}\varphi}_{f \boxtimes g} ( \cFb \overset{L}{\boxtimes} \cGb )$ has support contained in $V(f) \times V(g)$, and, in every open set in which such a containment holds, there are natural isomorphisms
\be\label{tsf2}
{{^p}\varphi}_{f \boxtimes g} ( \cFb \overset{L}{\boxtimes} \cGb ) \simeq \ell_! ( {{^p}\varphi}_{f} \cFb \overset{L}{\boxtimes} {{^p}\varphi}_{g} \cGb ) \simeq \ell_* ( {{^p}\varphi}_{f} \cFb \overset{L}{\boxtimes} {{^p}\varphi}_{g} \cGb ). 
\ee 
\et

\bc In the notations of the above theorem and with integer coefficients, there is an isomorphism 
\begin{multline}\label{tsigen}
\wti{H}^{i-1}(M_{f \boxtimes g,p}) \cong \bigoplus_{a+b=i} \left( \wti{H}^{a-1}(M_{f,pr_1(p)}) \otimes \wti{H}^{b-1}(M_{g,pr_2(p)}) \right) \\ \oplus \bigoplus_{c+d=i+1} \Tor\left(\wti{H}^{c-1}(M_{f,pr_1(p)}), \wti{H}^{d-1}(M_{g,pr_2(p)})\right),
\end{multline}
where $M_{f,x}$ denotes as usual the Milnor fiber of a function $f$ at $x$, and similarly for $M_{g,y}$.
\ec

\br Note that the Thom-Sebastiani Theorem \ref{tsvan}
 implies directly the {\it multiplicativity of normal Morse data} in terms of vanishing cycles with respect to {\it external products} as mentioned in
Example \ref{NMDproduct}.
\er

\bex [Brieskorn singularities and intersection cohomology] \index{Brieskorn singularities}
For $i=1,\ldots, n$, consider a $\bC$-local system $\cL_i$ of rank $r_i$ on $\bC^*$, with monodromy automorphism $h_i$, and denote the corresponding intersection cohomology complex on $\bC$ by $IC_{\bC}(\cL_i)$. The complex $IC_{\bC}(\cL_i)$ agrees with $\cL_i[1]$ on $\bC^*$, and has stalk cohomology at the origin concentrated in degree $-1$, where it is isomorphic to $\ker(id - h_i)$.
For positive integers $a_i$, consider the functions $f_i(x)=x^{a_i}$ on $\bC$. The complex ${{^p}\varphi}_{f_i} IC_{\bC}(\cL_i)$ is a perverse sheaf supported only at $0$; therefore, ${{^p}\varphi}_{f_i} IC_{\bC}(\cL_i)$ is non-zero only in degree zero, where it has dimension $a_i r_i - \dim \ker(id-h_i)$. 

Next, consider the $\bC$-local system  $\cL_1 {\boxtimes} \cdots {\boxtimes} \cL_n$ on $(\bC^*)^n$ with monodromy automorphism $h:=\boxtimes_{i=1}^n h_i$, and note that, as in Example \ref{ex1230},
$$IC_{\bC}(\cL_1) \overset{L}{\boxtimes} \cdots \overset{L}{\boxtimes} IC_{\bC}(\cL_n) \simeq IC_{\bC^n}(\cL_1 {\boxtimes} \cdots {\boxtimes} \cL_n).$$
The perverse sheaf $${{^p}\varphi}_{x_1^{a_1}+\cdots + x_n^{a_n}} IC_{\bC^n}(\cL_1 {\boxtimes} \cdots {\boxtimes} \cL_n)$$ is supported only at the origin, and hence is concentrated only in degree zero. In degree zero,  it can be seen by iterating the Thom-Sebastiani isomorphism that it has dimension equal to $$\prod_i \left( a_i r_i - \dim \ker(id-h_i) \right).$$
In the special case when $r_i=1$ and $h_i=1$ for all $i$, the above calculation recovers the classical result stating that the Milnor number of the isolated singularity at the origin of $x_1^{a_1}+\cdots + x_n^{a_n}=0$ is $\prod_i (a_i-1)$. 
\eex


\subsection{Gluing perverse sheaves via vanishing cycles} 
Assume $R=\bC$ (or, more generally, $R$ is an algebraically closed field; but only the case $R=\bC$ also nicely fits with the corresponding theory of (regular) holonomic D-modules). Let $f\colon  X\to \bC$ be a complex algebraic (or analytic) morphism with corresponding nearby and vanishing cycle functors $\psi_f$, $\varphi_f$. Recall that these two functors come equipped with monodromy automorphisms, both of which are denoted here by $h$.
For $\cFb \in D^b_c(X;\bC)$, the morphism $$can\colon \psi_f \cF^\bullet \lra \varphi_f \cF^\bullet$$ of (\ref{canmo}) is called the \index{canonical morphism} {\it canonical morphism}, and it is compatible with monodromy. Similarly for the \index{variation morphism} {\it variation morphism}
 $$var \colon \varphi_f \cF^\bullet \lra \psi_f \cF^\bullet$$ 
of (\ref{eq:dtrvar}).
In the above notations we also have
 \be can \circ var= h -id  \quad \text{and} \quad var \circ can= h - id \:.\ee
 
 The monodromy automorphisms acting on the nearby and vanishing cycle functors have Jordan decompositions $$h=h_u \circ h_s=h_s \circ h_u,$$ where $h_s$ is {\it semi-simple} (and locally of finite order) and  $h_u$ is  {\it unipotent}.  For any $\lambda \in \bC$ and $\cFb \in D^b_c(X;\bC)$, denote by 
$\psi_{f,\lambda}\cFb$ the generalized $\lambda$-eigenspace for $h$, and similarly for  $\varphi_{f,\lambda}\cFb$. By the definition of vanishing cycles, the canonical morphism $can$ induces morphisms $$can\colon \psi_{f,\lambda}\cFb \lra \varphi_{f,\lambda}\cFb,$$ which (since the monodromy acts trivially on $i^*\cFb$) are isomorphisms for $\lambda \neq 1$, and there is a distinguished triangle 
\be
i^*\cFb \overset{{sp}}{\lra} \psi_{f,1} \cFb \overset{can}{\lra} \varphi_{f,1} \cFb \overset{[1]}{\lra} \ .
\ee
There are decompositions
\be\label{678} \psi_{f}= \psi_{f,1} \oplus  \psi_{f,\neq1} \ \ {\rm and} \ \ \varphi_{f}= \varphi_{f,1} \oplus  \varphi_{f,\neq1} \ee so that $h_s=1$ on $\psi_{f,1}$ and $\varphi_{f,1}$, and $h_s$ has no $1$-eigenspace on $\psi_{f,\neq1}$ and $\varphi_{f,\neq1}$. Moreover, $can\colon \psi_{f,\neq1} \to \varphi_{f,\neq1}$ and $var\colon  \varphi_{f,\neq1} \to \psi_{f,\neq1}$ are isomorphisms.\\

The canonical and variation morphisms play an important role in the following \index{gluing of perverse sheaves} {\it gluing} of perverse sheaves.\\

Let $X$ be a complex  algebraic variety, with $i\colon  Z\hookrightarrow X$ the inclusion of a closed algebraic subvariety, and $j\colon  U\hookrightarrow X$ the inclusion
of the open complement $U:=X\backslash Z$. A natural question to address is if one can ``glue'' the categories $Perv(Z)$ and $Perv(U)$ to recover the category $Perv(X)$ of perverse sheaves on $X$. We discuss here only the case when $Z$ is a hypersurface, but see also \cite{Ve2} for a more general setup. 
 The gluing procedure, due to Beilinson \cite{Bei} and Deligne-Verdier \cite{Ve2}, establishes an equivalence of categories between perverse sheaves on the algebraic variety $X$ and a pair of perverse sheaves, one on $Z$, the other on $U$, together with a gluing datum.
A similar method was used by M. Saito for constructing his mixed Hodge modules \cite{Sa2}.

As a warm-up situation, consider $X=\bC$ with coordinate function $z$, $Z=\{0\}$ and $U=\bC^*$. Let $\cFb$ be a $\bC$-perverse sheaf  on $X$. One can form the diagram
$${^p\psi}_z \cFb
\overset{can}{\underset{var}\rightleftarrows}
{^p\varphi}_z \cFb$$
whose objects are perverse sheaves on $Z=\{0\}$, i.e., complex vector spaces. This leads to the following elementary description of the category of perverse sheaves on $\bC$, cf. \cite{Ve2}.
\bp
\label{DV}
The category of perverse sheaves on $\bC$ which are locally constant on $\bC^*$ is equivalent to the category of quivers (that is, diagrams of vector spaces) of the form
$${\psi} 
\overset{c}{\underset{v}\rightleftarrows}
\varphi$$
with $\psi, \varphi$ finite dimensional vector spaces, and $id+c \circ v$, $id+v \circ c$ invertible.
\ep

\bex\label{exp}
Let $\cL$ be a local system on $\bC^*$ with stalk $V$ and monodromy $h\colon V \to V$. The perverse sheaf $j_*\cL[1]$ (e.g., see Example \ref{ex221}) corresponds to 
$$V \overset{c}{\underset{v}\rightleftarrows} V/{\ker(h-id)},$$ 
where $c$ is the projection and $v$ is induced by $h-id$. Thus a quiver \index{quiver}
$${\psi} \overset{c}{\underset{v}\rightleftarrows} \varphi$$ 
with $c$ surjective arises from $j_*\cL[1]$, where $\cL_1=\psi$ is the stalk of $\cL$ and $h=id+v\circ c$.
\eex

More generally, let $f$ be a regular function on a smooth algebraic variety $X$, with $Z=f^{-1}(0)$ and $U=X\setminus Z$. Let $Perv(U,Z)_{gl}$ be the category whose objects are $(\cAb,\cBb,c,v)$, with $\cAb \in Perv(U)$, $\cBb \in Perv(Z)$, $c\in \Hom({^p\psi}_{f,1}\cAb,\cBb)$, $v \in \Hom(\cBb,{^p\psi}_{f,1}\cAb)$, and so that $id+v \circ c$ is invertible. Then one has the following.
\bt[Beilinson \cite{Bei}, Deligne-Verdier \cite{Ve2}]\label{tgl} \index{gluing of perverse sheaves}
There is an equivalence of categories
$$Perv(X) \simeq Perv(U,Z)_{gl}$$
defined by:
$$\cFb \mapsto (\cFb\vert_U,{^p\varphi}_{f,1}\cFb,can,var).$$
\et

See also \cite[Section 3.3, 3.4]{SW} for some other examples of quiver descriptions of perverse sheaves, especially also for the corresponding description of Verdier duality on the quiver side.


\subsection{Other applications}
There are many more applications of the nearby and vanishing cycles functors than we can mention in an expository paper, e.g., see \cite{Di}, \cite{Sc}, \cite{KS}, \cite{M2} and the references therein. Let us only mention here that vanishing cycles play an important role in the theory of characteristic classes for singular hypersurfaces (e.g., see \cite{CMSS}, \cite{MSS}, \cite{MSS2}, \cite{Sc2}), which have recently seen applications in birational geometry (cf. \cite{MSS2}), and they have also found applications to concrete optimization problems in applied algebraic geometry and algebraic statistics (e.g., see \cite{MRW2}, \cite{MRW3} and the survey \cite{M3}).

\medskip 

We finish this section with an application  to the calculus of \index{Grothendieck group}
Grothendieck groups of constructible sheaves and  their relation to the theory of  constructible functions (see, e.g., \cite[Section 6.0.6]{Sc} for a more general version).\\

\noindent
Let $f\colon  X\to \C$ be an algebraic (or analytic) morphism  with
$i\colon X_{0}:=\{f=0\}\hookrightarrow X$ and $j\colon U:=\{f \neq 0\}\hookrightarrow X$ the inclusions.
Then  one has the distinguished triangle~(\ref{eq:dtrwang})
$$\begin{CD}
i^{*}Rj_{*}j^{*}\cFb @>>> \psi_{f}\cFb @> h_{f} -id  >> \psi_{f}\cFb @>[1]>> 
\end{CD} $$
for any $\cFb\in D^b_c(X;R)$. But then also $\psi_f\cFb\in D^b_c(X_0;R)$ by
Corollary \ref{nv-weak} and Corollary \ref{nvc}, so that

\be \label{eq:K0}
[i^{*}Rj_{*}j^{*}\cFb]= \:0 \:\in K_{0}(D^b_c(X_0;R))
\ee
for the corresponding class in the 
Grothendieck group of constructible sheaf complexes on $X_0$.

Let more generally $X_{0}:=\{f_{i}=0,\:i=1,..,k+1\}$ be the common zero-set of
the algebraic (or holomorphic)  functions $f_{i}\colon X\to \bC$ ($i=1,\dots,k+1$).
Let $j\colon U:=X\backslash X_{0}\hookrightarrow X$ be again the inclusion.
Then the equality~(\ref{eq:K0}) is also true. This follows by induction from \index{Mayer-Vietoris triangle}
a {\it Mayer-Vietoris triangle} (see, e.g., \cite[p.114, (2.6.28)]{KS}, 
with $R\Gamma_{U}(\cdot)\simeq Rj_{*}j^{*}$) associated to the
open covering $\{U_{1},U_{2}\}$ of $U$, with 
$$U_{1}:=\{f_{k+1}\neq 0\} \quad \text{and} \quad
U_{2}:=X\backslash \{f_{i}=0,\:i=1,\dots,k\} \:.
$$
Note that 
$$U_{1}\cap U_{2}=
X\backslash \{f_{i}\cdot f_{k+1}=0,\:i=1,\dots,k\} $$
is the complement of the
common zero-set of $k$ algebraic (or holomorphic)  functions.

Let $i\colon  X_{0}\hookrightarrow X$ be the inclusion. By the distinguished triangle
$$\begin{CD}
Rj_{!}j^{*}\cFb @>>> Rj_{*}j^{*}\cFb @>>> i_{*}i^{*}Rj_{*}j^{*}\cFb @>[1]>> 
\end{CD} $$
we get
\be \label{eq:j*=j!K0}
Rj_{!}j^{*}=Rj_{*}j^{*}:\: K_{0}(D^b_c(X;R))\to K_{0}(D^b_c(X;R))\:.
\ee
Similarly, the distinguished triangle
$$ \begin{CD}
i^{!}\cFb @>>> i^{*} \cFb@>>> i^{*}Rj_{*}j^{*}\cFb @>[1]>> 
\end{CD} $$
implies
\be \label{eq:i*=i!K0}
i^{!}=i^{*}:\: K_{0}(D^b_c(X;R)) \to K_{0}(D^b_c(X_0;R))\:.
\ee

Using again a \index{Mayer-Vietoris triangle}
{\it Mayer-Vietoris triangle} (see, e.g., \cite[p.114, (2.6.28)]{KS}, 
with $R\Gamma_{U}(\cdot)\simeq Rj_{*}j^{*}$) associated to an
open affine covering, these equalities (\ref{eq:j*=j!K0}) and (\ref{eq:i*=i!K0}) are then available in the general complex algebraic context.
Consider in addition a proper morphism $h\colon X\to Y$. Then $Rh_{!}=Rh_{*}$ maps $D^b_c(X;R)$ into $D^b_c(X;R)$
and we get the following equality for the induced maps on the level of
Grothendieck groups:
\be \label{eq:*=!K0}
R(h\circ j)_{!}j^{*}=R(h\circ j)_{*}j^{*}:\: K_{0}(D^b_c(X;R)) \to K_{0}(D^b_c(Y;R))
\ee
in the complex algebraic context.
Let $f\colon X\to Y$ be a morphism of complex algebraic varieties.
By a theorem of Nagata (see, e.g., \cite{Luet, Con}) 
this can be (partially)  compactified to $f=\bar{f}\circ j$, with $j\colon  X\hookrightarrow \bar{X}$
an open inclusion and $\bar{f}\colon  \bar{X}\to Y$ a proper morphism of
complex algebraic varieties. 
Then the functors $Rf_{!}, Rf_{*}$ preserve algebraically constructible complexes
and one gets the equality (compare also with  \cite[p.210]{La} and \cite{Vir}): 
\be \label{eq:alg*=!K0}
Rf_{!}=Rf_{*}:\: K_{0}(D^b_c(X;R)) \to 
K_{0}(D^b_c(Y;R))\:,
\ee
which also implies for the induced group homomorphisms of  complex algebraically constructible functions
the equality
\be \label{eq:alg*=!CF}
f_{!}=f_{*}:\; CF(X) \to CF(X) \:.
\ee

Let us now consider the complex analytic context, with
$CF(X)$ the abelian group of 
complex analytically constructible  functions on $X$.
\bex\label{ex:analconstfct}
Let $i\colon X_{0}\hookrightarrow X$ be the inclusion
of a closed complex analytic subset, with $j\colon U:=X\backslash X_{0}\hookrightarrow X$ the
inclusion of the open complement. Fix also a proper holomorphic map $f\colon X\to
Y$ of complex analytic varieties. Then we get the following equalities 
for the corresponding group homomorphisms on the level of 
constructible functions:
\begin{enumerate}
\item $j_{!}j^{*}=j_{*}j^{*}:\; CF(X) \to CF(X)$,
especially $i^{*}j_{*}j^{*}=0$.
\item $i^{!}=i^{*}:\; CF(X) \to CF(X_{0})$.
\item $(f\circ j)_{!}j^{*}=(f\circ j)_{*}j^{*}:\; CF(X) \to CF(Y)$.
\end{enumerate}
\eex

Note that (1.) implies (3.). Moreover, (1.) and (2.) are local results so that we can assume
that $X_{0}$ is defined by finitely many holomorphic functions.
Then they follow from the corresponding result above on the level of
Grothendieck groups of constructible sheaf complexes.\\

A typically application is the following generalization of a classical result of
Sullivan for  $R$ a field and $X_0$ compact so that the constant map $f\colon X_0\to pt$ is proper. Then

\be \label{eq:link}
\chi(X_{0},i^{*}Rj_{*}j^{*}\cFb) = f_{*}i^{*}j_{*}j^{*}(\chi_{stalk}(\cFb)) = 0 
\ee
for any $\cFb\in D^b_c(X;R)$.

\br  $R\Gamma(X_{0}, i^{*}Rj_{*}j^{*}\cFb)$ calculates the \index{global
link cohomology} {\it global
link cohomology} of $\cFb$ as defined in \cite{DS}.
For example, for $X_{0}=\{x\}$ a point, one has
(for $0<\epsilon \ll 1$):
$$R\Gamma(X_{0}, i^{*}Rj_{*}j^{*}\cFb)\simeq 
R\Gamma(X\cap\{||z||=\epsilon\}, \cFb)
 \;,$$
with $||\cdot||$ defined by some local embedding $(X,x)\hookrightarrow
(\bC^{n},0)$. So, for $\cFb:=R_{X}$ the constant sheaf, we get that
the {\it link} \index{link}
$$X\cap\{||z||=\epsilon\} \quad \text{of $X$ in $x$}$$
has a {\it vanishing Euler characteristic}. 
For a proof of this classical result
of Sullivan \cite{Su}, based on the {\it Milnor fibration theorem}, compare with
\cite[Proposition 4.1]{BV}. 
\er


\section{Intersection cohomology, the decomposition theorem, applications}\label{sec:ih}

In this section, we overview properties of the intersection cohomology groups of complex algebraic varieties, which generalize the corresponding features of the cohomology groups of smooth varieties. These properties, consisting of Poincar\'e duality, Lefschetz theorems and the decomposition theorem, are collectively termed the {\it K\"ahler package} for intersection cohomology. 
In this section we work over $R=\bQ$.

Recall the following.
\bd
The (compactly supported) \index{intersection cohomology group} {\it intersection cohomology groups} of a pure complex $n$-dimensional complex algebraic variety $X$ are defined as:
$$IH^k(X;R) :={H}^{k-n}(X;IC_X) \ , \ \ 
IH_c^k(X;R) :={H}_c^{k-n}(X;IC_X).$$
\ed

First note that Poincar\'e duality for the intersection cohomology groups of $X$ is an immediate corollary of the self-duality of $IC_X$ (cf. Proposition \ref{sed}). Specifically, if $X$ is a pure $n$-dimensional complex algebraic variety, there is an isomorphism
\be
IH^{2n-k}(X;\bQ) \simeq IH^k_c(X;\bQ)^\vee,
\ee 
for any integer $k$.

Note also that if $X$ is a pure $n$-dimensional complex algebraic variety, there exists a natural morphism \be\label{natm} \alpha_X\colon {\Q}_X[n] \lra IC_X,\ee 
extending the natural quasi-isomorphism on the smooth locus of $X$. Moreover, $ \alpha_X$ becomes a quasi-isomorphism is $X$ is assumed to be a rational homology manifold (see, e.g., \cite[Theorem 6.6.3]{M}). 
Applying the hypercohomology functor to \eqref{natm} yields a morphism 
$$H^k(X;\bQ)\lra IH^k(X;\bQ),$$ which is an isomorphism if $X$ is a rational homology manifold.


\subsection{Lefschetz type results for intersection cohomology}
Lefschetz type results for intersection cohomology can be derived from sheaf theoretic statements about perverse sheaves. 
We begin with the following immediate consequence of the Artin vanishing Theorem \ref{Artinvan} for perverse sheaves: \index{weak Lefschetz theorem}
\bt[Weak Lefschetz Theorem for Perverse Sheaves]\label{wlps} 
If $X$ is a complex projective variety and $i\colon D \hookrightarrow X$ denotes the inclusion of a hyperplane section, then for every $\mathscr{F}^\bullet \in Perv(X;\bQ)$ the restriction map $$H^k(X;\cFb) \lra H^k(D;i^*\cFb)$$ is an isomorphism for $k<-1$ and is injective for $k=-1$.
\et

 \bex
 Let $X \subset \bC P^{n+1}$ be a complex projective hypersurface, with $i\colon D \hookrightarrow X$ the inclusion of a hyperplane section. By applying Theorem \ref{wlps} to the perverse sheaf $\cFb:={\bQ}_X[n]$ on $X$, one gets isomorphisms 
$$H^k(X;\bQ) \lra H^k(D;\bQ)$$ for all $k<n-1$ and a monomorphism for $k=n-1$.
 \eex

Theorem \ref{wlps} has the following important consequence (see also Example \ref{ex:weaklef}).
\bc[Lefschetz hyperplane section theorem for intersection cohomology]  
Let $X^n \subset \mathbb{C}{P}^N$ be a pure $n$-dimensional closed algebraic subvariety with a Whitney stratification $\mathscr{S}$. Let $H \subset \mathbb{C}{P}^N$ be a generic hyperplane (i.e., transversal to all strata of
 $\mathscr{S}$), with $i\colon D:=X\cap H \hookrightarrow X$ the inclusion of the corresponding hyperplane section. Then the natural homomorphism \[ IH^k(X;\bQ) \longrightarrow IH^k(D;\bQ) \] is an isomorphism for $0 \leq k \leq n-2$ and a monomorphism for $k= n-1$. 
\ec
\noindent (Indeed, by transversality, one gets that $i^*IC_X \simeq IC_{D}[1]$, e.g., see Example \ref{IC-transv}.)\\

The Hard Lefschetz theorem for intersection cohomology can also be deduced from a more general sheaf-theoretic statement. Let $f\colon X \to Y$ be a projective morphism and let $L \in H^2(X;\bQ)$ be the first Chern class of an $f$-ample line bundle on $X$. Then $L$ corresponds to a map of complexes $L\colon {\bQ}_X \to {\bQ}_X[2]$, which, after tensoring with $IC_X$ yields a map $$L\colon IC_X \lra IC_X[2].$$ This induces $L\colon Rf_*IC_X \to Rf_*IC_X[2]$, and after applying perverse cohomology one gets a map of perverse sheaves on $Y$:
$$L\colon {^p\cH}^{i}(Rf_*IC_X) \lra {^p\cH}^{i+2}(Rf_*IC_X).$$
Iterating, one gets maps of perverse sheaves 
$$L^i\colon {^p\cH}^{-i}(Rf_*IC_X) \lra {^p\cH}^{i}(Rf_*IC_X)$$ for every $i \geq 0$.
Then one has the following result, proved initially by positive characteristic methods \cite[Theorem 6.2.10]{BBD} (see also \cite{Sa1,Sa2}, or the more geometric approach of \cite{CM}):
\bt[Relative Hard Lefschetz] \label{hlih} \index{relative Hard Lefschetz theorem}
Let $f\colon X \to Y$ be a projective morphism of complex algebraic varieties with $X$ pure-dimensional, and let $L \in H^2(X;\bQ)$ be the first Chern class of an $f$-ample line bundle on $X$. For every $i > 0$, one has isomorphisms of perverse sheaves
$$L^i\colon {^p\cH}^{-i}(Rf_*IC_X) \overset{\simeq}{\lra} {^p\cH}^{i}(Rf_*IC_X).$$
\et

By taking $f$ in Theorem \ref{hlih} to be the constant map $f\colon X \to point$, one obtains as a consequence the Hard Lefschetz theorem for intersection cohomology groups:
\bc[Hard Lefschetz theorem for intersection cohomology]\label{hlic} \index{hard Lefschetz theorem for intersection cohomology}
Let $X$ be a complex projective variety of pure complex dimension $n$, with $L \in H^2(X;\bQ)$ the first Chern class of an ample line bundle on $X$. Then there are isomorphisms
\be\label{ihhl} L^i\colon  IH^{n-i}(X;\bQ) \overset{\simeq}{\lra} IH^{n+i}(X;\bQ)\ee
for every integer $i > 0$, induced by the cup product by $L^i$. In particular, the intersection cohomology Betti numbers of $X$ are unimodal, i.e., $\dim_\bQ IH^{i-2}(X;\bQ) \leq \dim_\bQ IH^{i}(X;\bQ)$ for all $i \leq n/2$.
\ec

\subsection{The decomposition theorem and immediate applications}
A great deal of information about intersection cohomology groups can be derived from the {\it BBD decomposition theorem}  \cite[Theorem 6.2.5]{BBD}, one of the most important results in the theory of perverse sheaves. It was conjectured by S. Gelfand and R. MacPherson, and proved soon after by Beilinson, Bernstein, Deligne and Gabber by reduction to positive characteristic. The proof given in \cite{BBD} ultimately rests on Deligne's proof of the Weil conjectures. Different proofs were given later on by M. Saito (as a consequence of his theory of mixed Hodge modules \cite{Sa1}) and, more recently, by de Cataldo and Migliorini \cite{CM} (involving only classical Hodge theory). A more general decomposition theorem (for semi-simple perverse sheaves) has been obtained by Mochizuki \cite{Mo1}, \cite{Mo2} (with substantial contributions of Sabbah \cite{Sab}), in relation to a conjecture of Kashiwara \cite{Kas}; see also \cite{dC2}.
For more topological versions of the decomposition theorem for self-dual perverse sheaves on the level of Witt and cobordism groups, see \cite{CS2, SW}.
In what follows, we explain the statement  of the decomposition theorem, together with a few applications. 

In its initial form of \cite{BBD}, the decomposition theorem calculates the derived pushforward of an $IC$-complex under a proper algebraic map. For simplicity, in this section we assume that all varieties are {\it irreducible}. 

Recall that every algebraic map $f\colon X \to Y$ of complex algebraic varieties can be {\it stratified}, i.e., there exist algebraic Whitney stratifications $\cS$ of $X$ and $\cT$ of $Y$ such that, given any connected component $T$ of a $\cT$-stratum on $Y$ one has the following properties:
\begin{itemize}
\item[(a)] \ $f^{-1}(T)$ is a union of connected components of strata of $\cS$, each of which is mapped submersively to $T$ by $f$;
\item[(b)] \ For every point $y \in T$, there is an Euclidean neighborhood $U$ of $y$ in $T$ and a stratum-preserving homeomorphism $h\colon U \times f^{-1}(y) \to f^{-1}(U)$ such that $f|_{f^{-1}(U)} \circ h$ is the projection to $U$.
\end{itemize}
Property $(b)$ is just {\it Thom's isotopy lemma}: for every stratum $T$ in $Y$, the restriction of $f$ over $T$, i.e., $f|_{f^{-1}(T)} \colon f^{-1}(T) \to T$, is a topologically locally trivial fibration. 
\bt[BBD decomposition theorem \cite{BBD}]\label{bbdgt}  \index{decomposition theorem}
Let $f \colon X \to Y$ be a proper map of complex
algebraic varieties. Then:
\begin{itemize}
\item[(i)] ({\it Decomposition}) \ There
is a (non-canonical) isomorphism in $D^b_c(Y;\bQ)$:
\be\label{bbdgi}
 Rf_\ast IC_X \simeq \bigoplus_i \ ^p\cH^i(Rf_\ast IC_X)[-i].
\ee
\item[(ii)] \index{semi-simplicity} ({\it Semi-simplicity}) \ 
Each $^p\cH^i(Rf_\ast IC_X)$ is a semi-simple object in $Perv(Y;\bQ)$, i.e., if $\cT$ is the set of connected
components of strata of $Y$ in a stratification of $f$, there is a
canonical isomorphism of perverse sheaves on $Y$:
\begin{equation}\label{ssp}
{^p\cH}^i( Rf_* IC_X) \simeq \bigoplus_{T \in \cT} \ IC_{\overline T}(\cL_{i,T})
\end{equation}
where the local systems 
$\cL_{i,T}$ 
are semi-simple.
\end{itemize}
 \et

\br
If $f$ is a projective submersion of smooth complex algebraic varieties, Theorem \ref{bbdgt} reduces to Deligne's decomposition theorem, see \cite{De1} and \cite[Theorem 4.2.6]{De2}.
\er

Standard facts in algebraic geometry (e.g., resolution of singularities and  Chow's lemma) reduce the proof of the BBD decomposition theorem (Theorem \ref{bbdgt})  to the case when $f\colon X \to Y$ is a projective morphism, with $X$ a smooth variety. 
 If the morphism $f$ is projective, then (\ref{bbdgi}) is a formal consequence of the relative Hard Lefschetz theorem (cf. \cite{De3}). So the heart of the BBD decomposition theorem consists of the semi-simplicity statement. After the above-mentioned reductions for $f$, the proof given in \cite{CM} is by induction on the pair of indices $(\dim Y, r(f))$, where $r(f)$ is the \index{degree of semi-smallness} {\it degree of semi-smallness} of $f$. The problem is then reduced to the case $r(f)=0$, i.e., that of a \index{semi-small map} {\it semi-small} map, when $Rf_*\bQ_X[\dim X] \simeq {^p\cH}^0(f_*\bQ_X[n])$ is perverse on $Y$. This case is handled via the non-degeneracy of a certain refined intersection pairing associated to the fibers over the most singular points of $f$.
 
One of the first consequences of the BBD decomposition theorem is that it gives a splitting of $IH^*(X;\bQ)$ in terms of twisted intersection cohomology groups of closures of strata in $Y$, namely:
\bc Under the assumptions and notations of Theorem \ref{bbdgt}, there is a splitting
 \be\label{ihspl} IH^j(X;\bQ) \simeq \bigoplus_{i \in \bZ} \bigoplus_{T \in \cT} \ IH^{j-\dim X +\dim T-i}({\overline T};\cL_{i,T}),\ee
 for every $j \in \bZ$.
 \ec
 
By applying Theorem \ref{bbdgt} to the case of a resolution of singularities one gets the following. 
\bc\label{ihh} The intersection cohomology of a complex algebraic variety is a direct summand of the cohomology of a resolution of singularities.
\ec

More generally, one has the following nice application of the BBDG decomposition theorem (e.g., see \cite[Section 4.5]{dC}) or \cite[Theorem 9.3.37]{M}):
\bt\label{dirs}
Let $f\colon X \to Y$ be a proper map of complex algebraic varieties, and let $Y':=f(X)$ be the image of $f$. Denote by $d=\dim X - \dim Y'$ the relative dimension of $f$. Then $IC_{Y'}[d]$ is a direct summand of $Rf_*IC_X$. In particular, $IH^j(Y';\bQ)$ is a direct summand of $IH^j(X;\bQ)$ for every integer $j$.
\et

An important application of the decomposition statement (\ref{bbdgi}) for $f\colon X \to Y$ is the $E_2$-degeneration of the corresponding {\it perverse Leray spectral sequence}:
$$E_2^{i,j} = H^i(Y;{^p\cH}^j(Rf_*IC_X)) \Longrightarrow H^{i+j}(Y;Rf_*IC_X)=IH^{i+j+\dim X}(X;\bQ).$$
This can be used to prove the following singular version of the classical global invariant cycle theorem for smooth projective maps, as well as a local version of it (see \cite[Corollary 6.2.8, Corollary 6.2.9]{BBD}):
\bt[Global and local invariant cycle theorems] \label{glli} \index{invariant cycle theorem}
Let $f\colon X \to Y$ be a proper map of complex algebraic varieties. Let $U \subseteq Y$ be a Zariski-open subset on which the sheaf $R^if_*IC_X$ is a local system. Then the following assertions hold:
\begin{itemize}
\item[(a)] ({\it Global}) The natural restriction map $$IH^i(X;\bQ) \lra H^0(U;R^if_*IC_X)$$ is surjective.
\item[(b)]	({\it Local}) Let $u \in U$ and $B_u \subset U$ be the intersection with a sufficiently small Euclidean ball (chosen with respect to a local embedding of $(Y,u)$ into a manifold) centered at $u$. Then the natural restriction/retraction map 
$$H^i(f^{-1}(u);IC_X) \simeq H^i(f^{-1}(B_u);IC_X) \lra H^0(B_u;R^if_*IC_X)$$
is surjective.
\end{itemize}
\et


\subsection{A recent application of the K\"ahler package for intersection cohomology}
In this section, we mention briefly a recent combinatorial application of the K\"ahler package for intersection cohomology. For more applications, the interested reader may consult \cite{CM2}, \cite{M} and the references therein.

Let $E=\{v_1,\cdots,v_n\}$ be a spanning subset of a $d$-dimensional complex vector space $V$, and let $w_i(E)$ be the number of $i$-dimensional subspaces spanned by subsets of $E$. 
In $1974$, Dowling and Wilson \cite{DW1,DW2} proposed the following conjecture (which is a special case of a more general conjecture for matroids; see \cite{Bra} for more recent developments):
\begin{conj}[Dowling--Wilson top-heavy conjecture]\label{thdw} \index{top-heavy conjecture}
	For all $i<d/2$ one has: 
 \be w_i(E) \leq w_{d-i}(E).\ee
\end{conj} 

\br
If $d=3$, de Bruijn--Erd\"os showed that $w_1(E)\leq w_2(E)$. More generally, Motzkin showed that $w_1(E)\leq w_{d-1}(E)$.
\er

Another conjecture concerning the numbers $w_i(E)$ was proposed by Rota \cite{R1,R2} in $1971$, and it can be formulated as follows:
\begin{conj}[Rota's unimodal conjecture]\label{Rota} \index{unimodal conjecture}
There is some $j$ so that
\be\label{Runim} w_0(E)\leq \cdots \leq w_{j-1}(E) \leq w_j(E) \geq w_{j+1}(E) \geq \cdots \geq w_d(E).\ee
\end{conj}

Huh-Wang \cite{HW} used the K\"ahler package on intersection cohomology groups to prove the Dowling-Wilson top-heavy conjecture, and the unimodality of the ``lower half'' of the sequence $\{w_i(E)\}$. Specifically, they showed the following:

\bt[Huh-Wang]\label{tthdw} For all $i<d/2$, the following properties hold:
\begin{itemize}
 \item[(a)]  (top heavy) \
 $ w_i(E) \leq w_{d-i}(E).$
 \item[(b)]  (unimodality) \ $w_i(E) \leq w_{i+1}(E).$\end{itemize}
\et

\begin{proof}[Sketch of proof] 
The key step in the proof is to show that the numbers $ w_i(E)$ are realized geometrically, i.e., 
there exists a complex $d$-dimensional projective variety $Y$ such that for every $0\leq i \leq d$ one has: \begin{center}
$H^{2i+1}(Y;\bQ)=0$ \ and \ $\dim_\bQ H^{2i}(Y;\bQ)=w_i(E).$\end{center}

To define the variety $Y$ one first uses $E=\{v_1,\cdots,v_n\}$ to construct a map $i_E\colon V^{\vee}\to \bC^n$ by regarding each $v_i\in E$ as a linear map on the dual vector space $V^{\vee}$. Precomposing $i_E$ with the open inclusion $\bC^n \hookrightarrow (\bC P^1)^n$ yields a map $$f\colon V^{\vee} \to (\bC P^1)^n.$$ Set $$Y:=\overline{\im(f)} \subset  (\bC P^1)^n.$$ Ardila-Boocher \cite{AB} showed that the variety $Y$ has an algebraic cell decomposition (i.e., it is paved by complex affine spaces), the number of $\bC^i$'s appearing in the decomposition of $Y$ being exactly $w_i(E)$. However, the space $Y$ is in this case highly singular, so its rational cohomology does not satisfy the K\"ahler package. Instead, one needs to use the corresponding intersection cohomology results. 

Next, note that for any complex projective variety $Y$ one has that $$\ker \left( H^i(Y;\bQ) \overset{\alpha}{\to} IH^i(Y;\bQ)\right) = W_{\leq i-1} H^i(Y;\bQ)$$ is the subspace of $H^i(Y;\bQ)$ consisting of classes of Deligne weight $\leq i-1$ (e.g., see \cite[Theorem 9.2]{Web}). Since the complex projective variety $Y$ constructed above has an algebraic cell decomposition, its cohomology $H^i(Y;\bQ)$ is pure of weight $i$. (This follows easily by induction using the fact that $H^i_c(\C^n;\Q)=0$ for $i \neq 2n$ and $H^{2n}_c(\C^n;\Q)=\C$ is pure of weight $2n$.)
Hence, the natural map $$\alpha\colon H^*(Y;\bQ) \to IH^*(Y;\bQ)$$ is injective in any degree $i$.

For $i<d/2$,  consider the following commutative diagram:
\[
\xymatrix{
H^{2i}(Y;\bQ) \ar[d]_{\smallsmile L^{d-2i}} \ar@{^(->}[r]^{\alpha} & IH^{2i}(Y;\bQ) \ar[d]^{\smallsmile L^{d-2i}}_{\simeq} \\
H^{2d-2i}(Y;\bQ) \ar@{^(->}[r]^{\alpha} & IH^{2d-2i}(Y;\bQ)  
}
\] 
where the right-hand vertical arrow is the Hard Lefschetz isomorphism for the intersection cohomology groups of $Y$ (see Corollary \ref{hlic}). Since the maps labelled by $\alpha$ are injective, it follows that
$$H^{2i}(Y;\bQ) \xrightarrow{\smallsmile L^{d-2i}} H^{2d-2i}(Y;\bQ)$$ is injective as well. In particular, 
$$w_i(E)=\dim_\bQ H^{2i}(Y;\bQ) \leq \dim_\bQ H^{2d-2i}(Y;\bQ)=w_{d-i}(E)$$ for every $i<d/2$, thus proving part $(a)$ of the theorem. 
Part $(b)$ follows similarly, by using the unimodality of the intersection cohomology Betti numbers (cf. Corollary \ref{hlic}).
\end{proof}


\section{Perverse sheaves on semi-abelian varieties}\label{semi-ab}
In this section, we survey recent developments in the study of perverse sheaves on semi-abelian varieties, with concrete applications in the study of homotopy types of complex algebraic manifolds (formulated in terms of their cohomology jump loci), as well as new topological characterizations of semi-abelian varieties. We begin by introducing and recalling the main ingredients needed to formulate our results. For complete details, the interested reader may consult  \cite{LMW3, LMW4, LMW5}.


\subsection{Cohomology jump loci} 
Let $X$ be a connected CW complex of finite type (e.g., a complex quasi-projective variety) with positive first Betti number, i.e., $b_1(X)>0$. The \index{character variety} {\it character variety} $\Char(X)$ of $X$ is the identity component of the moduli space of rank-one $\bC$-local systems on $X$, i.e., 
\begin{center}
$\Char(X):= 
\Hom (H_1(X;\bZ)/\text{Torsion}, \bC^*)\simeq (\bC^*)^{b_1(X)}.$
\end{center}

\bd The \index{cohomology jump loci} {\it $i$-th cohomology jump locus of $X$} is defined as: 
$${\sV^{i}(X)=\lbrace \rho\in \Char(X) \mid  H^{i}(X;L_{\rho})\neq 0 \rbrace},$$ 
where $L_{\rho}$ is the unique rank-one $\bC$-local system  on $X$ associated to the representation $\rho\in \Char(X) $.\ed

The jump loci $\sV^i(X)$ are closed subvarieties of $\Char(X)$ and homotopy invariants of $X$.  For projective varieties, they can be seen as topological counterparts of the Green-Lazarsfeld jump loci of topologically trivial line bundles  \cite{GLa,GrLa}. Cohomology jump loci emerged from work of Novikov on Morse theory for closed $1$-forms on manifolds, and provide a unifying framework for the study of a host of questions concerning homotopy types of complex algebraic varieties. 


 \subsection{Jump loci via constructible complexes}
The Albanese map construction allows one to interpret the cohomology jump loci of a smooth connected complex quasi-projective variety as cohomology jump loci of certain constructible complexes of sheaves (or even of perverse sheaves, if the Albanese map is proper) on a semi-abelian variety. This motivates the study of cohomology jump loci of such complexes, with a view towards characterizing important classes of objects (such as perverse sheaves) on semi-abelian varieties.
 
An \index{abelian variety} {\it abelian variety} of dimension $g$ is a compact complex torus $\bC^{g}/{\bZ^{2g}}$, which is also a complex projective variety.
A  \index{semi-abelian variety} {\it semi-abelian variety} $G$ is an abelian complex algebraic group, which is an extension
of an abelian variety by a complex affine torus. 
(Semi-)abelian varieties are naturally associated to smooth (quasi-)projective varieties via the {\it Albanese map} construction, see \cite{Iit}.
Specifically, if $X$ is a smooth  complex (quasi-)projective variety, the \index{Albanese map} {\it Albanese map} of $X$ is a morphism $$\alb\colon  X \to \Alb(X)$$ from $X$ to a (semi-)abelian variety $\Alb(X)$ such that for any morphism $f\colon  X\to G$ to a semi-abelian variety $G$, there exists a unique morphism $g\colon \Alb(X) \to G$ such that the following diagram commutes: 
\begin{center}
$\xymatrix{
X \ar[rd]^{f}  \ar[r]^{\alb\quad} & \Alb(X) \ar[d]^{g}\\
           &   G 
}$
\end{center} 
Here, $\Alb(X)$ is called the \index{Albanese variety}  {\it Albanese variety} associated to $X$.

The Albanese map induces an isomorphism on the torsion-free part of the first integral homology groups, i.e., 
\be\label{eq0}
H_1(X;\bZ)/\text{Torsion} \overset{\sim}\longrightarrow H_1(\Alb(X);\bZ).
\ee
In particular, there is an identification:
\be\label{eq1}
\Char(X) \simeq \Char(\Alb(X)),
\ee
and the equality of Betti numbers $b_1(X)=b_1(\Alb(X))$.

By using the projection formula, for any $\rho \in \Char(X)\simeq \Char(\Alb(X))$ one gets: 
\be\label{eq2}
H^i(X; L_\rho) \simeq H^i(X; \bC_X \otimes L_\rho) 
\simeq H^i\big(\Alb(X); (R\alb_* \bC_X)\otimes L_\rho\big).
\ee
If, moreover, the Albanese map $\alb\colon X \to \Alb(X)$ is proper (e.g., if $X$ is projective), then the BBD decomposition theorem \cite{BBD} 
yields that $R \alb_* \bC_X $ is a direct sum of (shifted)  
{perverse sheaves}. In view of (\ref{eq2}), this provides a description of the cohomology jump loci of $X$ in terms of cohomology jump loci of certain perverse sheaves on the semi-abelian variety $\Alb(X)$. This motivates the following.

\bd Let $\cFb \in D^b_c(G;\bC)$ be a bounded constructible complex of $\bC$-sheaves on a semi-abelian variety $G$. 
The {\it degree $i$ cohomology jump locus of $\cFb$} is defined as:
$$\sV^i(G,\cFb): = \{\rho \in \Char(G)\mid H^i(G; \cFb \otimes_\bC L_\rho) \neq 0\}.$$
\ed

Back to the cohomology jump loci of $X$,  we note that (\ref{eq2}) yields the following identification:
\be\label{eq3}
{\sV^i(X)= \sV^i(\Alb(X), R\alb_* \bC_X)}.
\ee


\subsection{Mellin transformation and applications}
Let $G$ be a semi-abelian variety defined by an extension
$$1 \to T \to G \to A \to 1,$$
where $A$ is an abelian variety of dimension $g$ and $T\simeq(\bC^*)^m$ is a complex affine torus of dimension $m$. Set 
$$\Gamma_G:=\bC[\pi_1(G)],$$ 
and note that $\Gamma_G$ is a Laurent polynomial ring in $m+2g$ variables. Moreover, $$\Char(G)\simeq \spec \Gamma_G.$$
Let $\cL_G$ be the (universal) rank-one local system of $\Gamma_G$-modules on $G$, defined by mapping the generators of $\pi_1(G)$ to the multiplication by the corresponding variables of  $\Gamma_G$.

\bd\cite{GL} The \index{Mellin transformation} {\it Mellin transformation}  $\sM_{\ast}\colon  D^{b}_{c}(G; \bC) \to D^{b}_{coh}(\Gamma_G)$ is given by
$${\sM_\ast(\cFb) := Ra_\ast( \cL_G\otimes_{\bC}\cFb)}, $$
where $a \colon G\to point$ is the constant map to a point space, and $D^{b}_{coh}(\Gamma_G)$ denotes  the bounded coherent complexes of $\Gamma_G$-modules. 
\ed

The Mellin transformation can be used to completely characterize perverse sheaves on complex affine tori. More precisely, one has the following result due to Gabber-Loeser \cite[Theorem 3.4.1 and Theorem 3.4.7]{GL} in the $l$-adic context, and extended to the present form in  \cite[Theorem 3.2]{LMW3}: 
\bt\label{GL}
 A constructible complex $\cFb \in D^b_c(T;\bC)$ on a complex affine torus $T$  is perverse if, and only if, \begin{center} $H^i(\sM_\ast(\cFb))=0$ \ for all $i \neq 0$.\end{center} 
\et

In the context of abelian varieties, the Mellin transformation was used in \cite{BSS} for proving generic vanishing results for the cohomology of perverse sheaves. Inducting on the dimension of the complex affine torus $T$, the result of \cite{BSS} was extended to the semi-abelian context as follows.

\bt\cite[Theorem 4.3]{LMW4} \label{Mellin}  For any $\bC$-perverse sheaf $\cFb$  on a  {semi-abelian variety} $G$, one has: \begin{center} $H^i( \sM_* (\cFb))=0$ \ for $i<0$, \end{center} 
and \begin{center}
$H^i (D_{\Gamma_G} (\sM_*(\cFb))) =0$ \ for $i<0$.
\end{center}
Here $D_{\Gamma_G}(-):=\rhomo_{\Gamma_G} (-,\Gamma_G) $ is the dualizing functor for $\Gamma_G$-modules.
\et

The Mellin transformation can be used to translate the question of understanding the cohomology jump loci of a constructible complex to a problem in commutative algebra.
Specifically, by the projection formula, cohomology jump loci of $\cFb \in D^b_c(G;\bC)$  are determined by those of the Mellin transformation $\sM_\ast(\cFb)$ of $\cFb$ as follows (see \cite{LMW3}, \cite{LMW4}):
\be\label{eq5}{\sV^i(G,\cFb)=\sV^i(\sM_\ast(\cFb))}.\ee
Here, if $R$ is a Noetherian domain and $E^{\centerdot}$ is a bounded complex of $R$-modules with finitely generated cohomology, the closed points of $ \sV^{i} (E^\centerdot)$ can be described as
$$ \sV^{i} (E^\centerdot):= \{ \mathfrak{p} \in \spec R \mid H^i (F^\centerdot \otimes _R R/\mathfrak{p}) \neq 0 \},$$ 
 with $F^\centerdot$ a bounded above finitely generated {\it free} resolution of $E^\centerdot$. One also has the following:

\bp \label{exact} \cite[Lemma 2.8]{BSS} Let $R$ be a regular Noetherian domain and $E^{\centerdot}$ a bounded complex of $R$-modules with finitely generated cohomology. Then $H^i(E^{\centerdot} )=0$ for $i<0$ if, and only if, $\codim \sV^{ -i}(E^{\centerdot}) \geq i$ for $i\geq 0$.
\ep
 
By using the identification (\ref{eq5}), together with  Proposition \ref{exact} and  standard facts from commutative algebra, the following result of \cite{LMW4} is a direct consequence of Theorem \ref{Mellin}:
\bt\cite[Theorem 4.7]{LMW4}\label{pro}   \index{propagation}  For any $\bC$-perverse sheaf $\cFb$ on a semi-abelian variety $G$, the cohomology jump loci of $\cFb$ satisfy the following properties:
\begin{enumerate}
\item[(i)] {\it Propagation}: \index{propagation}
$$
{ \sV^{-m-g}(G,\cFb) \subseteq \cdots   \subseteq \sV^{0}(G,\cFb) \supseteq \sV^{1}(G,\cFb) \supseteq \cdots \supseteq \sV^{g}(G,\cFb)}.$$
Moreover, $\sV^i(G,\cFb)= \emptyset$ if $i\notin [-m-g, g]$.
\item[(ii)] {\it Codimension lower bound}: \index{codimension lower bound} for all $i\geq 0$, $$ {\codim \sV^{ i}(G,\cFb) \geq i} \ \
{\rm and} \  \ {\codim \sV^{ -i}(G,\cFb) \geq i}.$$
\end{enumerate}
\et

Theorem \ref{pro} is inspired by, and can be viewed as a topological counterpart of, similar properties satisfied by the Green-Lazarsfeld jump loci of topologically trivial line bundles, see \cite{GLa, GrLa}.

\br
Let $\cFb$ be a $\bC$-perverse sheaf so that not all  $H^j(G; \cFb)$ are zero.  
The propagation property (i) can be restated as saying that the set of integers $j$ for which $H^j(G;\cFb)\neq 0$ form an {\it interval} of consecutive integers. Indeed, let 
\begin{center} $k_+:=\max\{j\mid H^j(G; \cFb)\neq 0\} \text{ and } k_-:=\min\{j \mid H^j(G; \cFb)\neq 0\}.$\end{center} 
Then it's easy to see that the propagation property (i) is equivalent to: $k_+\geq 0$, $k_-\leq 0$ and
$$H^j(G; \cFb)\neq 0 \ \iff \ k_-\leq j\leq k_+. $$
In this form, the result of Theorem \ref{pro}(i) provides a generalization of a result of Weissauer \cite[Corollary 1]{We16b} from the abelian context. 
\er

A nice consequence of Theorem \ref{pro} is the following \index{generic vanishing} generic vanishing result: 
\bc\label{c1} For any $\bC$-perverse sheaf $\cFb$ on a semi-abelian variety $G$,  $$H^{i}(G; \cFb \otimes_{\bC} L_\rho)=0$$ for any generic rank-one $\bC$-local system $L_\rho$ and all $i\neq 0$. 
In particular, \index{signed Euler characteristic property}
 \be\label{eq6}  \chi(G,\cFb) 
 \geq 0.\ee
Moreover, $ \chi(G,\cFb)=0$ if, and only if, $\sV^0(G,\cFb) \neq \Char(G)$.
\ec

The above generic vanishing statement was originally proved by other methods in 
 \cite[Theorem 2.1]{Kra} in the $l$-adic context and further generalized to arbitrary field coefficients in \cite[Theorem 1.1]{LMW2}. For abelian varieties, generic vanishing results were obtained in \cite[Theorem 1.1]{KW}, \cite[Corollary 7.5]{Sch}, \cite[Vanishing Theorem]{We16} or \cite[Theorem 1.1]{BSS}. The signed Euler characteristic property (\ref{eq6}) is originally due to Franecki and Kapranov \cite[Corollary 1.4]{FK} (and compare with the
Example \ref{cc-abv} for another approach to this signed Euler characteristic property).
 
\br \index{propagation package} The properties of perverse sheaves from Theorem \ref{pro} and Corollary \ref{c1} are collectively termed the {\it propagation package} for perverse sheaves on semi-abelian varieties.\er


\subsection{Characterization of perverse sheaves on semi-abelian varieties}

In this section, we discuss a complete (global) characterization of perverse sheaves on semi-abelian varieties; see \cite{LMW4} for complete details.  
Motivation is also provided by the following result:
\bt\cite{Sch} \label{Sch}
If $A$ is an abelian variety and $\cFb \in D^b_c(A;\bC)$, then  
$\cFb$ is perverse if, and only if,  for all $i \in \bZ$, 
$\codim \sV^{i}(A,\cFb) \geq \vert 2i \vert.$
\et

Furthermore, as a consequence of Theorem \ref{GL}, Proposition \ref{exact} and Artin's vanishing theorem \ref{Artinvan}, one gets the following:
\bt \cite[Corollary 6.8]{LMW4}\label{LMWaf}
$\cFb\in D^b_c(T; \bC)$ is perverse on a complex affine torus $T$ if, and only if, 
\begin{enumerate}
\item[(i)] For all $i>0$: $\sV^i(T, \cFb)=\emptyset$, and 
\item[(ii)] For all $i\leq 0$: $\codim \sV^i(T, \cFb)\geq -i$. 
\end{enumerate}
\et

In order to unify and generalize the results of Theorems \ref{Sch} and \ref{LMWaf} to the semi-abelian context, one can make use of the new notions of (co)dimension for the cohomology jump loci, which were introduced in \cite{LMW4}. 
First recall the following.
\bd\label{lin} A closed irreducible subvariety $V$ of $\Char(G)$ is called {\it linear} if
there exists a short exact sequence of semi-abelian varieties 
\begin{center}
$1 \to G''(V)\to G\overset{q}{\to} G'(V) \to 1$
\end{center} and  some $\rho\in \Char(G)$ such that
\begin{center}
$V:=\rho \cdot \im (q^\#\colon  \Char(G'(V)) \to \Char(G))$.
\end{center} 
Here $G''(V)$ and $G'(V)$ depend on $V$, and $q^\#$ is induced by $q\colon G\to G'(V)$.
\ed
With the above definition, one has the following important structure result:
\bt\cite{BW17} \label{bw}
For any $\cFb \in D^b_c(G;\bC)$, each jump locus $\sV^i(G,\cFb)$ is a finite union of  linear subvarieties of $\Char(G)$.
\et

\bd Let $G$ be a semi-abelian variety and let $V$ be an irreducible linear subvariety of $\Char(G)$. In the notations of Definition \ref{lin}, let  
$T''(V)$ and $A''(V)$ denote the complex affine torus and, resp., the  abelian variety part of $G''(V)$. Define:
\begin{itemize}
\item[] {\it abelian codimension}: $\codim_a V:=\dim A''(V),$ \index{abelian codimension}
\item[] {\it semi-abelian codimension}: $\codim_{sa} V:=\dim G''(V).$ \index{semi-abelian codimension}
\end{itemize}
Similar notions can be defined for reducible subvarieties by taking the minimum among all irreducible components. 
\ed

\br Let $V$ be a nonempty linear subvariety of $\spec \Gamma_G$. 
\begin{enumerate}
\item If $G=T$ is a complex affine torus, then: $\codim_{sa} V =\codim V$ and $\codim_{a} V=0$. 
\item If $G=A$ is an abelian variety, we have: $\codim_{sa}(V)=\codim_{a}(V)=\frac{1}{2}\codim(V)$. 
\end{enumerate}
\er

In the above notations, the following generalization of Schnell's result was obtained in \cite{LMW4}:
\bt\cite[Theorem 6.6]{LMW4} \label{main}
A constructible complex $\cFb\in D^b_c(G; \bC)$ is perverse on $G$ if, and only if, 
\begin{enumerate}
\item[(i)] $\codim_{a} \sV^i(G, \cFb) \geq i$ \ for all $i\geq 0$, and
\item[(ii)] $\codim_{sa} \sV^i(G, \cFb) \geq -i$ \ for all $i\leq 0$. 
\end{enumerate}
\et

\begin{proof} 
The ``only if'' part is proved by induction on $\dim T$, using Theorem \ref{Sch} as the beginning step of the induction process. For the ``if'' part, one shows that the two codimension lower bounds in the statement of Theorem \ref{main} are sharp. See \cite[Section 6]{LMW4} for complete details.
\end{proof}


\subsection{Application: cohomology jump loci of quasi-projective manifolds}
The results of Theorem \ref{pro}, Theorem \ref{main} and Corollary \ref{c1} can be directly applied for the study of cohomology jump loci $\sV^i(X) \subseteq \Char(X)=\Char(\Alb(X))$ of a smooth complex quasi-projective \index{quasi-projective}
 variety $X$. Specifically, one has the following.
\bt\label{jlp}
Let  $X$ be a  smooth quasi-projective variety of complex dimension $n$.  Assume that 
$R\alb_*( \bC_X[n])$ is a perverse sheaf on $\Alb(X)$ (e.g., $\alb$ is proper and semi-small, or $\alb$ is quasi-finite).
Then the cohomology jump loci $\sV^i(X) $ satisfy the following properties:
\begin{itemize}
\item[(1)] \ \index{propagation}  {\it Propagation property}:
$$ 
  \sV^{n}(X) \supseteq \sV^{n-1}(X)\supseteq \cdots \supseteq \sV^{0}(X) = \{ \bC_X \};
$$
$$ 
  \sV^{n}(X) \supseteq \sV^{n+1}(X) \supseteq \cdots \supseteq \sV^{2n}(X)  .
$$
\item[(2)] \ \index{codimension lower bound} {\it Codimension lower bound}: for all $ i\geq 0$, 
$$ \codim_{sa} \sV^{n-i}(X) \geq i \ \ {\rm and } \ \ \codim_a \sV^{n+i}(X) \geq i.$$
\item[(3)] \ \index{generic vanishing}  {\it Generic vanishing:} $H^{i}(X,  L_\rho)=0$ 
for generic $\rho \in \Char(X)$ and all $ i \neq n$.
\item[(4)] \  \index{signed Euler characteristic property} {\it Signed Euler characteristic property}:  $(-1)^n  \chi(X) \geq 0$.
\item[(5)] \ \index{Betti property} {\it Betti property}: $b_i(X)>0$ for all $i\in [0, n]$, and $b_1(X)\geq n$.
\end{itemize} 
\et

 \bex Let $X$ be a smooth closed subvariety of a semi-abelian variety $G$. The closed embedding  $i\colon X\hookrightarrow G$ is a proper semi-small map, and hence the Albanese map $\alb\colon X \to \Alb(X)$ is also proper and semi-small. Then $R \alb_* (\bC_X[\dim X])$ is a perverse sheaf on $\Alb(X)$ and the jump loci of $X$ satisfy the properties listed in Theorem \ref{jlp}.  \eex 
 
 
\subsection{Application: topological characterization of semi-abelian varieties}
The Structure Theorem \ref{bw} together with the propagation package of Theorem \ref{pro} and Corollary \ref{c1} can be used to give the following topological characterization of semi-abelian varieties (see \cite[Proposition 7.7]{LMW4}):
\bt\label{c3.5} \label{homotopy} Let $X$ be a smooth quasi-projective variety with proper Albanese map (e.g., $X$ is projective), and assume that $X$ is homotopy equivalent to a real torus. Then $X$ is isomorphic to a semi-abelian variety. 
\et
\begin{proof} Assume $X$ has complex dimension $n$.
By the BBD decomposition theorem \cite{BBD}, $R\alb_* (\bC_X[n])$ is a direct sum of shifted semi-simple perverse sheaves on $\Alb(X)$. Denote by $\mathcal{U}$ the collection of all simple summands $\cFb$ appearing (up to a shift) in $R\alb_* (\bC_X[n])$. Then, by using Theorem \ref{pro}(i) and the identification (\ref{eq3}), one gets that \be\label{un} \bigcup_{i=0}^{2n} \sV^i(X) = \bigcup_{\mathcal{U}} \sV^0(\Alb(X),\cFb).\ee  
Since $X$ is homotopy equivalent to a real torus, a direct calculation yields that  $\bigcup_{i=0}^{2n} \sV^i(X)$ is just an isolated point. Hence, by (\ref{un}), for every simple perverse sheaf $\cFb \in \mathcal{U}$, the jump locus $\sV^0(\Alb(X),\cFb)$ is  this isolated point, so Corollary \ref{c1} gives $\chi(\Alb(X),\cFb)=0$. 

Simple perverse sheaves with zero Euler characteristic on semi-abelian varieties are completely described in \cite[Theorem 5.5]{LMW4} by using the structure Theorem \ref{bw} and the propagation package for their cohomology jump loci. In particular, it follows that for all $\cFb \in \mathcal{U}$ one has that  $\cFb=\bC_{\Alb(X)}[\dim \Alb(X)]$.  So  $R\alb_* \bC_X$ is a direct sum of shifted rank-one constant sheaves on $\Alb(X)$. Since $X$ and $\Alb(X)$ are both homotopy equivalent to tori, and since $b_1(X)=b_1(\Alb(X))$, one gets that $b_i(X)=b_i(\Alb(X))$ for all $i$. Therefore, $$R\alb_* \bC_X\simeq \bC_{\Alb(X)}.$$ Since $\alb$ is proper, it follows that all fibers of $\alb$ are zero-dimensional. Then it can be seen easily that $\alb$ is in fact an isomorphism.
\end{proof}

Under the assumptions of Theorem \ref{homotopy}, it follows that the Albanese map $\alb\colon  X\to \Alb(X)$ is an isomorphism. One can similarly prove the following special case of a question of Bobadilla-Koll\'ar \cite{BK}, see \cite{LMW6} for more general results:
\bp
Let $X$ be a complex projective manifold, and denote by $X^{ab}$ the \index{universal free abelian cover}
universal free abelian cover of $X$, i.e., the covering associated with the homomorphism $\pi_1(X) \to H_1(X;\bZ)/Torsion$. If $X^{ab}$ is homotopy equivalent to a finite CW complex, then the sheaves $R^i\alb_*\bC_X$ are local systems on $\Alb(X)$ for all $i \geq 0$.
\ep

\br
The study of perverse sheaves on semi-abelian varieties has other interesting topological applications, e.g., strong finiteness properties for Alexander-type invariants, generic vanishing of Novikov and $L^2$-homology, the study of homological duality properties of smooth complex algebraic varieties, etc., see \cite{EHMW}, \cite{LMW2}, \cite{LMW3}, \cite{LMW4} for more details.
\er




\end{document}